\documentclass[12pt]{article}
\usepackage{amsmath,amsfonts,amssymb,amsthm}
\newcommand{\Hom}{{\rm Hom}\,}
\newcommand{\End}{{\rm End}\,}
\newcommand{\Res}{{\rm Res}\,}

\newcommand{\Aut}{{\rm Aut}\,}

\newcommand{\D}{\mathcal{D}}
\newcommand{\E}{\mathcal{E}}
\newcommand{\Y}{\mathcal{Y}}
\newcommand{\Z}{\mathbb{Z}}

\newcommand{\C}{\mathbb{C}}
\newcommand{\N}{\mathbb{N}}
\newcommand{\g}{\mathfrak{g}}
\newcommand{\h}{\mathfrak{h}}

\def \<{\langle} 
\def \>{\rangle}
\def \be{\begin{equation}\label}
\def \ee{\end{equation}}
\def \bex{\begin{exa}\label}
\def \eex{\end{exa}}
\def \bl{\begin{lem}\label}
\def \el{\end{lem}}
\def \bt{\begin{thm}\label}
\def \et{\end{thm}}
\def \bp{\begin{prop}\label}
\def \ep{\end{prop}}
\def \br{\begin{rem}\label}
\def \er{\end{rem}}
\def \bc{\begin{coro}\label}
\def \ec{\end{coro}}
\def \bd{\begin{de}\label}
\def \ed{\end{de}}

\newtheorem{thm}{Theorem}[section]
\newtheorem{prop}[thm]{Proposition}
\newtheorem{coro}[thm]{Corollary}

\newtheorem{exa}[thm]{Example}
\newtheorem{lem}[thm]{Lemma}
\newtheorem{rem}[thm]{Remark}
\newtheorem{de}[thm]{Definition}

\theoremstyle{definition}
\theoremstyle{remark}

\numberwithin{equation}{section}

\makeatletter
\@addtoreset{equation}{section}

\makeatother

\begin{document}
\begin{flushright}
\today
\\
Version $0.91$
\end{flushright}

\begin{Large}
\begin{center}
\textbf{A new construction of vertex algebras and quasi modules for vertex algebras}
\end{center}
\end{Large}

\begin{center}
{Haisheng Li\footnote{Partially supported by a NSA grant}\\
Department of Mathematical Sciences, Rutgers University, Camden, NJ 08102\\
and\\
Department of Mathematics, Harbin Normal University, Harbin, China}
\end{center}

\begin{abstract}
In this paper, a new construction of vertex algebras from
more general vertex operators is given and
a notion of quasi module for vertex algebras is introduced and studied.
More specifically, a notion of quasi local subset(space) of $\Hom
(W,W((x)))$ for any vector space $W$ is introduced and studied,
generalizing the notion of usual locality in the most possible way, and
it is proved that on any maximal quasi local subspace there exists a
natural vertex algebra structure and that any quasi local subset of
$\Hom (W,W((x)))$ generates a vertex algebra.  Furthermore, a notion of quasi
module for a vertex algebra is introduced and it is proved that $W$ is
a quasi module for each of the vertex algebras generated by quasi
local subsets of $\Hom (W,W((x)))$.  A notion of $\Gamma$-vertex
algebra is also introduced and studied, where $\Gamma$ is a subgroup
of the multiplicative group $\C^{\times}$ of nonzero complex numbers.
It is proved that any maximal quasi local subspace of $\Hom
(W,W((x)))$ is naturally a $\Gamma$-vertex algebra and that any quasi
local subset of $\Hom (W,W((x)))$ generates a $\Gamma$-vertex algebra.
It is also proved that a $\Gamma$-vertex algebra exactly amounts to a
vertex algebra equipped with a $\Gamma$-module structure which
satisfies a certain compatibility condition.  Finally, three families
of examples are given, involving twisted affine Lie algebras, certain
quantum Heisenberg algebras and certain quantum torus Lie algebras.
\end{abstract}

\section{Introduction}

Vertex (operator) algebras are often viewed as some kind of generalized
algebras, equipped with infinitely many compatible multiplications
parameterized by $\Z$.  Another viewpoint is that vertex operator
algebras are ``algebras'' of vertex operators, just as classical
(associative or Lie) algebras are algebras of linear operators.  
This particular viewpoint
was emphasized in [FLM] and the philosophy was deeply reflected in 
\cite{li-local} where a theory of representation for a vertex
algebra was developed and a general construction of vertex algebras
together with modules was given, by using vertex operators.

For a vector space $W$, a (weak) vertex operator on $W$
(\cite{li-local}, [LL]) is an element of $\Hom (W,W((x)))$ which is
alternatively denoted by $\E(W)$.  A subset (subspace) $U$ of $\E(W)$
is said to be {\em local} if for any $a(x),b(x)\in U$, there exists a
nonnegative integer $k$ such that
\begin{eqnarray}\label{e1.1}
(x_{1}-x_{2})^{k}a(x_{1})b(x_{2})=(x_{1}-x_{2})^{k}b(x_{2})a(x_{1}).
\end{eqnarray}
It was proved in \cite{li-local} that any maximal local subspace of
$\E(W)$ is naturally a vertex algebra with $W$ as a natural module,
where the identity operator $1_{W}$ is the vacuum vector and for
$a(x),b(x)\in \E(W),\; n\in \Z$,
\begin{eqnarray}\label{e1.4}
a(x)_{n}b(x)
=\Res_{x_{1}}\left((x_{1}-x)^{n}a(x_{1})b(x)
-(-x+x_{1})^{n}b(x)a(x_{1})\right).
\end{eqnarray}
Furthermore, it was proved therein that any local subset of $\E(W)$   
canonically generates a
vertex algebra with $W$ a module.  (Notice that Zone Lemma
simply implies that any local subset of $\E(W)$ is contained in some
maximal local subspace.)  If $W$ is taken to be a highest weight
module (or more generally a restricted module) for an {\em untwisted}
affine Kac-Moody Lie algebra or the Virasoro algebra, the
generating functions form a local subspace of $\E(W)$, so that one has
a natural vertex algebra with $W$ as a module.  In this way one has
vertex (operator) algebras and modules associated with affine Lie
algebras (including infinite-dimensional Heisenberg algebras) and
the Virasoro algebra.  (Historically, vertex operator algebras were
constructed differently (see [B1], [DL], [FLM], [FZ]).)

In the theory of Lie algebras, similar to untwisted affine Kac-Moody
Lie algebras, there are some other interesting
Lie algebras such as {\em twisted} affine Kac-Moody Lie algebras (see
[K1]), Lie algebras of $q$-pseudodifferential operators 
on the circle (sin-algebra) (see [G-K-L])
and some quantum torus Lie algebras (see [G-K-L], [BGT]).  
If $W$ is a ``restricted'' module for such a Lie algebra, the generating
functions, which are elements of $\E(W)$, 
are not local. Instead, they satisfy the relation
\begin{eqnarray}\label{e1.3}
p(x_{1}/x_{2})a(x_{1})b(x_{2})=p(x_{1}/x_{2})b(x_{2})a(x_{1})
\end{eqnarray}
(cf. (\ref{e1.1})) for some nonzero polynomial $p(x)$.  It has been a
fundamental question whether these Lie algebras (through restricted
modules) can be associated with some vertex-algebra-like
structures in the same (or similar) way as (or to) that 
untwisted affine Lie algebras are associated with vertex algebras.  
In [G-K-K], Golenishcheva-Kutuzova and Kac associated
these Lie algebras (without central extension)
to certain Lie-algebra-like structures called $\Gamma$-conformal
algebras. (In fact, as it was proved therein, 
$\Gamma$-conformal algebras are Lie algebras acted by $\Gamma$
by automorphisms.)  But the original
question  still remains open.

In this paper, we give a complete answer to this question by establishing
a conceptual result analogous to that of [Li2].
We found that these Lie algebras can also be associated with ordinary
vertex algebras in a certain natural way. Namely we proved that the
generating functions as elements of $\E(W)$ for a restricted module
$W$ still generate an ordinary vertex algebra in a certain {\em new}
and natural way.  This is rather unexpected. First, due to the general
locality (\ref{e1.3}) of vertex operators we started with, we would expect certain
vertex-algebra-like objects which are essentially defined by the
general locality.  (Motivated by our result we then prove (Proposition
\ref{pno-other}) that vertex-algebra-like objects defined by the
general locality are ordinary vertex algebras.)  Second, the known
constructions of vertex algebras from local (twisted) vertex operators
(see \cite{li-local}, \cite{li-twisted}, cf. \cite{li-G1}, [GL])
correspond to the notions of (twisted) module or (twisted)
representations for vertex algebras, but there is no notion other than
that of (twisted) module in the aspect of representation.  In reality,
those generating functions generate an ordinary vertex algebra, but
for this vertex algebra, the space $W$ equipped with the natural
action is neither a module nor a twisted module in the usual sense.
Then this new construction of vertex algebras from more general vertex
operators naturally leads us to a notion of what we call quasi module
for a vertex algebra.  Quasi modules naturally include modules and
certain ``deformation'' of twisted modules as we shall see through
examples.  Just as the notions of module and twisted module, this new
notion should also be of fundamental importance.

In the following we give an outline of this paper. 
We first start with what we call quasi local subspaces of $\E(W)$ for an
arbitrarily given vector space $W$.
A subset (subspace) $U$ of $\E(W)$ is {\em quasi local} if
for any $a(x),b(x)\in U$, there exists a nonzero polynomial $f(x_{1},x_{2})$
such that 
\begin{eqnarray}
f(x_{1},x_{2})a(x_{1})b(x_{2})=f(x_{1},x_{2})b(x_{2})a(x_{1}).
\end{eqnarray}
(Arguably, this is the most possible generalization of usual locality (\ref{e1.1}).)

We then consider appropriate actions of vertex operators on vertex operators.
It is important to point out that in this generality, the
actions of vertex operators on vertex operators defined by
(\ref{e1.4}) are not appropriate.  In fact, as we have seen in
\cite{li-local}, even in the (fermionic) supercase, (\ref{e1.4}) should
be modified; if
\begin{eqnarray}
(x_{1}-x_{2})^{k}a(x_{1})b(x_{2})=- (x_{1}-x_{2})^{k}b(x_{2})a(x_{1})
\end{eqnarray}
for some nonnegative integer $k$, $a(x)_{n}b(x)$ should be defined by
\begin{eqnarray}
a(x)_{n}b(x)
=\Res_{x_{1}}\left((x_{1}-x)^{n}a(x_{1})b(x)+ (-x+x_{1})^{n}b(x)a(x_{1})\right).
\end{eqnarray}
On the other hand, it was proved in [LL] (see also [R]) that
for any subspace $U$ of $\E(W)$ which contains the identity operator $1_{W}$ 
and is closed under the action defined by
(\ref{e1.4}), $U$ is a vertex algebra if and only if $U$ is local.

In \cite{li-G1}, motivated by [B2] we have studied suitably defined
``compatible'' subspaces of $\E(W)$ where compatibility generalizes
the usual locality in a certain direction.  We introduced a notion of
(axiomatic) $G_{1}$-vertex algebra and showed that compatible
subspaces of $\E(W)$ naturally give rise to $G_{1}$-vertex algebras
with $W$ as a module. In \cite{li-G1}, the actions of vertex operators
on vertex operations were essentially defined by using
the operator product expansion or weak associativity,
just as the product for classical
associative algebras are defined by the associativity $(ab)w=a(bw)$.
Informally, $a(x)_{n}b(x)$ for $n\in \Z$ were defined 
in terms of generating function
$\Y(a(x),x_{0})b(x)=\sum_{n\in \Z}a(x)_{n}b(x) x_{0}^{-n-1}$ by
\begin{eqnarray}
\<w^{*},(\Y (a(x),x_{0})b(x))w\>= \<w^{*}, a(x_{0}+x)b(x)w\>
\end{eqnarray}
for $w^{*}\in W^{*},\; w\in W$.
(This was precisely defined in terms of $\iota$-maps.)
In this paper, we extend this definition by 
naturally extending the $\iota$-maps explored
in [FLM], [FHL] and [LL].

Note that this definition requires that the product $a(x_{1})b(x_{2})$
be of a certain form, so that $\Y$ is only a {\em partial} map on
$\E(W)\otimes \E(W)$, unlike the usual case where $\Y$ was defined on
the whole space $\E(W)\otimes \E(W)$.  Nevertheless, for any quasi
local subspace $U$ of $\E(W)$, $\Y$ is a linear map from $U\otimes U$
to $\E(W)((x))$.  If $U$ is a maximal quasi local subspace of $\E(W)$,
we show that $U$ contains $1_{W}$, $\Y$ maps $U\otimes U$ into
$U((x))$ and that $(U,\Y,1_{W})$ carries the structure of a vertex
algebra.  For any quasi local subset $S$ of $\E(W)$, in view of Zorn
lemma, there exists a maximal quasi local subspace $U$ containing $S$.
Consequently, $S$ generates a (canonical) vertex algebra inside $\E(W)$.

There is a new feature of quasi locality. Similar to the map
$\Y$ we define a family of (partial) maps $\Y_{\alpha}$ 
for $\alpha\in \C^{\times}$ by
\begin{eqnarray}
\<w^{*},(\Y_{\alpha} (a(x),x_{0})b(x))w\>= \<w^{*}, a(x_{0}+\alpha x)b(x)w\>
\end{eqnarray}
for $w^{*}\in W^{*},\; w\in W$. It is proved that if $U$ is a maximal
quasi local subspace of $\E(W)$, for any nonzero complex number
$\alpha$, $\Y_{\alpha}$ maps $U\otimes U$ to $U((x))$, so that $U$ is
equipped with a family of linear maps $\Y_{\alpha}$ for $\alpha\in
\C^{\times}$. (If $U$ is a maximal
{\em local} subspace of $\E(W)$, this is not true in general.)
We furthermore obtain a Jacobi-like identity for
$(U,\{\Y_{\alpha}\},1_{W})$ (Theorem \ref{tkey-main}).  Motivated by
this and by the notion of $\Gamma$-conformal algebra in [G-K-K] we
introduce a notion of $\Gamma$-vertex algebra, where $\Gamma$ is any
group equipped with a group homomorphism from $\Gamma$ to
$\C^{\times}$ and we show that a $\Gamma$-vertex algebra amounts to an
(ordinary) vertex algebra equipped with a $\Gamma$-module structure
satisfying a certain compatibility condition.  In terms of this
notion, we have proved (Theorem \ref{tkey-main}) that for any maximal
quasi local subspace $U$ of $\E(W)$,
 $(U,\{\Y_{\alpha}\},1_{W})$ carries the structure of 
a $\Gamma$-vertex algebra with
$\Gamma$ any subgroup of $\C^{\times}$ and that any quasi local subset
of $\E(W)$ generates a $\Gamma$-vertex algebra.

Another important issue is about the module structure on $W$
for vertex algebras in $\E(W)$, generated by quasi local subsets.
Even though we still get 
ordinary vertex algebras from quasi local vertex operators on $W$,
$W$ equipped with the natural action of those vertex algebras
is {\em not} a module in the usual sense, unlike the usual case
\cite{li-local}, in which $W$ is naturally a module.  
Motivated by this, we formulate and study a notion of what we
call ``quasi module'' for a vertex algebra $V$.
A quasi $V$-module satisfies all the axioms
in the definition of the notion of module 
except that the Jacobi identity axiom is
replaced by the axiom that for $u,v\in V$,
there exists a nonzero polynomial $f(x_{1},x_{2})$ 
(depending on $(u,v)$) such that the Jacobi identity
for $(u,v)$ multiplied by $f(x_{1},x_{2})$ holds.
In terms of this notion, we have that
for any vector space $W$, any maximal quasi local subspace of $\E(W)$ 
is naturally a vertex algebra with $W$ a quasi module and 
any quasi local subset of $\E(W)$ generates a vertex algebra 
with $W$ as a quasi module.

In the last section, as an application of our general results we study
three families of examples.  In the first family, we consider
untwisted affine Lie algebra $\hat{\g}$ associated to any Lie algebra
$\g$ equipped with a nondegenerate symmetric invariant bilinear form
and twisted affine Lie algebra $\hat{\g}[\sigma]$ associated with an
automorphism $\sigma$ of $\g$ of (finite) order.  We show that on any
restricted $\hat{\g}[\sigma]$-module $W$ of level $\ell\in \C$, there
exists a unique quasi module structure for the vertex operator algebra
$V_{\hat{\g}}(\ell,0)$ (cf. [LL]). It has been known
(\cite{li-twisted}, cf. [FLM]) that on any restricted
$\hat{\g}[\sigma]$-module $W$ of level $\ell$, there exists a unique
$\sigma$-twisted module structure for $V_{\hat{\g}}(\ell,0)$.  This
suggests that for a vertex algebra $V$, twisted $V$-modules with
respect to finite order automorphisms should be connected with quasi
modules.  In a sequel \cite{li-qm-tm}, we shall study the relations
between $\sigma$-twisted $V$-modules and 
quasi $V$-modules, for a general vertex operator algebra $V$ and
for any finite order automorphism $\sigma$ of $V$.  We conjecture that
the category of $\sigma$-twisted $V$-modules is canonically isomorphic
to a subcategory of quasi $V$-modules.  In this way, quasi modules
unify modules and twisted modules. This connection may shed some light
on some difficult problems in the study of orbifold theory (see [DM]).

In the second family, we study a certain ``quantum'' Heisenberg Lie
algebra $\hat{\h}_{q}$ associated with
a vector space $\h$ equipped with a nondegenerate symmetric bilinear
form and with a nonzero complex number $q$.
We show that any restricted $\hat{\h}_{q}$-module 
is a quasi module for a certain Heisenberg vertex algebra.
In the third family, we study certain quantum torus Lie algebras 
(see [G-K-L], [BGT]) and we show that
any restricted module is a quasi module for a certain vertex algebra
associated with a certain affine Lie algebra $\hat{\g}$.
It was known  (see [BGT]) that some of quantum torus Lie algebras 
are essentially extended affine Lie algebras (without derivations added).
Previously, toroidal extended affine Lie algebras
have been related to vertex algebras (see [BBS], [BDT]).
It is our belief that all the extended affine Lie algebras
(see [AABGP]) can be linked to (general) untwisted affine Lie algebras
and then to vertex algebras in a similar way.

While one of our motivations for this paper is to find a solution to
the question mentioned earlier, our main motivation is to develop a
theory of quantum vertex algebras so that quantum affine algebras in
[FJ] will naturally give rise to quantum vertex algebras just as
(untwisted) affine Lie algebras naturally give rise to vertex
algebras.  The quantum vertex operators studied in [FJ] do
not satisfy quasi locality, but they are ``compatible'' in the sense
that is defined here, so that the actions of vertex operators on vertex
operators have been defined. In several
places our setting is more  general than we need for this
paper. This is exactly for our study on quantum vertex algebras in a
sequel \cite{li-qva}. 

The main results of this paper were reported in 
the workshop ``Conformal Field Theory and Vertex Algebras,''
Osaka, Japan (January 10-12, 2004). 
We would like to thank Professors Nagatomo and Tsuchiya
for invitation and for their hospitality.

This paper is organized as follows: In Section 2, we extend the usual
iota maps and we discuss certain associativity and cancelation
properties. In Section 3, we introduce certain basic operations on
$\E(W)$ and present some basic properties.  In Section 4, we study
quasi local subspaces of $\E(W)$ and prove the key results.  
In Section 5, we study vertex
algebras and quasi modules.  In Section 6, we introduce and study the
notions of $\Gamma$-vertex algebra and quasi module.  In Section 7, we
give three families of examples of vertex algebras and quasi
modules.

\newpage

\section{Some formal calculus}
In this section we define certain iota maps, which generalize
the usual iota maps introduced in [FLM], [FHL] and [LL]. We also discuss
the subtle issues of associativity and cancelation laws in formal calculus.

First, throughout this paper,
$x,y, x_{0},x_{1},x_{2}, x_{3},\dots$ are mutually commuting independent
formal variables. Vector spaces are considered to be over the field $\C$ of
complex numbers, though any algebraically closed field of
characteristic zero will work fine.

For a vector space $U$, $U((x))$
is the space of lower truncated (infinite) integral
power series of $x$ with coefficients in $U$, 
$U[[x]]$ is the space of (infinite) nonnegative integral
power series in $x$ with coefficients in $U$ and
$U[[x,x^{-1}]]$ is the space of doubly infinite
integral power series in $x$ with coefficients in $U$.
Multi-variable analogues of these spaces are defined in the obvious way.
For example, $U[[x_{1},x_{2},\dots,x_{r}]]$ is the space of 
(infinite) nonnegative integral power series in $x_{1},\dots,x_{r}$ with coefficients in $U$
and $U((x_{1},x_{2},\dots,x_{r}))$ 
is the space of lower truncated (infinite) integral
power series in $x_{1},\dots,x_{r}$ with coefficients in $U$:
\begin{eqnarray}
U((x_{1},\dots,x_{r}))=U[[x_{1},\dots,x_{r}]][x_{1}^{-1},\dots,x_{r}^{-1}].
\end{eqnarray}

\br{rdef-iotamaps}
{\em Let $F$ be any field. Then $F((x))$ is a field. In view of this,
$\C((x_{1}))$, $\C((x_{1}))((x_{2}))$,\dots,
$\C((x_{1}))((x_{2}))\cdots ((x_{r}))$ are fields.
Clearly,
\begin{eqnarray}
\C[x_{1},\dots,x_{r}]\subset \C[[x_{1},\dots,x_{r}]]
\subset \C((x_{1},\dots,x_{r}))
\subset \C((x_{1}))((x_{2}))\cdots ((x_{r})).
\end{eqnarray}
Denote by $\C_{*}(x_{1},\dots,x_{r})$ the localization of 
$\C[[x_{1},\dots,x_{r}]]$ at
$\C[x_{1},\dots,x_{r}]^{\times}$ (the nonzero polynomials):
\begin{eqnarray}
\C_{*}(x_{1},\dots,x_{r})
=\C[[x_{1},\dots,x_{r}]]_{\C[x_{1},\dots,x_{r}]^{\times}}.
\end{eqnarray}
The field $\C(x_{1},\dots,x_{r})$ of rational functions, 
i.e., the fraction field of the polynomial ring $\C[x_{1},\dots,x_{r}]$,
is a subring of $\C_{*}(x_{1},\dots,x_{r})$. 
Since $\C[[x_{1},\dots,x_{r}]]$ is a subring of 
$\C((x_{1}))((x_{2}))\cdots ((x_{r}))$,
$\C_{*}(x_{1},\dots,x_{r})$ can be naturally embedded into 
$\C((x_{1}))((x_{2}))\cdots ((x_{r}))$.
Denote this embedding by $\iota_{x_{1},\dots,x_{r}}$: 
\begin{eqnarray}
\iota_{x_{1},\dots,x_{r}}:
\C_{*}(x_{1},\dots,x_{r})\rightarrow \C((x_{1}))((x_{2}))\cdots
((x_{r})).
\end{eqnarray}
For any permutation $\sigma$ on $\{1,\dots,r\}$,
$\iota_{x_{\sigma(1)},\dots,x_{\sigma(r)}}$ denotes
the embedding
\begin{eqnarray}
\iota_{x_{\sigma(1)},\dots,x_{\sigma(r)}}:
\C_{*}(x_{1},\dots,x_{r})\rightarrow \C((x_{\sigma(1)}))\cdots
((x_{\sigma(r)})).
\end{eqnarray}
In particular, $\iota_{x_{2},x_{1}}$ denotes
the embedding of $\C_{*}(x_{1},x_{2})$ into $\C((x_{2}))((x_{1}))$.
For example, for any nonzero complex number $\alpha$ and for any integer $n$,
\begin{eqnarray}
& &\iota_{x_{1},x_{2}}(x_{1}-\alpha x_{2})^{n}
=\sum_{i\ge 0}\binom{n}{i}(-\alpha)^{i}x_{1}^{n-i}x_{2}^{i}=(x_{1}-\alpha x_{2})^{n},\\
& &\iota_{x_{2},x_{1}}(x_{1}-\alpha x_{2})^{n}
=\sum_{i\ge 0}\binom{n}{i}(-\alpha)^{n-i}x_{2}^{n-i}x_{1}^{i}=(-\alpha x_{2}+x_{1})^{n},
\end{eqnarray}
where we are using the usual binomial series expansion conventions.
Notice that partial derivative operators $\frac{\partial}{\partial x_{i}}$ act on
both the domains and codomains and that partial derivative operators 
commute with the iota-maps. Furthermore, all the iota-maps
are $\C((x_{1},\dots,x_{r}))$-homomorphisms and they are identity
on $\C((x_{1},\dots,x_{r}))$.}
\er

In formal calculus ([FLM], [FHL], [LL]), associativity and 
cancelation for products of formal series are subtle issues. 
Let $A$ be an associative algebra with identity 
(over $\C$) and let $U$ be an $A$-module. Let
$$F,G\in A[[x_{1}^{\pm 1},\dots,x_{n}^{\pm 1}]],\;\;
H\in U[[x_{1}^{\pm 1},\dots,x_{n}^{\pm 1}]].$$
The associativity $F(GH)=(FG)H$ does not hold in general 
and on the other hand, it does hold
if all the products $GH$, $FG$ and $FGH$ exist. 
As in associative algebra theory,
associativity sometimes leads to cancelation.
In the following we give a few useful cases.

\bl{lcancellation}
Let $U$ be any vector space. 
(a) Let $f(x_{1},x_{2})\in U[[x_{1}^{\pm 1},x_{2}^{\pm 1}]],\; m,n\in
\Z$. Then
$x_{1}^{m}x_{2}^{n}f(x_{1},x_{2})=0$ if and only if $f(x_{1},x_{2})=0$.

(b) For $f(x_{1},x_{2}),g(x_{1},x_{2})\in \C((x_{1}))((x_{2}))$ 
and for $H(x_{1},x_{2})\in U((x_{1}))((x_{2}))$, we have
\begin{eqnarray}
f(x_{1},x_{2})(g(x_{1},x_{2})H(x_{1},x_{2}))
=(f(x_{1},x_{2})g(x_{1},x_{2}))H(x_{1},x_{2}).
\end{eqnarray}

(c) Let $H(x_{1},x_{2}),K(x_{1},x_{2})\in
U((x_{1}))((x_{2}))$. If 
\begin{eqnarray}
f(x_{1},x_{2})H(x_{1},x_{2})=f(x_{1},x_{2})K(x_{1},x_{2})
\end{eqnarray}
for some $0\ne f(x_{1},x_{2})\in \C((x_{1}))((x_{2}))$,
then $H(x_{1},x_{2})=K(x_{1},x_{2})$.

(d) Let
$H_{1}\in U((x))((x_{1}))((x_{2}))$ and $H_{2}\in U((x))((x_{2}))((x_{1}))$.
If there exists a (nonzero) polynomial $f(x_{1},x_{2})$ with
$f(x_{1},x_{1})\ne 0$ such that
$$f(x_{1}+x,x_{2}+x)H_{1}=f(x_{1}+x,x_{2}+x)H_{2},$$
then $H_{1}=H_{2}$.
\el

\begin{proof} 
For (a), notice that $U[[x_{1}^{\pm 1},x_{2}^{\pm 1}]]$ is 
a $\C[x_{1}^{\pm 1},x_{2}^{\pm 1}]$-module
and  $x_{1}^{m}x_{2}^{n}$ for any $m,n\in \Z$ is invertible in 
$\C[x_{1}^{\pm 1},x_{2}^{\pm 1}]$.
Then for $f(x_{1},x_{2})\in U[[x_{1}^{\pm 1},x_{2}^{\pm 1}]],\; m,n\in \Z$,
$x_{1}^{m}x_{2}^{n}f(x_{1},x_{2})=0$ if and only if $f(x_{1},x_{2})=0$.

Both (b) and (c) follow from the fact
that $\C((x_{1}))((x_{2}))$ is a field and 
$U((x_{1}))((x_{2}))$ is a vector space over 
$\C((x_{1}))((x_{2}))$. 

For (d), write
$$f(x_{1}+x,x_{2}+x)=f(x,x)+R(x,x_{1},x_{2}),$$
where $R(x,x_{1},x_{2})\in \C[x,x_{1},x_{2}]$ with
$R(x,0,0)=0$.
For each nonnegative integer $i$, we view $f(x,x)^{-i-1}$ as a 
formal series (its formal Laurent
series expansion at zero) in $\C((x))$.
Set
\begin{eqnarray}
g(x,x_{1},x_{2})=\sum_{i\ge
0}(-1)^{i}f(x,x)^{-1-i}R(x,x_{1},x_{2})^{i}
\in \C((x))[[x_{1},x_{2}]].
\end{eqnarray}
Clearly, $f(x_{1}+x,x_{2}+x)g(x,x_{1},x_{2})=1$. Then it follows from
the associativity law.
\end{proof}

\br{r-delta-function}
{\em We recall some fundamental
delta-function properties from [FLM], [FHL] and [LL].
For any nonzero complex number $\alpha$, we have
\begin{eqnarray}\label{edelta-identity}
x_{0}^{-1}\delta\left(\frac{x_{1}-\alpha x_{2}}{x_{0}}\right)
-x_{0}^{-1}\delta\left(\frac{-\alpha x_{2} +x_{1}}{x_{0}}\right)
=x_{1}^{-1}\delta\left(\frac{\alpha x_{2} +x_{0}}{x_{1}}\right).
\end{eqnarray}
For any $g(x_{0},x_{1},x_{2})\in \C((x_{0},x_{1},x_{2}))$, we have
\begin{eqnarray}
& &x_{0}^{-1}\delta\left(\frac{x_{1}-\alpha x_{2}}{x_{0}}\right)g(x_{0},x_{1},x_{2})
=x_{0}^{-1}\delta\left(\frac{x_{1}-\alpha x_{2}}{x_{0}}\right)g(x_{0},x_{0}+\alpha
x_{2},x_{2}),\label{e-delta-1}\\
& &x_{1}^{-1}\delta\left(\frac{\alpha x_{2}+x_{0}}{x_{1}}\right)g(x_{0},x_{1},x_{2})
=x_{1}^{-1}\delta\left(\frac{\alpha x_{2} +x_{0}}{x_{1}}\right)
g(x_{0},\alpha x_{2}+x_{0},x_{2}).\label{e-delta-3}
\end{eqnarray}
Furthermore, if $g(x_{0},x_{1},x_{2})$ involves only nonnegative powers of $x_{1}$,
using these and (\ref{edelta-identity}) we get
\begin{eqnarray}\label{e-delta-2}
x_{0}^{-1}\delta\left(\frac{-\alpha x_{2} +x_{1}}{x_{0}}\right)g(x_{0},x_{1},x_{2})
=x_{0}^{-1}\delta\left(\frac{-\alpha x_{2}+x_{1}}{x_{0}}\right)g(x_{0},\alpha
x_{2}+x_{0},x_{2}),
\end{eqnarray}
noticing that $g(x_{0},\alpha x_{2}+x_{0},x_{2})=g(x_{0},x_{0}+\alpha x_{2},x_{2})$.}
\er

The following is a reformulation of Proposition 3.4.2 of
[LL] with a slightly different proof (cf. \cite{li-local}, \cite{li-G1}):

\bl{lformaljacobiidentity}
Let $U$ be a vector space, let $\alpha$ be a nonzero complex number and let
\begin{eqnarray}
A(x_{1},x_{2})&\in& U((x_{1}))((x_{2})),\\
B(x_{1},x_{2})&\in& U((x_{2}))((x_{1})),\\
C(x_{0},x_{2})&\in& U((x_{2}))((x_{0})).
\end{eqnarray}
Then 
\begin{eqnarray}\label{eformaljacobiABC}
& &x_{0}^{-1}\delta\left(\frac{x_{1}-\alpha x_{2}}{x_{0}}\right)A(x_{1},x_{2})
-x_{0}^{-1}\delta\left(\frac{\alpha x_{2}-x_{1}}{-x_{0}}\right)B(x_{1},x_{2})
\nonumber\\
&=&x_{1}^{-1}\delta\left(\frac{\alpha x_{2}+x_{0}}{x_{1}}\right)C(x_{0},x_{2})
\end{eqnarray}
if and only if there exist nonnegative integers 
$k$ and $l$ such that
\begin{eqnarray}
(x_{1}-\alpha x_{2})^{k}A(x_{1},x_{2})
&=&(x_{1}-\alpha x_{2})^{k}B(x_{1},x_{2}),
\label{eformalA=B}\\
(x_{0}+\alpha x_{2})^{l}A(x_{0}+\alpha x_{2},x_{2})
&=&(x_{0}+\alpha x_{2})^{l}C(x_{0},x_{2}).
\label{eformalA=C}
\end{eqnarray}
\el

\begin{proof} In the case $\alpha=1$ it was proved in \cite{li-G1}.
The general case follows from the special case after
the substitution $x_{2}\rightarrow \alpha^{-1}x_{2}$.
\end{proof}

We shall also use the following result (cf. [DLM], [G-K-K], [K2], \cite{li-local}, [LL]):

\bl{lpa=0}
Let $U$ be a vector space, let $A(x_{1},x_{2})\in U[[x_{1}^{\pm 1},x_{2}^{\pm 1}]]$
and let $p(x)$ be a nonzero polynomial. Then
\begin{eqnarray}\label{epa=0}
p(x_{1}/x_{2})A(x_{1},x_{2})=0
\end{eqnarray}
if and only if
\begin{eqnarray}\label{ea=sum}
A(x_{1},x_{2})= \sum_{i=1}^{r}\sum_{j=0}^{k_{i}-1}
\left(\frac{\alpha_{i}^{-j}}{j!}\left(\frac{\partial}{\partial x_{2}}\right)^{j}
x_{1}^{-1}\delta\left(\frac{\alpha_{i}x_{2}}{x_{1}}\right)\right)B_{i,j}(x_{2})
\end{eqnarray}
for some $B_{i,j}(x)\in U[[x,x^{-1}]]$, where $\alpha_{1},\dots,\alpha_{r}$ 
are the nonzero (distinct) roots of $p(x)$ with multiplicities $k_{1},\dots, k_{r}$.
Furthermore, these $B_{i,j}(x)$ in (\ref{ea=sum}) are uniquely determined by $A(x_{1},x_{2})$.
\el

\begin{proof} From \cite{li-local} (cf. [LL]) we have
$$(x_{1}-x_{2})^{m}\left(\frac{\partial}{\partial x_{2}}\right)^{n}
x_{1}^{-1}\delta\left(\frac{x_{2}}{x_{1}}\right)=0\;\;\;\mbox{ for }m>n\ge 0.$$
For any nonzero complex number $\beta$, by changing variable $x_{2}$ we have
\begin{eqnarray}
(x_{1}-\beta x_{2})^{m}\left(\frac{\partial}{\partial x_{2}}\right)^{n}
x_{1}^{-1}\delta\left(\frac{\beta x_{2}}{x_{1}}\right)=0\;\;\;\mbox{ for }m>n\ge 0.
\end{eqnarray}
Using this we see immediately that (\ref{ea=sum}) implies (\ref{epa=0}).

Conversely, assume that (\ref{epa=0}) holds.
In view of Lemma \ref{lcancellation} (d),
we may assume that $p(x)$ is monic with $p(0)\ne 0$. Then
$$p(x)=(x-\alpha_{1})^{k_{1}}\cdots (x-\alpha_{r})^{k_{r}}.$$
For $1\le i\le r$, set
\begin{eqnarray}
p_{i}(x)=p(x)(x-\alpha_{i})^{-k_{i}}\in \C[x].
\end{eqnarray}
Since $p_{1}(x),\dots,p_{r}(x)$ are relatively prime, there are polynomials
$q_{1}(x),\dots,q_{r}(x)$ such that
\begin{eqnarray}\label{1=products}
1=p_{1}(x)q_{1}(x)+\cdots +p_{r}(x)q_{r}(x).
\end{eqnarray}
We have
\begin{eqnarray}\label{eA=sum}
A(x_{1},x_{2})=p_{1}(x_{1}/x_{2})q_{1}(x_{1}/x_{2})A(x_{1},x_{2})+\cdots 
+p_{r}(x_{1}/x_{2})q_{r}(x_{1}/x_{2})A(x_{1},x_{2}).
\end{eqnarray}
For $1\le i\le r$, since
$$(x_{1}/x_{2}-\alpha_{i})^{k_{i}}(p_{i}(x_{1}/x_{2})q_{i}(x_{1}/x_{2})A(x_{1},x_{2}))
=p(x_{1}/x_{2})q_{i}(x_{1}/x_{2})A(x_{1},x_{2})=0,$$
{}from [DLM] (cf. [K2]) we have
\begin{eqnarray}\label{esub-expre}
p_{i}(x_{1}/x_{2})q_{i}(x_{1}/x_{2})A(x_{1},x_{2})
=\sum_{j=0}^{k_{i}-1}\left(\frac{\alpha_{i}^{-j}}{j!}\left(\frac{\partial}{\partial x_{2}}\right)^{j}
x_{1}^{-1}\delta\left(\frac{\alpha_{i}x_{2}}{x_{1}}\right)\right)B_{i,j}(x_{2})
\end{eqnarray}
for some $B_{i,j}(x)\in U[[x,x^{-1}]]$ $(j=0,\dots,k_{i}-1)$. Now, combining (\ref{eA=sum}) with
(\ref{esub-expre}), we obtain (\ref{ea=sum}). 

For the uniqueness, let us assume (\ref{ea=sum}).
Notice that for $s\ne i,\; 0\le j\le k_{i}-1$, we have
$$p_{s}(x_{1}/x_{2})\left(\frac{\partial}{\partial x_{2}}\right)^{j}
x_{1}^{-1}\delta\left(\frac{\alpha_{i}x_{2}}{x_{1}}\right)=0,
$$
since $p_{s}(x)$ is a multiple of $(x-\alpha_{i})^{k_{i}}$.
Using this (twice) and (\ref{1=products}), for each $i$ we get
\begin{eqnarray}
& &p_{i}(x_{1}/x_{2})q_{i}(x_{1}/x_{2})A(x_{1},x_{2})\nonumber\\
&=&p_{i}(x_{1}/x_{2})q_{i}(x_{1}/x_{2})
\sum_{j=0}^{k_{i}-1}\left(\frac{\alpha_{i}^{-j}}{j!}\left(\frac{\partial}{\partial x_{2}}\right)^{j}
x_{1}^{-1}\delta\left(\frac{\alpha_{i}x_{2}}{x_{1}}\right)\right)B_{i,j}(x_{2})\nonumber\\
&=&\left(1-\sum_{s\ne i}p_{s}(x_{1}/x_{2})q_{s}(x_{1}/x_{2})\right)
\sum_{j=0}^{k_{i}-1}\left(\frac{\alpha_{i}^{-j}}{j!}\left(\frac{\partial}{\partial x_{2}}\right)^{j}
x_{1}^{-1}\delta\left(\frac{\alpha_{i}x_{2}}{x_{1}}\right)\right)B_{i,j}(x_{2})\nonumber\\
&=&\sum_{j=0}^{k_{i}-1}\left(\frac{\alpha_{i}^{-j}}{j!}\left(\frac{\partial}{\partial x_{2}}\right)^{j}
x_{1}^{-1}\delta\left(\frac{\alpha_{i}x_{2}}{x_{1}}\right)\right)B_{i,j}(x_{2}).
\end{eqnarray}
Then it follows from \cite{li-local} (cf. [LL]) that $B_{i,j}(x)$ are uniquely determined. 
\end{proof}

\section{Basic  operations on $\E(W)$}
In this section, for a vector space $W$
we introduce a natural action of the multiplicative group $\C^{\times}$ 
(of nonzero complex numbers) on $\E(W)$, a notion of compatibility 
for an ordered pair in $\E(W)$ and certain partial operations on $\E(W)$.
We then present certain properties analogous to the $\D$-bracket and 
$\D$-derivative properties for vertex algebras..

Let $W$ be a vector space, fixed throughout this section.
Following [LL] we set
\begin{eqnarray}
\E(W)=\Hom (W,W((x))).
\end{eqnarray}
We consider $\End W$ as a subspace of $\E(W)$ where each
linear operator on $W$ is considered as a constant series in $x$.
The identity operator on $W$, denoted by $1_{W}$,
is a special element of $\E(W)$.
Clearly, $\E(W)$ is closed under the formal derivative operator
$\frac{d}{dx}$. Set
\begin{eqnarray}
D={d\over dx}.
\end{eqnarray}
Just as in calculus, for $a(x)\in \E(W)$,
we use $a'(x)$ for the formal derivative of $a(x)$:
\begin{eqnarray}
a'(x)=\frac{d}{dx}a(x)=Da(x).
\end{eqnarray}
For any $\alpha\in \C^{\times}$ (nonzero complex numbers)
and $a(x)\in \E(W)$, it is clear that $a(\alpha x)\in \E(W)$.
For $\alpha\in \C^{\times}$, we define $R_{\alpha}\in \End (\E(W))$ by
\begin{eqnarray}
R_{\alpha}a(x)=a(\alpha x)\;\;\;\mbox{ for }a(x)\in \E(W)
\end{eqnarray}
(cf. [G-K-K]).
Alternatively we have
\begin{eqnarray}
R_{\alpha}=\alpha^{x{d\over dx}}.
\end{eqnarray}
It is clear that the map $R: \C^{\times}\rightarrow \End (\E(W))$
sending $\alpha$ to $R_{\alpha}$ 
is a representation of $\C^{\times}$ on $\E (W)$.
We have
\begin{eqnarray}\label{e-DRalpha}
D R_{\alpha}=\alpha R_{\alpha}D\;\;\;\mbox{ for }\alpha\in \C^{\times}.
\end{eqnarray}

We introduce the following notion of compatibility
(cf. \cite{li-G1}, [B2]):

\bd{dcompatibility-pair}
{\em  An ordered pair $(a(x),b(x))$ in $\E(W)$ 
is said to be {\em compatible} if there exists a nonzero polynomial
$f(x_{1},x_{2})$ such that
\begin{eqnarray}\label{ecompatibility-definition}
f(x_{1},x_{2})a(x_{1})b(x_{2})
\in \Hom (W,W((x_{1},x_{2}))).
\end{eqnarray}
A subset(space) $U$ of $\E(W)$ is said to be {\em pairwise compatible}
if every (ordered) pair in $U$ is compatible.}
\ed

\br{rdifference}
{\em Note that this compatibility is more general than 
the one defined in \cite{li-G1}
where polynomials $f(x_{1},x_{2})$ are 
of the special form $(x_{1}-x_{2})^{k}$.}
\er

\bl{ldefinition-prepare}
Let $(a(x),b(x))$ be a compatible (ordered) pair in $\E(W)$ and let
$f(x_{1},x_{2})$ be any nonzero polynomial such that
\begin{eqnarray}\label{efab-condition}
f(x_{1},x_{2})a(x_{1})b(x_{2})\in \Hom (W,W((x_{1},x_{2}))).
\end{eqnarray}
For any nonzero complex number $\alpha$, the expression
\begin{eqnarray}\label{eF-1-Fab}
\iota_{x,x_{0}}\left(f(x_{0}+\alpha x,x)^{-1}\right)
x_{1}^{-1}\delta\left(\frac{\alpha x+x_{0}}{x_{1}}\right)
\left(f(x_{1},x)a(x_{1})b(x)\right)
\end{eqnarray}
exists in $(\End W)[[x^{\pm 1},x_{0}^{\pm 1},x_{1}^{\pm 1}]]$ and 
in fact it lies in $\left(\Hom (W,W((x))((x_{0})))\right)[[x_{1}^{\pm 1}]]$.
Furthermore, it is independent of
the choice of $f(x_{1},x_{2})$ in (\ref{efab-condition}).
\el

\begin{proof} Set 
$$H(x_{1},x_{2})=f(x_{1},x_{2})a(x_{1})b(x_{2})\in \Hom (W,W((x_{1},x_{2}))).
$$
{}From Remark \ref{r-delta-function} we have
$$x_{1}^{-1}\delta\left(\frac{\alpha x+x_{0}}{x_{1}}\right)H(x_{1},x)
=x_{1}^{-1}\delta\left(\frac{\alpha x+x_{0}}{x_{1}}\right)H(\alpha x+x_{0},x),$$
which lies in $\left(\Hom (W,W((x))[[x_{0}]])\right)[[x_{1}^{\pm 1}]]$.
As $\iota_{x,x_{0}}\left(f(x_{0}+\alpha x,x)^{-1}\right)\in
\C((x))((x_{0}))$, we have
$$\iota_{x,x_{0}}\left(f(x_{0}+\alpha x,x)^{-1}\right)
x_{1}^{-1}\delta\left(\frac{\alpha x+x_{0}}{x_{1}}\right)H(\alpha
x+x_{0},x)\in \left(\Hom (W,W((x))((x_{0})))\right)[[x_{1}^{\pm 1}]].$$

Assume that $g(x_{1},x_{2})$ is another nonzero polynomial such that
$$g(x_{1},x_{2})a(x_{1})b(x_{2})\in \Hom (W,W((x_{1},x_{2}))).$$
Using Lemma \ref{lcancellation} we have
\begin{eqnarray*}
& &\iota_{x,x_{0}}\left(f(x_{0}+\alpha x,x)^{-1}
g(x_{0}+\alpha x,x)^{-1}\right)
x_{1}^{-1}\delta\left(\frac{\alpha x+x_{0}}{x_{1}}\right)
\left(f(x_{1}, x)g(x_{1},x)a(x_{1})b(x)\right)\\
&=&\iota_{x,x_{0}}\left(f(x_{0}+\alpha x,x)^{-1}
g(x_{0}+\alpha x,x)^{-1}\right)
x_{1}^{-1}\delta\left(\frac{\alpha x+x_{0}}{x_{1}}\right)
f(x_{1},x)\left(g(x_{1},x)a(x_{1})b(x)\right)\\
&=&\iota_{x,x_{0}}\left(f(x_{0}+\alpha x,x)^{-1}\right)
\iota_{x,x_{0}}\left(g(x_{0}+\alpha x,x)^{-1}\right)
x_{1}^{-1}\delta\left(\frac{\alpha x+x_{0}}{x_{1}}\right)f(x_{0}+\alpha x,x)\\
& &\ \ \cdot\left(g(x_{1},x)a(x_{1})b(x)\right)\\
&=&\iota_{x,x_{0}}\left(g(x_{0}+\alpha x,x)^{-1}\right)
x_{1}^{-1}\delta\left(\frac{\alpha x+x_{0}}{x_{1}}\right)
\left(g(x_{1},x)a(x_{1})b(x)\right),
\end{eqnarray*}
noticing that all the assumptions for associativity are satisfied.
Symmetrically, we have
\begin{eqnarray*}
& &\iota_{x,x_{0}}\left(f(x_{0}+\alpha x,x)^{-1}g(x_{0}+\alpha x,x)^{-1}\right)
x_{1}^{-1}\delta\left(\frac{\alpha x+x_{0}}{x_{1}}\right)
\left(f(x_{1},x)g(x_{1},x)a(x_{1})b(x)\right)\\
&=&\iota_{x,x_{0}}\left(f(x_{0}+\alpha x,x)^{-1}\right)
x_{1}^{-1}\delta\left(\frac{\alpha x+x_{0}}{x_{1}}\right)
\left(f(x_{1},x)a(x_{1})b(x)\right).
\end{eqnarray*}
It now follows that the expression in (\ref{eF-1-Fab}) is independent of 
the choice of $f(x_{1},x_{2})$.
\end{proof}

Now, we introduce the following partial operations on $\E(W)$.

\bd{d-operation}
{\em Let $(a(x),b(x))$ be a compatible (ordered) pair
in $\E(W)$. For $\alpha\in \C^{\times},\; n\in \Z$, we define
$a(x)_{(\alpha,n)}b(x)\in (\End W)[[x,x^{-1}]]$ in terms of generating
function
\begin{eqnarray}
\Y_{\alpha}(a(x),x_{0})b(x)=\sum_{n\in \Z}(a(x)_{(\alpha,n)}b(x)) x_{0}^{-n-1}
\in (\End W)[[x_{0}^{\pm 1},x^{\pm 1}]]
\end{eqnarray}
by
\begin{eqnarray}
& &\Y_{\alpha}(a(x),x_{0})b(x)\nonumber\\
&=&\Res_{x_{1}}\iota_{x,x_{0}}\left(f(x_{0}+\alpha x,x)^{-1}\right)
x_{1}^{-1}\delta\left(\frac{\alpha x+x_{0}}{x_{1}}\right)
\left(f(x_{1},x)a(x_{1})b(x)\right)\\
&=&\iota_{x,x_{0}}\left(f(x_{0}+\alpha x,x)^{-1}\right)
\left(f(x_{1},x)a(x_{1})b(x)\right)|_{x_{1}=\alpha x+x_{0}},\label{edef-3.15}
\end{eqnarray}
where $f(x_{1},x_{2})$ is any nonzero polynomial such that
$$f(x_{1},x_{2})a(x_{1})b(x_{2})\in \Hom (W,W((x_{1},x_{2}))).$$
We particularly set
\begin{eqnarray}
& &\Y(a(x),x_{0})b(x)=\Y_{1}(a(x),x_{0})b(x),\\
& &a(x)_{n}b(x)=a(x)_{(1,n)}b(x)\;\;\;\mbox{  for }n\in \Z.
\end{eqnarray}}
\ed

\br{rsubtlty}
{\em In (\ref{edef-3.15}), we are not allowed to write
$$\left(f(x_{1},x)a(x_{1})b(x)\right)|_{x_{1}=\alpha x+x_{0}}
=f(\alpha x+x_{0},x)a(\alpha x+x_{0})b(x),$$ 
since $a(\alpha x+x_{0})b(x)$ in general does not exist
in $(\End W)[[x^{\pm 1},x_{0}^{\pm 1}]]$. On the other hand,
$f(\alpha x+x_{0},x)a(x_{0}+\alpha x)b(x)$ exists.
But, 
\begin{eqnarray*}
\left(f(x_{1},x)a(x_{1})b(x)\right)|_{x_{1}=\alpha x+x_{0}}
&\ne& \left(f(x_{1},x)a(x_{1})b(x)\right)|_{x_{1}=x_{0}+\alpha x}\\
& &\left(=f(x_{0}+\alpha x,x)a(x_{0}+\alpha x)b(x)\right),
\end{eqnarray*}
unless $f(x_{1},x)a(x_{1})b(x)$ involves only nonnegative powers of $x_{1}$.}
\er

First, we have:

\bp{pexistence} 
Let $(a(x),b(x))$ be a compatible (ordered) pair
in $\E(W)$. Then
\begin{eqnarray}\label{eew-closed}
a(x)_{(\alpha,n)}b(x)\in \E (W)\;(=\Hom (W,W((x))))
\;\;\;\mbox{ for }\alpha\in \C^{\times},\; n\in \Z.
\end{eqnarray}
Let $f(x_{1},x_{2})$ be any nonzero polynomial such that
\begin{eqnarray}\label{efabcondition}
f(x_{1},x_{2})a(x_{1})b(x_{2})\in \Hom(W,W((x_{1},x_{2})))
\end{eqnarray}
and let $k$ be an integer such that
$$x_{0}^{k}\iota_{x_{2},x_{0}}
\left(f(x_{0}+\alpha x_{2},x_{2})^{-1}\right)\in \C((x_{2}))[[x_{0}]].$$
Then
\begin{eqnarray}\label{e-truncation-ab}
a(x)_{(\alpha,n)}b(x)=0\;\;\;\mbox{ for }n\ge k,
\end{eqnarray}
Furthermore, for $w\in W$, let $l$ be a nonnegative integer depending on $w$ such that
$$x_{1}^{l}f(x_{1},x_{2})a(x_{1})b(x_{2})w\in W[[x_{1},x_{2}]][x_{2}^{-1}],$$
then
\begin{eqnarray}\label{eweakassoc-def}
& &(x_{0}+\alpha x_{2})^{l}f(x_{0}+\alpha x_{2},x_{2})
(\Y_{\alpha}(a(x_{2}),x_{0})b(x_{2}))w\nonumber\\
&=&(x_{0}+\alpha x_{2})^{l}f(x_{0}+\alpha x_{2},x_{2})
a(x_{0}+\alpha x_{2})b(x_{2})w.
\end{eqnarray}
\ep

\begin{proof} For any nonzero polynomial $f(x_{1},x_{2})$ satisfying
(\ref{efabcondition}), by Lemma \ref{ldefinition-prepare}
we have
$$(f(x_{1},x)a(x_{1})b(x))|_{x_{1}=\alpha x+x_{0}}
\in \Hom (W,W((x))[[x_{0}]])=(\Hom (W,W((x))))[[x_{0}]].$$
With 
$\iota_{x,x_{0}}\left(f(x_{0}+\alpha x,x)^{-1}\right)\in
\C((x))((x_{0}))$,
we have
$$\Y_{\alpha}(a(x),x_{0})b(x)\in \left(\Hom (W,W((x)))\right)((x_{0})),$$
proving (\ref{eew-closed}). 
This also proves (\ref{e-truncation-ab}), as
\begin{eqnarray*}
x_{0}^{k}\Y_{\alpha}(a(x),x_{0})b(x)
&=&x_{0}^{k}\iota_{x,x_{0}}\left(f(x_{0}+\alpha x,x)^{-1}\right)
(f(x_{1},x)a(x_{1})b(x))|_{x_{1}=\alpha x+x_{0}}\\
&\in&(\Hom (W,W((x))))[[x_{0}]].
\end{eqnarray*}

For $w\in W$, let $\ell\in \N$ (depending on $w$) such that
$$x_{1}^{l}f(x_{1},x_{2})a(x_{1})b(x_{2})w\in W[[x_{1},x_{2}]][x_{2}^{-1}].$$
Then
\begin{eqnarray*}
& &\left(x_{1}^{l}f(x_{1},x_{2})a(x_{1})b(x_{2})w\right)|_{x_{1}=\alpha
x+x_{0}}\nonumber\\
&=&\left(x_{1}^{l}f(x_{1},x_{2})a(x_{1})b(x_{2})w\right)|_{x_{1}=x_{0}+\alpha
x}\nonumber\\
&=&(x_{0}+\alpha x)^{l}f(x_{0}+\alpha x,x)a(x_{0}+\alpha x)b(x)w.
\end{eqnarray*}
Then
\begin{eqnarray}\label{eYalpha-formula}
& &(x_{0}+\alpha x)^{l}f(x_{0}+\alpha x,x)
\left(\Y_{\alpha}(a(x),x_{0})b(x)\right)w\nonumber\\
&=&\Res_{x_{1}}(x_{0}+\alpha x)^{l}
x_{1}^{-1}\delta\left(\frac{\alpha x+x_{0}}{x_{1}}\right)
\left(f(x_{1},x)a(x_{1})b(x)w\right)\nonumber\\
&=&\Res_{x_{1}}x_{1}^{l}
x_{1}^{-1}\delta\left(\frac{\alpha x+x_{0}}{x_{1}}\right)
\left(f(x_{1},x)a(x_{1})b(x)w\right)\nonumber\\
&=&\left(x_{1}^{l}f(x_{1},x)a(x_{1})b(x)w\right)|_{x_{1}=\alpha
x+x_{0}}\nonumber\\
&=&\left(x_{1}^{l}f(x_{1},x)a(x_{1})b(x)w\right)|_{x_{1}=x_{0}+\alpha
x}\nonumber\\
&=&(x_{0}+\alpha x)^{l}f(x_{0}+\alpha x,x)
a(x_{0}+\alpha x)b(x)w,
\end{eqnarray}
proving (\ref{eweakassoc-def}). 
Now, the proof is complete.
\end{proof}

The following result says that we can express $\Y_{\alpha}(a(x),x_{0})b(x)$
in terms of matrix-coefficients and iota-maps:

\bl{lmatrix-coeff}
Let $(a(x),b(x))$ be a compatible (ordered) pair in $\E(W)$. 
Then for any nonzero complex number $\alpha$ and
for any $w^{*}\in W^{*},\; w\in W$, the formal series
$$\<w^{*},a(x_{0}+\alpha x_{2})b(x_{2})w\>,$$
an element of $\C((x_{0}))((x_{2}))$, 
lies in the image of $\iota_{x_{0},x_{2}}$. 
Furthermore,
\begin{eqnarray}
\<w^{*},(\Y_{\alpha}(a(x),x_{0})b(x))w\>
=\iota_{x,x_{0}}\iota_{x_{0},x}^{-1}
\<w^{*},a(x_{0}+\alpha x)b(x)w\>.
\end{eqnarray}
\el

\begin{proof} It follows from (\ref{eweakassoc-def}) as in the usual case 
(cf. [FLM], [FHL], [DL], [LL]).
\end{proof}

The following are analogues of the creation property and the vacuum property
in the definition of the notion of vertex algebra:

\bp{pvacuum}
For any $a(x)\in \E(W)$, the sequences $(1_{W}, a(x))$ and
$(a(x), 1_{W})$ are compatible and for $\alpha\in \C^{\times}$ we have
\begin{eqnarray}
& &\Y_{\alpha}(1_{W},x_{0})a(x)=a(x),\label{e-vacuum-calculus}\\
& &\Y_{\alpha}(a(x),x_{0})1_{W}=a(\alpha x+x_{0})
=e^{\alpha^{-1}x_{0}{d\over dx}}a(\alpha x)
=e^{\alpha^{-1}x_{0}{d\over dx}}R_{\alpha}a(x).\label{e-creation-calculus}
\end{eqnarray}
In particular,
\begin{eqnarray}
& &a(x)_{(\alpha,-1)}1_{W}=R_{\alpha}a(x),\label{ea-alpha-1-R}\\
& &a(x)_{-2}1_{W}\;\left(=a(x)_{(1,-2)}1_{W}\right)
=\frac{d}{dx}a(x)=Da(x).\label{ea-2=da}
\end{eqnarray}
\ep

\begin{proof} The compatibility assertion is clear.
For $w^{*}\in W^{*},\;w\in W$, using Lemma \ref{lmatrix-coeff} we have
\begin{eqnarray}
\< w^{*},\left(\Y_{\alpha}(1_{W},x_{0})a(x)\right)w\>
&=&\iota_{x,x_{0}}\iota_{x_{0},x}^{-1}\<w^{*},1_{W}(x_{0}+\alpha x)a(x)w\>
\nonumber\\
&=&\iota_{x,x_{0}}\iota_{x_{0},x}^{-1}\<w^{*},a(x)w\>\nonumber\\
&=&\<w^{*},a(x)w\>,
\end{eqnarray}
proving (\ref{e-vacuum-calculus}). Similarly,
for $w^{*}\in W^{*},\;w\in W$, we have
\begin{eqnarray}
\<w^{*},\left(\Y_{\alpha}(a(x),x_{0})1_{W}\right)w\>
&=&\iota_{x,x_{0}}\iota_{x_{0},x}^{-1}
\<w^{*},a(x_{0}+\alpha x)w\>\nonumber\\
&=&\<w^{*},a(\alpha x+x_{0})w\>,
\end{eqnarray}
proving (\ref{e-creation-calculus}).
\end{proof}

We also have the following properties, analogous to the $\D$-properties 
for vertex algebras (cf. [LL]):

\bp{pbasic-property2}
Let $(a(x),b(x))$ be a compatible (ordered) pair in $\E(W)$.
Then the ordered pairs $(a'(x),b(x)), (a(x),b'(x))$ and $(a(\alpha x), b(\beta x))$
for $\alpha,\beta\in \C^{\times}$ are compatible.
Furthermore, for $\alpha,\beta\in \C^{\times}$, we have
\begin{eqnarray}
& &\Y_{\alpha}(Da(x),x_{0})b(x)={\partial\over \partial x_{0}}
\Y_{\alpha}(a(x),x_{0})b(x),\label{e-Dderivative-calculus}\\
& & D\Y_{\alpha}(a(x),x_{0})b(x)-\Y_{\alpha}(a(x),x_{0})D b(x)
=\alpha {\partial\over \partial x_{0}}
\Y_{\alpha}(a(x),x_{0})b(x),\label{e-Dbracket-calculus}\\
& &R_{\alpha}\Y_{\beta}(a(x),x_{0})b(x)
=\Y_{\alpha\beta}(a(x),x_{0})R_{\alpha}b(x),
\label{e-conjugation-calculus}\\
& &\Y_{\alpha}(R_{\beta}a(x),\beta^{-1}x_{0})b(x)
=\Y_{\alpha\beta}(a(x),x_{0})b(x),
\label{e-conjugation-calculus-2}
\end{eqnarray}
In particular,
\begin{eqnarray}\label{e-conjugation-calculus-3}
\Y_{\alpha}(a(x),x_{0})b(x)=R_{\alpha}\Y_{1}(a(x),x_{0})R_{\alpha^{-1}}b(x)
=\Y_{1}(R_{\alpha}a(x),\alpha^{-1}x_{0}).
\end{eqnarray}
\ep

\begin{proof} Let $0\ne f(x_{1},x_{2})\in \C[x_{1},x_{2}]$ such that
$$f(x_{1},x_{2})a(x_{1})b(x_{2})\in \Hom (W,W((x_{1},x_{2}))).$$
We immediately have
$$f(\alpha x_{1},\beta x_{2})a(\alpha x_{1})b(\beta x_{2})\in \Hom (W,W((x_{1},x_{2}))).$$
We also have
\begin{eqnarray}
f(x_{1},x_{2})^{2}a(x_{1})b(x_{2}),\;\;
f(x_{1},x_{2})f_{x_{1}}(x_{1},x_{2})a(x_{1})b(x_{2})\in \Hom (W,W((x_{1},x_{2}))).
\end{eqnarray}
Then
\begin{eqnarray}
& &f(x_{1},x_{2})^{2}a'(x_{1})b(x_{2})
=\frac{\partial}{\partial x_{2}}\left(f(x_{1},x_{2})^{2}a'(x_{1})b(x_{2})\right)
-2f(x_{1},x_{2})f_{x_{1}}(x_{1},x_{2})a(x_{1})b(x_{2})\nonumber\\
&\in& \Hom (W,W((x_{1},x_{2}))).
\end{eqnarray}
This proves the compatibility assertions.
For simplicity, we locally use $\partial_{x}$ 
for $\frac{\partial}{\partial x}$ in this proof.
Let $w^{*}\in W^{*},\; w\in W$. Using Lemma \ref{lmatrix-coeff} we have
\begin{eqnarray}
\<w^{*},(\Y_{\alpha}(a'(x),x_{0})b(x))w\>
&=&\iota_{x,x_{0}}\iota_{x_{0},x}^{-1}
\<w^{*},a'(x_{0}+\alpha x)b(x)w\>\nonumber\\
&=&\iota_{x,x_{0}}\iota_{x_{0},x}^{-1}\partial_{x_{0}}
\<w^{*},a(x_{0}+\alpha x)b(x)w\>\nonumber\\
&=&\partial_{x_{0}}\iota_{x,x_{0}}\iota_{x_{0},x}^{-1}
\<w^{*},a(x_{0}+\alpha x)b(x)w\>\nonumber\\
&=&\partial_{x_{0}}\<w^{*},(\Y_{\alpha}(a(x),x_{0})b(x))w\>,
\end{eqnarray}
proving (\ref{e-Dderivative-calculus}). Similarly, we have
\begin{eqnarray}
& &\<w^{*},(D\Y_{\alpha}(a(x),x_{0})b(x))w\>\nonumber\\
&=&\partial_{x}\<w^{*},(\Y_{\alpha}(a(x),x_{0})b(x))w\>\nonumber\\
&=&\partial_{x}\iota_{x,x_{0}}\iota_{x_{0},x}^{-1}
\<w^{*},a(x_{0}+\alpha x)b(x)w\>\nonumber\\
&=&\iota_{x,x_{0}}\iota_{x_{0},x}^{-1}
\<w^{*},\alpha a'(x_{0}+\alpha x)b(x)+a(x_{0}+\alpha x)b'(x)w\>\nonumber\\
&=&\iota_{x,x_{0}}\iota_{x_{0},x}^{-1}\alpha \partial_{x_{0}}
\<w^{*},a(x_{0}+\alpha x)b(x)w\>
+\iota_{x,x_{0}}\iota_{x_{0},x}^{-1}
\<w^{*},a(x_{0}+\alpha x)b'(x)w\>\nonumber\\
&=&\alpha \partial_{x_{0}}\<w^{*},(\Y_{\alpha}(a(x),x_{0})b(x))w\>
+\<w^{*},(\Y_{\alpha}(a(x),x_{0})b'(x))w\>,
\end{eqnarray}
proving (\ref{e-Dbracket-calculus}).
For $\beta\in \C^{\times}$, let $R_{\beta}$ naturally act on $\E(W)[[x_{0}^{\pm 1}]]$.
We have
\begin{eqnarray}
& &\<w^{*},(R_{\alpha}\Y_{\beta}(a(x),x_{0})b(x))w\>\nonumber\\
&=&\alpha^{x\frac{d}{dx}}\<w^{*},(\Y_{\beta}(a(x),x_{0})b(x))w\>\nonumber\\
&=&\alpha^{x\frac{d}{dx}}\iota_{x,x_{0}}\iota_{x_{0},x}^{-1}
\<w^{*},a(x_{0}+\beta x)b(x)w\>\nonumber\\
&=&\iota_{x,x_{0}}\iota_{x_{0},x}^{-1}
\<w^{*},a(x_{0}+\alpha\beta x)b(\alpha x)w\>\nonumber\\
&=&\iota_{x,x_{0}}\iota_{x_{0},x}^{-1}
\<w^{*},a(x_{0}+\alpha\beta x)(R_{\alpha}b(x))w\>\nonumber\\
&=&\<w^{*},(\Y_{\alpha\beta}(a(x),x_{0})R_{\alpha}b(x))w\>.
\end{eqnarray}
This proves (\ref{e-conjugation-calculus}). The identity (\ref{e-conjugation-calculus-2})
also holds, since
\begin{eqnarray}
& &\<w^{*},(\Y_{\alpha}(R_{\beta}(a(x),\beta^{-1}x_{0})b(x))w\>\nonumber\\
&=&\iota_{x,x_{0}}\iota_{x_{0},x}^{-1}
\<w^{*},(R_{\beta}a)(\beta^{-1}x_{0}+\alpha x)b(x)w\>\nonumber\\
&=&\iota_{x,x_{0}}\iota_{x_{0},x}^{-1}
\<w^{*},a(x_{0}+\alpha\beta x)b(x))w\>\nonumber\\
&=&\<w^{*},(\Y_{\alpha\beta}(a(x),x_{0})b(x))w\>.
\end{eqnarray}
Specializing $\alpha,\beta$ appropriately we obtain the particular case;
specialize $\beta=1$ in (\ref{e-conjugation-calculus}) to get the first equality
and substitute $(\alpha,\beta)$ with $(1,\alpha)$ in (\ref{e-conjugation-calculus-2})
to get the outside equality.
\end{proof}

The following result gives Jacobi-like identities:

\bp{pjacobi-definition}
Let $a(x),b(x)\in \E(W)$. Suppose that there exist
$$0\ne f(x_{1},x_{2})\in \C[x_{1},x_{2}],\; K(x_{2},x_{1})\in \Hom (W,W((x_{2}))((x_{1})))$$
such that
\begin{eqnarray}\label{e-quantum-relation}
f(x_{1},x_{2})a(x_{1})b(x_{2})=f(x_{1},x_{2})K(x_{2},x_{1}).
\end{eqnarray}
Then $(a(x),b(x))$ is compatible and for any $\alpha\in \C^{\times}$,
\begin{eqnarray}\label{e-jacobi-calculus}
& &x_{1}^{-1}\delta\left(\frac{\alpha x+x_{0}}{x_{1}}\right)
f(x_{1},x)\Y_{\alpha}(a(x),x_{0})b(x)\nonumber\\
&=&x_{0}^{-1}\delta\left(\frac{x_{1}-\alpha x}{x_{0}}\right)
f(x_{1},x)a(x_{1})b(x)
-x_{0}^{-1}\delta\left(\frac{\alpha x-x_{1}}{-x_{0}}\right)f(x_{1},x)K(x,x_{1}).
\end{eqnarray}
\ep

\begin{proof} Noticing that the expression on the left hand side of
(\ref{e-quantum-relation}) lies in $\Hom(W,W((x_{1}))((x_{2})))$ and
the expression on the right hand side lies in
$\Hom(W,W((x_{2}))((x_{1})))$, we have
$$f(x_{1},x_{2})a(x_{1})b(x_{2})\in \Hom (W,W((x_{1},x_{2}))),$$
proving that $(a(x),b(x))$ is compatible. Let $w\in W$ be arbitrarily fixed. 
There exists $l\in \N$ such that
$$x_{1}^{l}f(x_{1},x_{2})a(x_{1})b(x_{2})w\in W[[x_{1},x_{2}]][x_{2}^{-1}].$$
By Proposition \ref{pexistence} we have
\begin{eqnarray}
& &(x_{0}+\alpha x)^{l}f(x_{0}+\alpha x,x)
\left(\Y_{\alpha}(a(x),x_{0})b(x)\right)w\nonumber\\
&=&(x_{0}+\alpha x)^{l}f(x_{0}+\alpha x,x)a(x_{0}+\alpha x)b(x)w.
\end{eqnarray}
By Lemma \ref{lformaljacobiidentity} with
$$A=f(x_{1},x_{2})a(x_{1})b(x_{2})w,\;\;
B=f(x_{1},x_{2})K(x_{2},x_{1})w$$
and
$$C=f(x_{0}+\alpha x,x)a(x_{0}+\alpha x)b(x)w$$
we have
\begin{eqnarray}
& &x_{1}^{-1}\delta\left(\frac{\alpha x+x_{0}}{x_{1}}\right)
f(x_{0}+\alpha x,x)
\left(\Y_{\alpha}(a(x),x_{0})b(x)\right)w\nonumber\\
&=&x_{0}^{-1}\delta\left(\frac{x_{1}-\alpha x}{x_{0}}\right)
f(x_{1},x)a(x_{1})b(x)w
-x_{0}^{-1}\delta\left(\frac{-\alpha x+
x_{1}}{x_{0}}\right)f(x_{1},x)K(x,x_{1})w.\hspace{1cm}
\end{eqnarray}
Since (by (\ref{e-delta-2}))
$$x_{1}^{-1}\delta\left(\frac{\alpha x+x_{0}}{x_{1}}\right)
f(x_{0}+\alpha x,x)
=x_{1}^{-1}\delta\left(\frac{\alpha x+x_{0}}{x_{1}}\right)
f(x_{1},x),$$
we immediately have (\ref{e-jacobi-calculus}).
\end{proof}

We also have the following result:

\bp{pcommutator-formula}
Let $a(x),b(x)\in \E(W)$ be such that
\begin{eqnarray}
p(x_{1}/x_{2})a(x_{1})b(x_{2})=p(x_{1}/x_{2})K(x_{2},x_{1})
\end{eqnarray}
for some nonzero polynomial $p(x)$ and for some
$K(x_{2},x_{1})\in \Hom (W,W((x_{2}))((x_{1})))$. Then $(a(x),b(x))$ is compatible and
\begin{eqnarray}\label{e-commutator-calculus-general}
a(x_{1})b(x_{2})-K(x_{2},x_{1})
=\sum_{\alpha\in \C^{\times}} \Res_{x_{0}}
x_{1}^{-1}\delta\left(\frac{\alpha x_{2}+x_{0}}{x_{1}}\right)
\Y_{\alpha}(a(x_{2}),x_{0})b(x_{2}),
\end{eqnarray}
which is a finite sum over only nonzero roots $\alpha$ of $p(x)$.
\ep

\begin{proof} The ordered pair $(a(x),b(x))$ is compatible because
$$x_{2}^{\deg p(x)}p(x_{1}/x_{2})a(x_{1})b(x_{2})
=x_{2}^{\deg p(x)}p(x_{1}/x_{2})K(x_{2},x_{1})$$
with $0\ne x_{2}^{\deg p(x)}p(x_{1}/x_{2})\in \C[x_{1},x_{2}]$.
If $p(x)=x^{k}q(x)$ where $k\in \N,\; q(x)\in \C[x]$,
then $q(x_{1}/x_{2})a(x_{1})b(x_{2})=q(x_{1}/x_{2})K(x_{2},x_{1})$.
In view of this, we may assume that $p(0)\ne 0$.
Let $\alpha_{1},\dots,\alpha_{r}$ be the distinct (nonzero) roots of $p(x)$ of multiplicities
$k_{1},\dots,k_{r}$.
For $1\le i\le r$, let $p_{i}(x)$ be the (unique) polynomial such that
\begin{eqnarray}
p(x)=p_{i}(x)(x-\alpha_{i})^{k_{i}}\;\;\;(\mbox{with }\;p_{i}(\alpha_{i})\ne 0).
\end{eqnarray}
Notice that (\ref{e-jacobi-calculus}) with $f(x_{1},x_{2})$
being replaced by $p(x_{1}/x_{2})$ holds.
For $1\le i\le r$, multiplying by $(x_{2}/x_{0})^{k_{i}}$ and then
taking $\Res_{x_{0}}$ (of (\ref{e-jacobi-calculus})) we get
\begin{eqnarray}\label{emodified-comm}
& &p_{i}(x_{1}/x_{2})(a(x_{1})b(x_{2})-K(x_{2},x_{1}))\nonumber\\
&=&\Res_{x_{0}}p_{i}(x_{1}/x_{2})
x_{1}^{-1}\delta\left(\frac{\alpha_{i} x_{2}+x_{0}}{x_{1}}\right)
\Y_{\alpha_{i}}(a(x_{2}),x_{0})b(x_{2}).
\end{eqnarray}
Since $p_{1}(x),\dots, p_{r}(x)$ 
are relatively prime, there are polynomials $q_{1}(x),\dots, q_{r}(x)$ such that
\begin{eqnarray}\label{e1=sum}
1=q_{1}(x)p_{1}(x)+\cdots +q_{r}(x)p_{r}(x).
\end{eqnarray}
Using this and (\ref{emodified-comm}) we get
\begin{eqnarray}
& &a(x_{1})b(x_{2})-K(x_{2},x_{1})\nonumber\\
&=&\sum_{i=1}^{r} 
q_{i}(x_{1}/x_{2})p_{i}(x_{1}/x_{2})(a(x_{1})b(x_{2})-K(x_{2},x_{1}))\nonumber\\
&=&\sum_{i=1}^{r}\Res_{x_{0}}q_{i}(x_{1}/x_{2})p_{i}(x_{1}/x_{2})
x_{1}^{-1}\delta\left(\frac{\alpha_{i} x_{2}+x_{0}}{x_{1}}\right)
\Y_{\alpha_{i}}(a(x_{2}),x_{0})b(x_{2}).
\end{eqnarray}
In view of Proposition \ref{pexistence} we have
\begin{eqnarray}\label{e-lowest-power}
x_{0}^{k_{i}}\Y_{\alpha_{i}}(a(x),x_{0})b(x)\in \E(W)[[x_{0}]]
\;\;\mbox{ for }i=1,\dots,r. 
\end{eqnarray}
For $i\ne j$, since
$p_{i}(\alpha_{j}+x_{0}/x)\in x_{0}^{k_{j}}\C[x,x^{-1},x_{0}]$,
using (\ref{e-delta-3}) and (\ref{e-lowest-power}) we get
\begin{eqnarray}\label{ei-not-j}
& &\Res_{x_{0}}q_{i}(x_{1}/x_{2})p_{i}(x_{1}/x_{2})
x_{1}^{-1}\delta\left(\frac{\alpha_{j} x_{2}+x_{0}}{x_{1}}\right)
\Y_{\alpha_{j}}(a(x_{2}),x_{0})b(x_{2})\nonumber\\
&=&\Res_{x_{0}}q_{i}(\alpha_{j}+x_{0}/x_{2})
p_{i}(\alpha_{j} +x_{0}/x_{2})
x_{1}^{-1}\delta\left(\frac{\alpha_{j} x_{2}+x_{0}}{x_{1}}\right)
\Y_{\alpha_{j}}(a(x_{2}),x_{0})b(x_{2})\nonumber\\
&=&0.
\end{eqnarray}
Using (\ref{e1=sum}) and (\ref{ei-not-j}) we get
\begin{eqnarray}
& &\Res_{x_{0}}q_{j}(x_{1}/x_{2})p_{j}(x_{1}/x_{2})
x_{1}^{-1}\delta\left(\frac{\alpha_{j} x_{2}+x_{0}}{x_{1}}\right)
\Y_{\alpha_{j}}(a(x_{2}),x_{0})b(x_{2})\nonumber\\
&=&\sum_{i=1}^{r}\Res_{x_{0}}q_{i}(x_{1}/x_{2})p_{i}(x_{1}/x_{2})
x_{1}^{-1}\delta\left(\frac{\alpha_{j} x_{2}+x_{0}}{x_{1}}\right)
\Y_{\alpha_{j}}(a(x_{2}),x_{0})b(x_{2})\nonumber\\
&=&\Res_{x_{0}}
x_{1}^{-1}\delta\left(\frac{\alpha_{j} x_{2}+x_{0}}{x_{1}}\right)
\Y_{\alpha_{j}}(a(x_{2}),x_{0})b(x_{2}),
\end{eqnarray}
so that
\begin{eqnarray}
a(x_{1})b(x_{2})-K(x_{2},x_{1})
=\sum_{i=1}^{r}\Res_{x_{0}}
x_{1}^{-1}\delta\left(\frac{\alpha_{i} x_{2}+x_{0}}{x_{1}}\right)
\Y_{\alpha_{i}}(a(x_{2}),x_{0})b(x_{2}).
\end{eqnarray}
Since $\Y_{\alpha}(a(x),x_{0})b(x)\in \E(W)[[x_{0}]]$ for $\alpha\ne
\alpha_{1},\dots,\alpha_{r}$,
we immediately have (\ref{e-commutator-calculus-general}).
\end{proof}

The following is an important consequence
of Proposition \ref{pcommutator-formula}:

\bc{cquantum-commutator-formula}
Let $a(x),b(x)\in \E(W)$. Suppose that there are $0\ne p(x)\in \C[x]$,
$f_{i}(x)\in \C((x))$ and $c_{i}(x),d_{i}(x)\in \E(W)$ for $i=1,\dots,r$ 
such that
\begin{eqnarray}
p(x_{1}/x_{2})a(x_{1})b(x_{2})=p(x_{1}/x_{2})\sum_{i=1}^{r}f_{i}(x_{1}/x_{2})c_{i}(x_{2})d_{i}(x_{1}).
\end{eqnarray}
Then $(a(x),b(x))$ is compatible and
\begin{eqnarray}\label{e-diff-calculus}
& &a(x_{1})b(x_{2})-\sum_{i=1}^{r}f_{i}(x_{1}/x_{2})c_{i}(x_{2})d_{i}(x_{1})\nonumber\\
&=&\sum_{\alpha\in \C^{\times}} \Res_{x_{0}}
x_{1}^{-1}\delta\left(\frac{\alpha x_{2}+x_{0}}{x_{1}}\right)
\Y_{\alpha}(a(x_{2}),x_{0})b(x_{2}),
\end{eqnarray}
which is a finite sum over only nonzero roots of $p(x)$.
In particular, if
\begin{eqnarray}
p(x_{1}/x_{2})a(x_{1})b(x_{2})=p(x_{1}/x_{2})b(x_{2})a(x_{1})
\end{eqnarray}
for some nonzero polynomial $p(x)$, then $(a(x),b(x))$ is compatible and
\begin{eqnarray}\label{e-commutator-calculus}
[a(x_{1}),b(x_{2})]
=\sum_{\alpha\in \C^{\times}} \Res_{x_{0}}
x_{1}^{-1}\delta\left(\frac{\alpha x_{2}+x_{0}}{x_{1}}\right)
\Y_{\alpha}(a(x_{2}),x_{0})b(x_{2}),
\end{eqnarray}
which is a finite sum over only nonzero roots of $p(x)$.
\ec

In the same spirit we have the following result which
will be very useful in determining the structure of certain vertex algebras
later: 

\bp{pdecomposition}
Let $W$ be a vector space as before and let $a(x),b(x)\in \E(W)$. Assume that
there exists $K(x_{2},x_{1})\in \Hom(W,W((x_{2}))((x_{1})))$ such that
\begin{eqnarray}\label{eab-commuator-decomposition}
a(x_{1})b(x_{2})-K(x_{2},x_{1})=\sum_{i=1}^{r}\sum_{j=0}^{k_{i}}c_{i,j}(x_{2})\frac{1}{j!} 
\left(\alpha_{i}^{-1}\frac{\partial}{\partial x_{2}}\right)^{j}
x_{1}^{-1}\delta\left(\frac{\alpha_{i}x_{2}}{x_{1}}\right),
\end{eqnarray}
where $\alpha_{1},\dots,\alpha_{r}$ are finitely many distinct nonzero complex numbers,
$k_{i}$ are nonnegative integers
and $c_{i,j}(x)\in \E(W)$. Then $(a(x),b(x))$ is compatible and
\begin{eqnarray}
& &a(x)_{(\alpha,n)}b(x)=0\;\;\;\ \ \ \ \ \ \ 
\mbox{ for }n\ge 0,\; \alpha\ne \alpha_{1},\dots,\alpha_{r},\\
& &a(x)_{(\alpha_{i},n)}b(x)=c_{i,n}(x)\;\;\;\mbox{ for }0\le n\le k_{i},\\
& &a(x)_{(\alpha_{i},n)}b(x)=0\;\;\;\ \ \ \ \ \ \mbox{ for }n> k_{i}.
\end{eqnarray}
\ep

\begin{proof} From (\ref{eab-commuator-decomposition}) we have
\begin{eqnarray*}
(x_{1}-\alpha_{1}x_{2})^{m}\cdots (x_{1}-\alpha_{r}x_{2})^{m}(a(x_{1})b(x_{2})-K(x_{2},x_{1})=0
\end{eqnarray*}
for any nonnegative integer $m\ge {\rm max}\{k_{1}+1,\dots,k_{r}+1\}$.
By Proposition \ref{pcommutator-formula}, $(a(x),b(x))$ is compatible and
\begin{eqnarray}\label{eab=yalpha}
& &a(x_{1})b(x_{2})-K(x_{2},x_{1})\nonumber\\
&=&\sum_{i=1}^{r}\Res_{x_{0}}
x_{1}^{-1}\delta\left(\frac{\alpha_{i} x_{2}+x_{0}}{x_{1}}\right)
\Y_{\alpha_{i}}(a(x_{2}),x_{0})b(x_{2})\nonumber\\
&=&\sum_{i=1}^{r}\sum_{j=0}^{m}\Res_{x_{0}}a(x_{2})_{(\alpha_{i},j)}b(x_{2})
\frac{1}{j!}\left(\alpha_{i}^{-1}\frac{\partial}{\partial x_{2}}\right)^{j}
x_{1}^{-1}\delta\left(\frac{\alpha_{i} x_{2}}{x_{1}}\right).
\end{eqnarray}
Combining (\ref{eab-commuator-decomposition}) with (\ref{eab=yalpha}), by
Lemma \ref{lpa=0} we have all the rest assertions.
\end{proof}

\br{rpartial-operation} 
{\em Notice that the defined maps
$(a(x),b(x))\mapsto a(x)_{(\alpha,n)}b(x)$ are just partial operations
on $\E(W)$. If $U$ is a pairwise compatible subspace of $\E(W)$ (in the sense that
every (ordered) pair in $U$ is compatible), then
these maps are bilinear maps from $U\otimes U$ to $\E(W)$.
Furthermore, if $U$ contains $1_{W}$ and if $U$ is closed under all these
bilinear operations $(a(x),b(x))\mapsto a(x)_{(\alpha,n)}b(x)$ 
for $\alpha\in \C^{\times},\; n\in \Z$, we would like to know
the algebraic structure that $(U,\{\Y_{\alpha}\},1_{W})$ carries.}  
\er

\section{Algebraic structures on quasi-local subspaces of $\E(W)$}
In this section,  we investigate the algebraic structures
on ``quasi local'' subspaces of $\E(W)$ for any given vector space $W$. 
We show that for any quasi local subset $S$ of $\E(W)$, there exists
a ``closed'' quasi local subspace $U$ which contains $S$ and $1_{W}$
and we furthermore show that $(U,\Y_{1},1_{W})$ carries 
the structure of a vertex algebra. For a general ``closed'' quasi
local subspace $U$
we establish a Jacobi-like identity for the adjoint vertex operators
$\Y_{\alpha}(a(x),x_{0})$ for $\alpha\in \C^{\times}$.
The key and most difficult results are Propositions \ref{pgamma-locality-key} 
and \ref{p-local-adjoint}, 
whose proofs involve rather complicated formal calculus.
The main results are summarized in Theorem \ref{tkey-main}.

As in Section 3, let $W$ be a vector space fixed throughout this section.
The following notion generalizes the notion of mutual locality
([DL], [Li2], cf. [G-K-K]):

\bd{d-quasi-local}
{\em Formal series $a(x),b(x)\in \E(W)$ are said to
be {\em mutually quasi local} if there exists 
a nonzero polynomial $f(x_{1},x_{2})$ such that
\begin{eqnarray}\label{egamma-locality}
f(x_{1},x_{2})a(x_{1})b(x_{2})=f(x_{1},x_{2})b(x_{2})a(x_{1}).
\end{eqnarray}}
\ed

A subset (subspace) $U$ of $\E(W)$ is said to be {\em quasi-local} if
any $a(x), b(x)\in U$ (where $a(x)$ and $b(x)$ may be the same) are
mutually quasi local. Clearly, any quasi-local subset of $\E(W)$
(linearly) spans a quasi local subspace and any quasi-local subset of
$\E(W)$ is contained in some maximal quasi local subspace of $\E(W)$.

We also have the following notion (cf. [G-K-K]):

\bd{d-gammalocal}
{\em  Let $\Gamma$ be a subgroup of $\C^{\times}$.
Formal series $a(x),b(x)\in \E(W)$ are said to be {\em mutually $\Gamma$-local} 
if there exists a (nonzero) polynomial 
\begin{eqnarray}
f(x_{1},x_{2})\in \< (x_{1}-\alpha x_{2})\;|\; \alpha\in \Gamma\>\subset \C[x_{1},x_{2}]
\end{eqnarray}
such that (\ref{egamma-locality}) holds.
The notion of $\Gamma$-local subset(space) of $\E(W)$ is defined in the obvious way.}
\ed

The following is the first key result:

\bp{pgamma-locality-key}
Let $a(x),b(x),c(x)\in \E (W)$. Assume that 
\begin{eqnarray}
&&f(x_{1},x_{2})a(x_{1})b(x_{2})=f(x_{1},x_{2})b(x_{2})a(x_{1}),
\label{efab}\\
&&g(x_{1},x_{2})a(x_{1})c(x_{2})=\tilde{g}(x_{1},x_{2})c(x_{2})a(x_{1}),
\label{egca}\\
&&h(x_{1},x_{2})b(x_{1})c(x_{2})=\tilde{h}(x_{1},x_{2})c(x_{2})b(x_{1}),
\label{ehcb}
\end{eqnarray}
where $f(x,y), g(x,y),\tilde{g}(x,y), h(x,y), \tilde{h}(x,y)$ are nonzero polynomials.
Then for any $\alpha\in \C^{\times}$ and for any $n\in \Z$, 
there exists $k\in \N$ (depending on $n$) such that
\begin{eqnarray}
f(x_{3},\alpha x)^{k}g(x_{3},x)a(x_{3})(b(x)_{(\alpha,n)}c(x))
=f(x_{3},\alpha x)^{k}\tilde{g}(x_{3},x)(b(x)_{(\alpha,n)}c(x))a(x_{3}).
\end{eqnarray}
\ep

\begin{proof} Let $n\in \Z$ be arbitrarily fixed. Let $k$ be a nonnegative integer 
depending on $n$ such that
\begin{eqnarray}\label{exk+n}
x_{0}^{k+n}\iota_{x,x_{0}}(h(\alpha x+x_{0},x)^{-1})\in \C((x))[[x_{0}]].
\end{eqnarray}
(Recall that $\iota_{x,x_{0}}(h(\alpha x+x_{0},x)^{-1})\in \C((x))((x_{0}))$.)
Noticing that $\left(\frac{\partial}{\partial x_{1}}\right)^{i}f(x_{3},x_{1})^{k}$
is a multiple of $f(x_{3},x_{1})$ for $0\le i<k$, from (\ref{efab}) we have
\begin{eqnarray}\label{e4.9proof}
\left(\left(\frac{\partial}{\partial x_{1}}\right)^{i}
f(x_{3},x_{1})^{k}\right)a(x_{3})b(x_{1})
=\left(\left(\frac{\partial}{\partial x_{1}}\right)^{i}
f(x_{3},x_{1})^{k}\right)b(x_{1})a(x_{3})
\;\;\mbox{ for }0\le i<k.
\end{eqnarray}
Using (\ref{exk+n}), (\ref{e4.9proof}) and (\ref{egca})
we obtain
\begin{eqnarray}
& &f(x_{3},\alpha x)^{k}g(x_{3},x)a(x_{3})\left(b(x)_{(\alpha,n)}c(x)\right)\nonumber\\
&=&\Res_{x_{0}}x_{0}^{n}f(x_{3},\alpha x)^{k}g(x_{3},x)a(x_{3})
\left(\Y_{\alpha}(b(x),x_{0})c(x)\right)\nonumber\\
&=&\Res_{x_{0}}\Res_{x_{1}}x_{0}^{n}f(x_{3},\alpha x)^{k}g(x_{3},x)a(x_{3})
\iota_{x,x_{0}}(h(\alpha x+x_{0},x)^{-1})\nonumber\\
& &\ \ \cdot x_{1}^{-1}\delta\left(\frac{\alpha x+x_{0}}{x_{1}}\right)
(h(x_{1},x)b(x_{1})c(x))
\nonumber\\
&=&\Res_{x_{0}}\Res_{x_{1}}x_{0}^{n}f(x_{3},x_{1}-x_{0})^{k}g(x_{3},x)a(x_{3})
\iota_{x,x_{0}}(h(\alpha x+x_{0},x)^{-1})\nonumber\\
& &\ \ \cdot x_{1}^{-1}\delta\left(\frac{\alpha x+x_{0}}{x_{1}}\right)
(h(x_{1},x)b(x_{1})c(x))
\nonumber\\
&=&\Res_{x_{0}}\Res_{x_{1}}x_{0}^{n}
\left(e^{-x_{0}\frac{\partial}{\partial x_{1}}}f(x_{3},x_{1})^{k}\right)g(x_{3},x)a(x_{3})
\iota_{x,x_{0}}(h(\alpha x+x_{0},x)^{-1})\nonumber\\
& &\ \ \cdot x_{1}^{-1}\delta\left(\frac{\alpha x+x_{0}}{x_{1}}\right)
(h(x_{1},x)b(x_{1})c(x))
\nonumber\\
&=&\Res_{x_{0}}\Res_{x_{1}}\sum_{i=0}^{k-1}\frac{(-1)^{i}}{i!}x_{0}^{n+i}
\left(\left(\frac{\partial}{\partial
x_{1}}\right)^{i}f(x_{3},x_{1})^{k}\right)
g(x_{3},x)a(x_{3})
\iota_{x,x_{0}}(h(\alpha x+x_{0},x)^{-1})\nonumber\\
& &\ \ \cdot x_{1}^{-1}\delta\left(\frac{\alpha x+x_{0}}{x_{1}}\right)
(h(x_{1},x)b(x_{1})c(x))
\nonumber\\
&=&\Res_{x_{0}}\Res_{x_{1}}\sum_{i=0}^{k-1}\frac{(-1)^{i}}{i!}x_{0}^{n+i}
\left(\left(\frac{\partial}{\partial x_{1}}\right)^{i}f(x_{3},x_{1})^{k}\right)\tilde{g}(x_{3},x)
\iota_{x,x_{0}}(h(\alpha x+x_{0},x)^{-1})\nonumber\\
& &\ \ \cdot x_{1}^{-1}\delta\left(\frac{\alpha x+x_{0}}{x_{1}}\right)(h(x_{1},x)b(x_{1})c(x))a(x_{3})
\nonumber\\
&=&\Res_{x_{0}}\Res_{x_{1}}x_{0}^{n}
\left(e^{-x_{0}\frac{\partial}{\partial x_{1}}}f(x_{3},x_{1})^{k}\right)\tilde{g}(x_{3},x)
\iota_{x,x_{0}}(h(\alpha x+x_{0},x)^{-1})\nonumber\\
& &\ \ \cdot x_{1}^{-1}\delta\left(\frac{\alpha x+x_{0}}{x_{1}}\right)(h(x_{1},x)b(x_{1})c(x))
a(x_{3})\nonumber\\
&=&\Res_{x_{0}}\Res_{x_{1}}x_{0}^{n}
f(x_{3},x_{1}-x_{0})^{k}\tilde{g}(x_{3},x)
\iota_{x,x_{0}}(h(\alpha x+x_{0},x)^{-1})\nonumber\\
& &\ \ \cdot x_{1}^{-1}\delta\left(\frac{\alpha x+x_{0}}{x_{1}}\right)(h(x_{1},x)b(x_{1})c(x))
a(x_{3})\nonumber\\
&=&\Res_{x_{0}}x_{0}^{n}f(x_{3},\alpha x)^{k}\tilde{g}(x_{3},x)
\left(\Y_{\alpha}(b(x),x_{0})c(x)\right)a(x_{3})\nonumber\\
&=&f(x_{3},\alpha x)^{k}\tilde{g}(x_{3},x)\left(b(x)_{(\alpha,n)}c(x)\right)a(x_{3}).
\end{eqnarray}
This completes the proof.
\end{proof}

Recall Remark \ref{rpartial-operation} that if $U$ is a pairwise compatible
subspace of $\E(W)$, 
we have linear maps $\Y_{\alpha}$ for $\alpha\in \C^{\times}$ from $U$ to
$\Hom (U, \E(W)((x)))$. 
In view of Propositions \ref{pexistence} and \ref{pjacobi-definition},
any quasi local subspace $U$ of $\E(W)$ is pairwise compatible so that
we have these linear maps.

\bd{dclosed-space}
{\em Let $\Gamma$ be a subgroup of $\C^{\times}$.
A pairwise compatible subspace $U$ of $\E(W)$ is said to be {\em $\Y_{\Gamma}$-closed} if
\begin{eqnarray}
a(x)_{(\alpha,n)}b(x)\in U\;\;\;\mbox{ for }a(x),b(x)\in U,\; \alpha\in \Gamma,\;,n\in \Z.
\end{eqnarray}
In the case $\Gamma=\{1\}$, we say $U$ is $\Y_{1}$-closed.}
\ed

Using Proposition \ref{pgamma-locality-key} we have (cf. [LL]):

\bc{clocality-maximal}
Let $W$ be any vector space.
Every maximal quasi local subspace $U$ of $\E(W)$ contains $1_{W}$ and it
is $\Y_{\Gamma}$-closed for any subgroup $\Gamma$ of $\C^{\times}$. 
Furthermore, for any quasi local subset $S$ of $\E(W)$ and
for any subgroup $\Gamma$ of $\C^{\times}$,
there exists a $\Y_{\Gamma}$-closed quasi local subspace of $\E(W)$, which contains $S$ 
and $1_{W}$.
\ec

\begin{proof} Clearly, the linear span of $U$ and $1_{W}$ is quasi local.
With $U$ being maximal we must have $1_{W}\in U$.
Let $a(x),b(x)\in U,\; \alpha\in \Gamma,\; n\in \Z$.
By Proposition \ref{pgamma-locality-key},  for any $c(x)\in U$,
$a(x)_{(\alpha,n)}b(x)$ and $c(x)$ are quasi local.
In particular, $a(x)_{(\alpha,n)}b(x)$ is quasi local with
$a(x)$ and $b(x)$. Using Proposition \ref{pgamma-locality-key} again, we see that
$a(x)_{(\alpha,n)}b(x)$ is quasi local with itself. 
Thus $U+ \C a(x)_{(\alpha,n)}b(x)$ is quasi local.
Since $U$ is maximal, we must have $U+ \C a(x)_{(\alpha,n)}b(x)=U$, which implies that 
$a(x)_{(\alpha,n)}b(x)\in U$. This proves that $U$ is $\Y_{\Gamma}$-closed.

For any quasi local subset $S$ of $\E(W)$, 
the linear span of $S$ is also quasi local. It follows from Zorn's lemma that
there exists a maximal 
quasi local subspace $U$ of $\E(W)$, which contains $S$. 
By the first assertion, $U$ is $\Y_{\Gamma}$-closed and it contains $1_{W}$. 
Thus $U$ is a $\Y_{\Gamma}$-closed quasi local subspace which
contains $S$ and $1_{W}$.
\end{proof}

In view of Corollary \ref{clocality-maximal},
for any quasi local subset $S$ of $\E(W)$ and 
for any subgroup $\Gamma$ of $\C^{\times}$, 
the smallest $\Y_{\Gamma}$-closed
quasi local subspace, containing $S$ and $1_{W}$, of $\E(W)$ exists, which is
the intersection of all $\Y_{\Gamma}$-closed
quasi local subspaces containing $S$ and $1_{W}$.

\bd{dS-generating}
{\em Let $S$ be any quasi local subspace of $\E(W)$.
For any subgroup $\Gamma$ of $\C^{\times}$, we define
$\<S\>_{\Gamma}$ to be the smallest $\Y_{\Gamma}$-closed
quasi local subspace, containing $S$ and $1_{W}$, of $\E(W)$.}
\ed

We also have the following results similar to 
Corollary \ref{clocality-maximal}: 

\bc{cgamma-local-closed}
Let $\Gamma$ be any subgroup of $\C^{\times}$.
Every maximal $\Gamma$-local subspace $U$ of $\E(W)$ contains $1_{W}$ 
and is $\Y_{\Gamma}$-closed. 
For any $\Gamma$-local subset $S$ of $\E(W)$,
$\<S\>_{\Gamma}$ is $\Gamma$-local and it is
the smallest $\Y_{\Gamma}$-closed $\Gamma$-local subspace 
containing $S$ and $1_{W}$ of $\E(W)$.
\ec

\begin{proof} From the proof of Corollary \ref{clocality-maximal}, 
using Proposition \ref{pgamma-locality-key} we immediately see that
every maximal $\Gamma$-local subspace $U$ of $\E(W)$ contains $1_{W}$ 
and is $\Y_{\Gamma}$-closed. For any $\Gamma$-local subset $S$ of
$\E(W)$, by Zorn Lemma $S$ is contained in some
maximal $\Gamma$-local subspace $U$, which has been proved to be
$\Y_{\Gamma}$-closed. Then $\<S\>_{\Gamma}\subset U$.
Consequently, $\<S\>_{\Gamma}$ must be $\Gamma$-local.
The rest is clear.
\end{proof}

The following is the second key result:

\bp{p-local-adjoint}
Let $\Gamma$ be a subgroup of $\C^{\times}$ and let $U$ be a $\Y_{\Gamma}$-closed 
quasi local subspace of $\E (W)$.
For any $a(x),b(x),c(x)\in U$, let $f(x_{1},x_{2})$ be 
a nonzero polynomial such that 
\begin{eqnarray}
f(x_{1},x_{2})a(x_{1})b(x_{2})=f(x_{1},x_{2})b(x_{2})a(x_{1}).
\end{eqnarray}
Then for $\alpha,\beta\in \Gamma$,
\begin{eqnarray}
(x_{1}-\alpha \beta^{-1}x_{2})^{s}
[\Y_{\alpha}(a(x),x_{1}),\Y_{\beta}(b(x),x_{2})]=0,
\end{eqnarray}
where $s$ is the order of the zero of $f(x_{1},x_{2})$ 
at $x_{1}=\alpha \beta^{-1}x_{2}$.
\ep

\begin{proof} Let $g(x,y), h(x,y)$ be nonzero polynomials such that
\begin{eqnarray*}
& &g(x_{1},x_{2})a(x_{1})c(x_{2})=g(x_{1},x_{2})c(x_{2})a(x_{1}),\\
& &h(x_{1},x_{2})b(x_{1})c(x_{2})=h(x_{1},x_{2})c(x_{2})b(x_{1}).
\end{eqnarray*}
Let $\alpha,\beta\in \Gamma$ and let $n\in \Z$ be arbitrarily fixed. 
Since $b(x)_{(\beta,m)}c(x)=0$ for $m$ sufficiently large
(by Proposition \ref{pexistence}), from Proposition \ref{pgamma-locality-key}, 
there exists a nonnegative integer $k$ depending on $n$ such that
for {\em all} $m\ge n$,
\begin{eqnarray}
f(x_{1},\beta x)^{k}g(x_{1},x)a(x_{1})(b(x)_{(\beta,m)}c(x))
=f(x_{1},\beta x)^{k}g(x_{1},x)(b(x)_{(\beta,m)}c(x))a(x_{1}).
\end{eqnarray}
Let $w\in W$ be arbitrarily fixed and let $l$ be a
nonnegative integer such that $x^{l}a(x)w\in W[[x]]$. By Proposition \ref{pexistence} we have
\begin{eqnarray}
& &f(x_{1}+\alpha x,\beta x)^{k}g(x_{1}+\alpha x,x)(x_{1}+\alpha x)^{l}
\Y_{\alpha}(a(x),x_{1})(b(x)_{(\beta,m)}c(x))w\nonumber\\
&=&f(x_{1}+\alpha x,\beta x)^{k}g(x_{1}+\alpha x,x)(x_{1}+\alpha x)^{l}
a(x_{1}+\alpha x)(b(x)_{(\beta,m)}c(x))w
\end{eqnarray}
for {\em all} $m\ge n$.
That is, for {\em all} $i\ge 0$,
\begin{eqnarray}
& &\Res_{x_{2}}x_{2}^{n+i}f(x_{1}+\alpha x,\beta x)^{k}
g(x_{1}+\alpha x,x)(x_{1}+\alpha x)^{l}
\left(\Y_{\alpha}(a(x),x_{1})\Y_{\beta}(b(x),x_{2})c(x)\right)w
\nonumber\\
&=&\Res_{x_{2}}x_{2}^{n+i}f(x_{1}+\alpha x,\beta x)^{k}
g(x_{1}+\alpha x,x)(x_{1}+\alpha x)^{l}
a(x_{1}+\alpha x)\left(\Y_{\beta}(b(x),x_{2})c(x)\right)w.\nonumber\\
& &\mbox{}
\end{eqnarray}
Since $f(x_{1}+\alpha x,x_{2}+\beta x)$ is a polynomial in $x_{1},x_{2},x$,  we immediately
have that for {\em all} $j\ge 0$,
\begin{eqnarray}\label{ekey-stone}
& &\Res_{x_{2}}x_{2}^{n+j}f(x_{1}+\alpha x,\beta x)^{k}
g(x_{1}+\alpha x,x)(x_{1}+\alpha x)^{l}f(x_{1}+\alpha x,x_{2}+\beta x)\nonumber\\
& &\ \ \cdot \left(\Y_{\alpha}(a(x),x_{1})\Y_{\beta}(b(x),x_{2})c(x)\right)w
\nonumber\\
&=&\Res_{x_{2}}x_{2}^{n+j}f(x_{1}+\alpha x,\beta x)^{k}
g(x_{1}+\alpha x,x)(x_{1}+\alpha x)^{l}f(x_{1}+\alpha x,x_{2}+\beta x)\nonumber\\
& &\ \ \cdot a(x_{1}+\alpha x)\left(\Y_{\beta}(b(x),x_{2})c(x)\right)w.
\end{eqnarray}
Now, let $r$ be a nonnegative integer such that $x^{r}b(x)w\in W[[x]]$. 
By Proposition \ref{pexistence} we have
$$h(x_{2}+\beta x,x)(x_{2}+\beta x)^{r}
\left(\Y_{\beta}(b(x),x_{2})c(x)\right)w
=h(x_{2}+\beta x,x)(x_{2}+\beta x)^{r}b(x_{2}+\beta x)c(x)w.$$
Using this, from (\ref{ekey-stone}) we obtain
\begin{eqnarray}\label{eto-be-cancel}
& &\Res_{x_{2}}x_{2}^{n}f(x_{1}+\alpha x,\beta x)^{k}
g(x_{1}+\alpha x,x)(x_{1}+\alpha x)^{l}h(x_{2}+\beta x,x)(x_{2}+\beta x)^{r}
\nonumber\\
& &\ \ \cdot f(x_{1}+\alpha x,x_{2}+\beta x)
\left(\Y_{\alpha}(a(x),x_{1})\Y_{\beta}(b(x),x_{2})c(x)\right)w
\nonumber\\
&=&\Res_{x_{2}}x_{2}^{n}f(x_{1}+\alpha x,\beta x)^{k}g(x_{1}+\alpha x,x)
(x_{1}+\alpha x)^{l}h(x_{2}+\beta x,x)(x_{2}+\beta x)^{r}\nonumber\\
& &\ \ \cdot f(x_{1}+\alpha x,x_{2}+\beta x)
a(x_{1}+\alpha x)b(x_{2}+\beta x)c(x)w.
\end{eqnarray}
Noticing that
$$g(x_{1}+\alpha x,x)(x_{1}+\alpha x)^{l}a(x_{1}+\alpha x)c(x)w\in W((x))[[x_{1}]],$$
we have
\begin{eqnarray*}
& &g(x_{1}+\alpha x,x)(x_{1}+\alpha x)^{l}f(x_{1}+\alpha x,x_{2}+\beta x)
a(x_{1}+\alpha x)b(x_{2}+\beta x)c(x)w\nonumber\\
&=&g(x_{1}+\alpha x,x)(x_{1}+\alpha x)^{l}f(x_{1}+\alpha x,x_{2}+\beta x)
b(x_{2}+\beta x)a(x_{1}+\alpha x)c(x)w,
\end{eqnarray*}
which lies in $W((x))((x_{1}))[[x_{2},x_{2}^{-1}]]$,
a vector space over $\C((x))((x_{1}))$. The expression on the left hand side
of (\ref{eto-be-cancel}) also lies in $W((x))((x_{1}))[[x_{2},x_{2}^{-1}]]$.
In view of Lemma \ref{lcancellation} 
we can cancel the factor $f(x_{1}+\alpha x,\beta x)^{k}$, to obtain
\begin{eqnarray}
& &\Res_{x_{2}}x_{2}^{n}
g(x_{1}+\alpha x,x)(x_{1}+\alpha x)^{l}h(x_{2}+\beta x,x)(x_{2}+\beta x)^{r}
\nonumber\\
& &\ \ \cdot f(x_{1}+\alpha x,x_{2}+\beta x)
\left(\Y_{\alpha}(a(x),x_{1})\Y_{\beta}(b(x),x_{2})c(x)\right)w
\nonumber\\
&=&\Res_{x_{2}}x_{2}^{n}g(x_{1}+\alpha x,x)
(x_{1}+\alpha x)^{l}h(x_{2}+\beta x,x)(x_{2}+\beta x)^{r}\nonumber\\
& &\ \ \cdot f(x_{1}+\alpha x,x_{2}+\beta x)
a(x_{1}+\alpha x)b(x_{2}+\beta x)c(x)w.
\end{eqnarray}
Notice that only $k$ here depends on $n$ (all the other components
do not depend on $n$) and $n$ is arbitrarily fixed. 
Then we can drop $\Res_{x_{2}}x_{2}^{n}$ to get
\begin{eqnarray}\label{eyaybexp}
& &
g(x_{1}+\alpha x,x)(x_{1}+\alpha x)^{l}h(x_{2}+\beta x,x)(x_{2}+\beta x)^{r}
\nonumber\\
& &\ \ \cdot f(x_{1}+\alpha x,x_{2}+\beta x)
\left(\Y_{\alpha}(a(x),x_{1})\Y_{\beta}(b(x),x_{2})c(x)\right)w
\nonumber\\
&=&g(x_{1}+\alpha x,x)
(x_{1}+\alpha x)^{l}h(x_{2}+\beta x,x)(x_{2}+\beta x)^{r}\nonumber\\
& &\ \ \cdot f(x_{1}+\alpha x,x_{2}+\beta x)
a(x_{1}+\alpha x)b(x_{2}+\beta x)c(x)w.
\end{eqnarray}
Symmetrically, there are nonnegative integers $l',r'$ such that
\begin{eqnarray}\label{eybyaexp}
& &
g(x_{1}+\alpha x,x)(x_{1}+\alpha x)^{l'}h(x_{2}+\beta x,x)(x_{2}+\beta x)^{r'}
\nonumber\\
& &\ \ \cdot f(x_{1}+\alpha x,x_{2}+\beta x)
\left(\Y_{\beta}(b(x),x_{2})\Y_{\alpha}(a(x),x_{1})c(x)\right)w
\nonumber\\
&=&g(x_{1}+\alpha x,x)
(x_{1}+\alpha x)^{l}h(x_{2}+\beta x,x)(x_{2}+\beta x)^{r}\nonumber\\
& &\ \ \cdot f(x_{1}+\alpha x,x_{2}+\beta x)
b(x_{2}+\beta x)a(x_{1}+\alpha x)c(x)w.
\end{eqnarray}
Since
$$f(x_{1}+\alpha x,x_{2}+\beta x)a(x_{1}+\alpha x)b(x_{2}+\beta x)=
f(x_{1}+\alpha x,x_{2}+\beta x)b(x_{2}+\beta x)a(x_{1}+\alpha x),$$
combining (\ref{eyaybexp}) with (\ref{eybyaexp}), 
then setting $l''=l+l'$ and $r''=r+r$ we get
\begin{eqnarray}
& &g(x_{1}+\alpha x,x)(x_{1}+\alpha x)^{l''}h(x_{2}+\beta x,x)(x_{2}+\beta x)^{r''}
\nonumber\\
& &\ \ \cdot f(x_{1}+\alpha x,x_{2}+\beta x)
\left(\Y_{\alpha}(a(x),x_{1})\Y_{\beta}(b(x),x_{2})c(x)\right)w
\nonumber\\
&=&g(x_{1}+\alpha x,x)(x_{1}+\alpha x)^{l''}h(x_{2}+\beta x,x)(x_{2}+\beta x)^{r''}
\nonumber\\
& &\ \ \cdot f(x_{1}+\alpha x,x_{2}+\beta x)
\left(\Y_{\beta}(b(x),x_{2})\Y_{\alpha}(a(x),x_{1})c(x)\right)w.
\end{eqnarray}
Noting that 
\begin{eqnarray*}
& &\left(\Y_{\alpha}(a(x),x_{1})\Y_{\beta}(b(x),x_{2})c(x)\right)w
\in W((x))((x_{1}))((x_{2}))
\subset W((x))((x_{1}))[[x_{2},x_{2}^{-1}]]\\
& &\left(\Y_{\beta}(b(x),x_{2})\Y_{\alpha}(a(x),x_{1})c(x)\right)w
\in W((x))((x_{2}))((x_{1}))\subset W((x))((x_{1}))[[x_{2},x_{2}^{-1}]]
\end{eqnarray*}
and that $W((x))((x_{1}))[[x_{2},x_{2}^{-1}]]$ is a vector space over
$\C((x))((x_{1}))$, in view of Lemma \ref{lcancellation} 
we can cancel the factor $g(x_{1}+\alpha x,x)(x_{1}+\alpha x)^{l''}$ 
and similarly we can cancel the factor
$h(x_{2}+\beta x,x)(x_{2}+\beta x)^{r''}$.
(But we cannot cancel the factor $f(x_{1}+\alpha x,x_{2}+\beta x)$.)
Then we get
\begin{eqnarray}
& &f(x_{1}+\alpha x,x_{2}+\beta x)
\left(\Y_{\alpha}(a(x),x_{1})\Y_{\beta}(b(x),x_{2})c(x)\right)w\nonumber\\
&=&f(x_{1}+\alpha x,x_{2}+\beta x)
\left(\Y_{\beta}(b(x),x_{2})\Y_{\alpha}(a(x),x_{1})c(x)\right)w.
\end{eqnarray}
Noticing that $f$ does not depend on $w$, we have 
\begin{eqnarray}\label{efabc}
& &f(x_{1}+\alpha x,x_{2}+\beta x)
\Y_{\alpha}(a(x),x_{1})\Y_{\beta}(b(x),x_{2})c(x)\nonumber\\
&=&f(x_{1}+\alpha x,x_{2}+\beta x)
\Y_{\beta}(b(x),x_{2})\Y_{\alpha}(a(x),x_{1})c(x).
\end{eqnarray}

Write
$$f(x_{1},x_{2})=(x_{1}-\alpha\beta^{-1}x_{2})^{s}f_{0}(x_{1},x_{2})),$$
where $s$ is a nonnegative integer and $f_{0}(x_{1},x_{2})\in
\C[x_{1},x_{2}]$
with $f_{0}(\alpha\beta^{-1}x_{2}, x_{2})\ne 0$.
We have
$$f(x_{1}+\alpha x,x_{2}+\beta x)
=(x_{1}-\alpha\beta^{-1}x_{2})^{s}f_{0}(x_{1}+\alpha x,x_{2}+\beta x)$$
with $f_{0}(\alpha x,\beta x)\ne 0$.
In view of Lemma \ref{lcancellation} we can cancel
the factor $f_{0}(x_{1}+\alpha x,x_{2}+\beta x)$ from 
(\ref{efabc}). After cancelation we get
\begin{eqnarray}
& &(x_{1}-\alpha \beta^{-1}x_{2})^{s}
\Y_{\alpha}(a(x),x_{1})\Y_{\beta}(b(x),x_{2})c(x)
\nonumber\\
&=&(x_{1}-\alpha \beta^{-1}x_{2})^{s}
\Y_{\beta}(b(x),x_{2})\Y_{\alpha}(a(x),x_{1})c(x).
\end{eqnarray}
This completes the proof. 
\end{proof}

Now we are in a position to present our main and key result:

\bt{tkey-main}
Let $W$ be a vector space, $\Gamma$ a subgroup of $\C^{\times}$ and 
$U$ a $\Y_{\Gamma}$-closed quasi local subspace of $\E(W)$, containing $1_{W}$,
in the sense that
\begin{eqnarray}
a(x)_{(\alpha,n)}b(x)\in U\;\;\;\mbox{ for }a(x),b(x)\in U,\; 
\alpha\in \Gamma,\; n\in \Z.
\end{eqnarray}
Then 
\begin{eqnarray}
& &\Y_{\alpha}(1_{W},x_{0})a(x)=a(x),\\
& &\Y_{\alpha}(a(x),x_{0})1_{W}=e^{\alpha^{-1}x_{0}D}R_{\alpha}a(x)
\;\;\;\mbox{ for }a(x)\in U,
\end{eqnarray}
and the following identity holds for any 
$a(x),b(x)\in U,\; \alpha,\beta\in \Gamma$:
\begin{eqnarray}\label{ejacobi-end}
& &x_{0}^{-1}\delta\left(\frac{x_{1}-\alpha\beta^{-1} x_{2}}{x_{0}}\right)
\Y_{\alpha}(a(x),x_{1})\Y_{\beta}(b(x),x_{2})\nonumber\\
& &-x_{0}^{-1}\delta\left(\frac{-\alpha\beta^{-1} x_{2}+
x_{1}}{x_{0}}\right)
\Y_{\beta}(b(x),x_{2})\Y_{\alpha}(a(x),x_{1})\nonumber\\
&=&x_{1}^{-1}\delta\left(\frac{\alpha\beta^{-1}x_{2}+x_{0}}{x_{1}}\right)
\Y_{\beta}(\Y_{\alpha\beta^{-1}}(a(x),x_{0})b(x),x_{2}).
\end{eqnarray}
In particular, $(U,\Y_{1},1_{W})$ carries the structure of a vertex algebra.
Furthermore, for any quasi local subset $S$ of $\E(W)$,
$(\<S\>_{\Gamma},\Y_{1},1_{W})$ carries the structure of a vertex algebra,
where $\<S\>_{\Gamma}$ is the smallest $\Y_{\Gamma}$-closed 
quasi local subspace containing $S$ and $1_{W}$.
\et

\begin{proof} It is clear that all the assertions follow
if the Jacobi-like identity (\ref{ejacobi-end}) is proved.
With Propositions \ref{pvacuum} and \ref{pbasic-property2} and \ref{p-local-adjoint}, 
it follows from \cite{li-local} ([LL], cf. [DL], [FHL]) that $(U,\Y_{1},1_{W})$ carries 
the structure of a vertex algebra. Thus (\ref{ejacobi-end}) holds for $\alpha=\beta=1$.
The general case follows from this special case and 
the relation (\ref{e-conjugation-calculus-3}):
For $u,v\in U,\; \alpha,\beta\in \Gamma$, we have
\begin{eqnarray}
& &x_{0}^{-1}\delta\left(\frac{\alpha^{-1}x_{1}-\beta^{-1}x_{2}}
{x_{0}}\right)
\Y_{1}(R_{\alpha}u,\alpha^{-1}x_{1})\Y_{1}(R_{\beta}v,\beta^{-1}x_{2})\nonumber\\
& &-x_{0}^{-1}\delta\left(
\frac{-\beta^{-1}x_{2}+\alpha^{-1}x_{1}}{x_{0}}\right)
\Y_{1}(R_{\beta}v,\beta^{-1}x_{2})\Y_{1}(R_{\alpha}u,\alpha^{-1}x_{1})\nonumber\\
&=&\alpha x_{1}^{-1}\delta\left(
\frac{\beta^{-1}x_{2}+x_{0}}{\alpha^{-1}x_{1}}\right)
\Y_{1}(\Y_{1}(R_{\alpha}u,x_{0})R_{\beta}v,\beta^{-1}x_{2}).
\end{eqnarray}
By rewriting the three delta functions we get
\begin{eqnarray}
& &(\alpha x_{0})^{-1}\delta\left(\frac{x_{1}-\alpha\beta^{-1}x_{2}}
{\alpha x_{0}}\right)
\Y_{1}(R_{\alpha}u,\alpha^{-1}x_{1})\Y_{1}(R_{\beta}v,\beta^{-1}x_{2})\nonumber\\
& &-(\alpha x_{0})^{-1}\delta\left(
\frac{-\alpha\beta^{-1}x_{2}+x_{1}}{\alpha x_{0}}\right)
\Y_{1}(R_{\beta}v,\beta^{-1}x_{2})\Y_{1}(R_{\alpha}u,\alpha^{-1}x_{1})\nonumber\\
&=&x_{1}^{-1}\delta\left(
\frac{\alpha\beta^{-1}x_{2}+\alpha x_{0}}{x_{1}}\right)
\Y_{1}(\Y_{1}(R_{\alpha}u,x_{0})R_{\beta}v,\beta^{-1}x_{2}).
\end{eqnarray}
In view of (\ref{e-conjugation-calculus-3}), we have
$$\Y_{1}(\Y_{1}(R_{\alpha}u,x_{0})R_{\beta}v,\beta^{-1}x_{2})=
\Y_{\beta}(\Y_{\beta^{-1}\alpha}(u,\alpha x_{0})v,x_{2}).$$
Replacing $x_{0}\rightarrow \alpha^{-1}x_{0}$ we get
the desired Jacobi identity
\end{proof}

Similar to the situation for vertex subalgebras 
generated by subsets (cf. [LL]) we have:

\bp{pgenerating-property}
Let $S$ be a quasi local subset of $\E(W)$. Then
for any subgroup $\Gamma$ of $\C^{\times}$,
$\<S\>_{\Gamma}$ is linearly spanned by the vectors
\begin{eqnarray}\label{egenerating-iterate}
a^{(1)}(x)_{(\alpha_{1},n_{1})}\cdots a^{(r)}(x)_{(\alpha_{r},n_{r})}1_{W}
\end{eqnarray}
for $r\ge 0,\; a^{(i)}(x)\in S,\;\alpha_{i}\in \Gamma,\; n_{i}\in \Z$.
\ep

\begin{proof} Let $J$ be the subspace of $\E(W)$, linearly spanned by
those vectors in (\ref{egenerating-iterate}). It follows from the Jacobi-like identity
(\ref{ejacobi-end}) and induction that $J$ is $\Y_{\Gamma}$-closed. Clearly,
$J$ contains $S$ and $1_{W}$. Thus, $\<S\>_{\Gamma}\subset J$.
On the other hand, any $\Y_{\Gamma}$-closed quasi local subspace 
containing $S$ and $1_{W}$ 
must contain $J$, so that $J\subset \<S\>_{\Gamma}$.
This proves $\<S\>_{\Gamma}=J$, completing the proof.
\end{proof}

\bd{dva}
{\em Let $S$ be a quasi local subset of $\E(W)$.
We set 
\begin{eqnarray}
\<S\>=\<S\>_{\{1\}},
\end{eqnarray}
the smallest $\Y_{1}$-closed quasi local subspace, 
containing $S$ and $1_{W}$, of $\E(W)$.
We call $\<S\>$ the {\em vertex algebra generated by $S$}.}
\ed

{}From definition we have $\<S\>_{\Gamma}\subset \<S\>_{\Gamma'}$
if $\Gamma\subset \Gamma'$ are subgroups of $\C^{\times}$.
Then $\<S\>$ is a vertex subalgebra of 
$\<S\>_{\Gamma}$ for any subgroup $\Gamma$ of $\C^{\times}$
and it is the vertex subalgebra of 
$\<S\>_{\Gamma}$ generated by $S$.

\bl{lfor-application}
Let $W$ be a vector space, $\Gamma$ a subgroup of $\C^{\times}$ and
$U$ a $\Gamma$-local subspace of $\E(W)$. If $U$ is also a
$\Gamma$-submodule of $\E(W)$, then $\<U\>=\<U\>_{\Gamma}$. 
\el

\begin{proof} With $\<U\>=\<U\>_{\{1\}}\subset\<U\>_{\Gamma}$ we must prove 
$\<U\>_{\Gamma}\subset \<U\>$.
Recall from (\ref{e-conjugation-calculus}) (Proposition \ref{pbasic-property2}) that
\begin{eqnarray*}
R_{\alpha}\Y_{1}(a(x),x_{0})=\Y_{1}(R_{\alpha}a(x),\alpha^{-1}x_{0})R_{\alpha^{-1}}.
\end{eqnarray*}
Since $U$ is a $\Gamma$-submodule of $\E(W)$, it follows from 
Proposition \ref{pgenerating-property} (and induction) that 
$\<U\>$ is a $\Gamma$-submodule of $\E(W)$.
{}From (\ref{e-conjugation-calculus-2}) 
(Proposition \ref{pbasic-property2}) we have
\begin{eqnarray*}
\Y_{\alpha}(a(x),x_{0})=\Y_{1}(R_{\alpha}a(x),\alpha^{-1}x_{0})
\end{eqnarray*}
It follows that $\<U\>$ is $\Y_{\Gamma}$-closed.
Then we have $\<U\>_{\Gamma}\subset \<U\>$, completing the proof.
\end{proof}

\section{Vertex algebras and their quasi modules}
In this section we formulate and study 
the notion of quasi module for a vertex algebra.
We strengthen some of the assertions of Theorem \ref{tkey-main}, 
showing that for any vector space $W$ and 
for any $\Y_{\Gamma}$-closed quasi local
subspace $U$ of $\E(W)$ containing $1_{W}$, $W$ is a
quasi module for the vertex algebra $(U,\Y_{1},1_{W})$.
We also prove certain results analogous to those for modules.

First we study the possibility of a certain natural generalization of
the notion of vertex algebra. It has been well known that
the locality (namely weak commutativity) property plays 
a very important role in the theory of vertex algebras.
With the notion of quasi-locality, one naturally considers to generalize the
notion of vertex algebra by replacing the usual locality with quasi
locality. The following proposition says that quasi
locality does not give rise to something more general than vertex
algebras.

\bp{pno-other}
Let $V$ be a vector space equipped with a linear map
$Y$ from $V$ to $\Hom (V,V((x)))$, a distinguished vector ${\bf 1}$
and a linear operator $\D$ on $V$ such that all the following
conditions hold: 
\begin{eqnarray}
& &\D ({\bf 1})=0,\label{ed1=0}\\
& &Y({\bf 1},x)=1,\\
& &Y(v,x){\bf 1}\in V[[x]]\;\;\mbox{ and }\;\lim_{x\rightarrow
0}Y(v,x){\bf 1}=v,\label{eno-other-creation}\\
& &[\D, Y(v,x)]={d\over dx}Y(v,x)\label{eno-other-dbracket}
\end{eqnarray}
for $v\in V$ and such that for $u,v\in V$, there exists a nonzero polynomial
$f(x_{1},x_{2})$ such that
\begin{eqnarray}\label{epuv=pvu}
f(x_{1},x_{2})Y(u,x_{1})Y(v,x_{2})=f(x_{1},x_{2})Y(v,x_{2})Y(u,x_{1}).
\end{eqnarray}
Then $(V,Y,{\bf 1})$ is an (ordinary) vertex algebra.
\ep

\begin{proof} First, (\ref{eno-other-dbracket}) together with 
the Taylor theorem immediately
gives the following conjugation formula
\begin{eqnarray}
e^{x\D}Y(v,x_{0})e^{-x\D}=Y(v,x_{0}+x)\;\;\;\mbox{ for }v\in V.
\end{eqnarray}
Applying this to ${\bf 1}$, using (\ref{ed1=0}) and (\ref{eno-other-creation}) 
(by setting $x_{0}=0$) we get
$$Y(v,x){\bf 1}=e^{x\D}v\;\; \;\;\;\mbox{ for }v\in V.$$
Let $u,v,w\in V$ and let $f(x_{1},x_{2})$ be a nonzero polynomial
such that (\ref{epuv=pvu}) holds.
We have
\begin{eqnarray}
Y(Y(u,x_{1})Y(v,x_{2})w,x){\bf 1}
&=&e^{x\D}Y(u,x_{1})Y(v,x_{2})w\nonumber\\
&=&Y(u,x_{1}+x)Y(v,x_{2}+x)e^{x\D}w,\\
Y(Y(v,x_{2})Y(u,x_{1})w,x){\bf 1}
&=&Y(v,x_{2}+x)Y(u,x_{1}+x)e^{x\D}w.
\end{eqnarray}
Using these relations and (\ref{epuv=pvu}) we get
\begin{eqnarray}
& &f(x_{1}+x,x_{2}+x)Y(Y(u,x_{1})Y(v,x_{2})w,x){\bf 1}\nonumber\\
&=&f(x_{1}+x,x_{2}+x)Y(Y(v,x_{2})Y(u,x_{1})w,x){\bf 1}.
\end{eqnarray}
Factor $f(x_{1},x_{2})$ as 
$$f(x_{1},x_{2})=(x_{1}-x_{2})^{k}f_{0}(x_{1},x_{2}),$$
where $k\ge 0,\; f_{0}(x_{1},x_{2})\in \C[x_{1},x_{2}]$ with
$f_{0}(x_{1},x_{1})\ne 0$. Then we have
\begin{eqnarray}\label{eproof-precomm}
& &(x_{1}-x_{2})^{k}f_{0}(x_{1}+x,x_{2}+x)
Y(Y(u,x_{1})Y(v,x_{2})w,x){\bf 1}\nonumber\\
&=&(x_{1}-x_{2})^{k}f_{0}(x_{1}+x,x_{2}+x)Y(Y(v,x_{2})Y(u,x_{1})w,x){\bf 1}.
\end{eqnarray}
Notice that
\begin{eqnarray*}
& &Y(Y(u,x_{1})Y(v,x_{2})w,x){\bf 1}\in V[[x]]((x_{1}))((x_{2})),\\
& &Y(Y(v,x_{2})Y(u,x_{1})w,x){\bf 1}\in V[[x]]((x_{2}))((x_{1})).
\end{eqnarray*}
In view of Lemma \ref{lcancellation}, from (\ref{eproof-precomm})
we get
\begin{eqnarray}
(x_{1}-x_{2})^{k}Y(Y(u,x_{1})Y(v,x_{2})w,x){\bf 1}
=(x_{1}-x_{2})^{k}Y(Y(v,x_{2})Y(u,x_{1})w,x){\bf 1}.
\end{eqnarray}
In view of (\ref{eno-other-creation}), setting $x=0$ we obtain
\begin{eqnarray}
(x_{1}-x_{2})^{k}Y(u,x_{1})Y(v,x_{2})w
=(x_{1}-x_{2})^{k}Y(v,x_{2})Y(u,x_{1})w.
\end{eqnarray}
This proves that the usual locality (weak commutativity) holds. 
Then it follows from \cite{li-local} ([LL], cf. [DL], [FHL]) 
that $V$ is a vertex algebra.
\end{proof}

Even though the notion of quasi locality does not lead to
a generalization of the notion of vertex algebra,
it does lead to a new notion of what we call quasi module
for an (ordinary) vertex algebra.

\bd{dquasi-module-va}
{\em Let $V$ be a vertex algebra. A {\em quasi $V$-module} 
is a vector space $W$ equipped with
a linear map $Y_{W}$ from $V$ to $\Hom (W,W((x)))$ such that
\begin{eqnarray}
Y_{W}({\bf 1},x)=1_{W}
\end{eqnarray}
and such that for any $u,v\in V$, there exists a nonzero polynomial
$f(x_{1},x_{2})$ such that
\begin{eqnarray}\label{epjacobi}
& &x_{0}^{-1}\delta\left(\frac{x_{1}-x_{2}}{x_{0}}\right)
f(x_{1},x_{2})Y_{W}(u,x_{1})Y_{W}(v,x_{2})\nonumber\\
& &-x_{0}^{-1}\delta\left(\frac{x_{2}-x_{1}}{-x_{0}}\right)f(x_{1},x_{2})
Y_{W}(v,x_{2})Y_{W}(u,x_{1})\nonumber\\
&=&x_{2}^{-1}\delta\left(\frac{x_{1}-x_{0}}{x_{2}}\right)
f(x_{1}, x_{2})Y_{W}(Y(u,x_{0})v,x_{2}).
\end{eqnarray}}
\ed

For convenience we refer the identity (\ref{epjacobi}) as
the {\em quasi Jacobi identity for $(u,v)$}.

\bl{lquasi-module-dproperty}
Let $(W,Y_{W})$ be a quasi module for a vertex algebra $V$.
Then
\begin{eqnarray}\label{edproperty-quasi}
Y_{W}(\D v,x)=\frac{d}{dx}Y_{W}(v,x)\;\;\;\mbox{ for }v\in V.
\end{eqnarray}
\el

\begin{proof} It is similar to the proof for usual modules
(cf. [LL]). Let $v\in V$. There exists 
$0\ne f(x_{1},x_{2})\in \C[x_{1},x_{2}]$ such that the quasi Jacobi 
identity for $(v,{\bf 1})$ holds.
Taking $\Res_{x_{1}}$ and using $Y_{W}({\bf 1},x)=1_{W}$, we get
\begin{eqnarray}
& &f(x_{2}+x_{0},x_{2})Y_{W}(Y(v,x_{0}){\bf 1},x_{2})\nonumber\\
&=&\Res_{x_{1}}\left(
x_{0}^{-1}\delta\left(\frac{x_{1}-x_{2}}{x_{0}}\right)
f(x_{1},x_{2})Y_{W}(v,x_{1})
-x_{0}^{-1}\delta\left(\frac{x_{2}-x_{1}}{-x_{0}}\right)
f(x_{1},x_{2})Y_{W}(v,x_{1})\right)
\nonumber\\
&=&\Res_{x_{1}}
x_{1}^{-1}\delta\left(\frac{x_{2}+x_{0}}{x_{1}}\right)
f(x_{1},x_{2})Y_{W}(v,x_{1})
\nonumber\\
&=&f(x_{2}+x_{0},x_{2})Y_{W}(v,x_{2}+x_{0}).
\end{eqnarray}
Noticing that both $Y_{W}(Y(v,x_{0}){\bf 1},x_{2})$ and $Y_{W}(v,x_{2}+x_{0})$
involve only nonnegative powers of $x_{0}$, we can apply the
cancelation law (Lemma \ref{lcancellation}), to get
$$Y_{W}(Y(v,x_{0}){\bf
1},x_{2})=Y_{W}(v,x_{2}+x_{0})
=e^{x_{0}\frac{\partial}{\partial x_{2}}}Y_{W}(v,x_{2}).$$
{}From this we immediately have (\ref{edproperty-quasi}).
\end{proof}

With the notion of quasi module we have the following important result:

\bt{tvertex-algebra-case}
Let $W$ be any vector space, let $\Gamma$ be a subgroup of $\C^{\times}$
and let $U$ be a $\Y_{\Gamma}$-closed quasi local subspace of $\E(W)$, 
containing $1_{W}$. Then
$(W,Y_{W})$ carries the structure of a quasi module for the vertex
algebra $(U,\Y_{1},1_{W})$ obtained in Theorem \ref{tkey-main}, where
\begin{eqnarray}
Y_{W}(a(x),x_{0})=a(x_{0})\;\;\;\mbox{ for }a(x)\in U. 
\end{eqnarray}
Furthermore, 
for any quasi local subset $S$ of $\E(W)$, $(W,Y_{W})$ 
carries the structure of a quasi module for the vertex algebra 
$\<S\>_{\Gamma}$.
\et

\begin{proof} Clearly, we only need to prove the first assertion.
For $a(x),b(x)\in U$, let $f(x_{1},x_{2})$ be 
a nonzero polynomial such that
\begin{eqnarray}\label{efab=fba-proof}
f(x_{1},x_{2})a(x_{1})b(x_{2})=f(x_{1},x_{2})b(x_{2})a(x_{1}).
\end{eqnarray}
We have
\begin{eqnarray}
f(x_{1},x_{2})Y_{W}(a(x),x_{1})Y_{W}(b(x),x_{2})
&=&f(x_{1},x_{2})a(x_{1})b(x_{2})\nonumber\\
&=&f(x_{1},x_{2})b(x_{2})a(x_{1})\nonumber\\
&=&f(x_{1},x_{2})Y_{W}(b(x),x_{2})Y_{W}(a(x),x_{1}).
\end{eqnarray}
{}From (\ref{efab=fba-proof}) we have
that $f(x_{1},x_{2})a(x_{1})b(x_{2})\in \Hom (W,W((x_{1},x_{2})))$.
Let $w\in W$ and let $l$ be a nonnegative integer such that
$x^{l}a(x)w\in W[[x]]$, so that
$$x_{1}^{l}f(x_{1},x_{2})a(x_{1})b(x_{2})w\in W[[x_{1},x_{2}]][x_{2}^{-1}].$$
With
\begin{eqnarray*}
& &Y_{W}(a(x),x_{0}+x_{2})Y_{W}(b(x),x_{2})w
=a(x_{0}+x_{2})b(x_{2})w,\\
& &Y_{W}(\Y_{1}(a(x),x_{0})b(x),x_{2})w
=\left(\Y_{1}(a(x),x_{0})b(x)\right)w|_{x=x_{2}},
\end{eqnarray*}
{}from Proposition \ref{pexistence}  we have
\begin{eqnarray}
& &(x_{0}+x_{2})^{l}f(x_{0}+x_{2},x_{2})
Y_{W}(a(x),x_{0}+x_{2})
Y_{W}(b(x),x_{2})w\nonumber\\
&=&(x_{0}+x_{2})^{l}f(x_{0}+x_{2},x_{2})
Y_{W}(\Y_{1}(a(x),x_{0})b(x),x_{2})w.
\end{eqnarray}
In view of Lemma \ref{lformaljacobiidentity}, 
the required quasi Jacobi identity for $(u,v)$ holds. By definition $Y_{W}(1_{W},x_{0})=1_{W}$. 
This proves that $W$ is a quasi module.
\end{proof}

The following result says that for a vertex algebra $V$,
the notion of quasi $V$-module is canonically
equivalent to a notion of ``quasi representation'':

\bt{trepresentation}
Let $V$ be a vertex algebra and let $W$ be a vector space
equipped with a linear map $Y_{W}$ from $V$ to $\E(W)$.
Then $(W,Y_{W})$ carries the structure of a quasi $V$-module
if and only if $\{ Y_{W}(v,x)\;|\; v\in V \}$ is a $\Y_{1}$-closed quasi local
subspace containing $1_{W}$ of $\E(W)$ and 
$Y_{W}$ is a vertex algebra homomorphism.
\et

\begin{proof} Assume that $(W,Y_{W})$ is a quasi $V$-module. 
Set
\begin{eqnarray*}
U=\{ Y_{W}(v,x)\;|\; v\in V\}\subset \E(W).
\end{eqnarray*}
For any $u,v\in V$, there exists 
$0\ne f(x_{1},x_{2})\in \C[x_{1},x_{2}]$ such that
\begin{eqnarray}
& &x_{0}^{-1}\delta\left(\frac{x_{1}-x_{2}}{x_{0}}\right)f(x_{1},x_{2})
Y_{W}(u,x_{1})Y_{W}(v,x_{2})\nonumber\\
& &\ \ \ -x_{0}^{-1}\delta\left(\frac{x_{2}-x_{1}}{-x_{0}}\right)
f(x_{1},x_{2})Y_{W}(v,x_{2})Y_{W}(u,x_{1})\nonumber\\
&=&x_{2}^{-1}\delta\left(\frac{x_{1}-x_{0}}{x_{2}}\right)
f(x_{1},x_{2})Y_{W}(Y(u,x_{0})v,x_{2}),
\end{eqnarray}
which implies
\begin{eqnarray}\label{ekpuv=0}
(x_{1}-x_{2})^{k}f(x_{1},x_{2})[Y_{W}(u,x_{1}),Y_{W}(v,x_{2})]=0
\end{eqnarray}
for some nonnegative integer $k$ with $x^{k}Y(u,x)v\in V[[x]]$.
This shows that $U$ is a quasi local
subspace of $\E(W)$. Furthermore, with (\ref{ekpuv=0}), 
using Proposition \ref{pjacobi-definition} we get
\begin{eqnarray}
& &x_{1}^{-1}\delta\left(\frac{x_{2}+x_{0}}{x_{1}}\right)
f(x_{1},x_{2})(x_{1}-x_{2})^{k}
\Y_{1}(Y_{W}(u,x_{2}),x_{0})Y_{W}(v,x_{2})\nonumber\\
&=&x_{0}^{-1}\delta\left(\frac{x_{1}-x_{2}}{x_{0}}\right)
f(x_{1},x_{2})(x_{1}-x_{2})^{k}Y_{W}(u,x_{1})Y_{W}(v,x_{2})\nonumber\\
& &-x_{0}^{-1}\delta\left(\frac{x_{2}-x_{1}}{-x_{0}}\right)
f(x_{1},x_{2})(x_{1}-x_{2})^{k}Y_{W}(v,x_{2})Y_{W}(u,x_{1})\nonumber\\
&=&x_{1}^{-1}\delta\left(\frac{x_{2}+x_{0}}{x_{1}}\right)f(x_{1},x_{2})
(x_{1}-x_{2})^{k}Y_{W}(Y(u,x_{0})v,x_{2}).
\end{eqnarray}
Taking $\Res_{x_{1}}$ we get
\begin{eqnarray}
f(x_{2}+x_{0},x_{2})x_{0}^{k}
\Y_{1}(Y_{W}(u,x_{2}),x_{0})Y_{W}(v,x_{2})
=f(x_{2}+x_{0},x_{2})x_{0}^{k}Y_{W}(Y(u,x_{0})v,x_{2}).
\end{eqnarray}
Since both $\Y_{1}(Y_{W}(u,x_{2}),x_{0})Y_{W}(v,x_{2})$ and 
$Y_{W}(Y(u,x_{0})v,x_{2})$ involve only finitely many negative powers
of $x_{0}$, using the cancelation law (Lemma \ref{lcancellation})
we get
\begin{eqnarray}\label{ehom-property}
\Y_{1}(Y_{W}(u,x_{2}),x_{0})Y_{W}(v,x_{2})
=Y_{W}(Y(u,x_{0})v,x_{2}).
\end{eqnarray}
This implies that $U$ is $\Y_{1}$-closed. The space $U$ contains $1_{W}$
since $Y_{W}({\bf 1},x)=1_{W}$. {}From (\ref{ehom-property}),
$(U,\Y_{1},1_{W})$ is a vertex algebra and
the linear map $Y_{W}$ is a vertex algebra homomorphism. 

Conversely, assume that $Y_{W}$ is a (vertex algebra)
homomorphism from $V$ to a $\Y_{1}$-closed
quasi local subspace containing $1_{W}$ of $\E(W)$. We have
$Y_{W}({\bf 1},x)=1_{W}$.
For $u,v\in V$, let $0\ne f(x_{1},x_{2})\in
\C[x_{1},x_{2}]$ such that
$$f(x_{1},x_{2})Y_{W}(u,x_{1})Y_{W}(v,x_{2})
=f(x_{1},x_{2})Y_{W}(v,x_{2})Y_{W}(u,x_{1}).$$
By Proposition \ref{pjacobi-definition} we have
\begin{eqnarray}
& &x_{0}^{-1}\delta\left(\frac{x_{1}-x_{2}}{x_{0}}\right)f(x_{1},x_{2})
Y_{W}(u,x_{1})Y_{W}(v,x_{2})\nonumber\\
& &-x_{0}^{-1}\delta\left(\frac{x_{2}-x_{1}}{-x_{0}}\right)f(x_{1},x_{2})
Y_{W}(v,x_{2})Y_{W}(u,x_{1})\nonumber\\
&=&x_{2}^{-1}\delta\left(\frac{x_{1}-x_{0}}{x_{2}}\right)f(x_{1},x_{2})
Y_{1}(Y_{W}(u,x_{2}),x_{0})Y_{W}(v,x_{2})\nonumber\\
&=&x_{2}^{-1}\delta\left(\frac{x_{1}-x_{0}}{x_{2}}\right)f(x_{1},x_{2})
Y_{W}(Y(u,x_{0})v,x_{2}).
\end{eqnarray}
Thus $(W,Y_{W})$ is a quasi $V$-module.
\end{proof}

The following is an analogue of
Theorem 5.7.6 of [LL]:

\bt{tmodule-extension}
Let $V$ be a vertex algebra with a generating
subspace $U$ and let $W$ be a vector space equipped with a linear map
$Y_{W}^{0}$ from $U$ to $\Hom (W,W((x)))$.
Then $Y_{W}^{0}$ can be extended to a linear map
$Y_{W}$ from $V$ to $\Hom (W,W((x)))$
such that $(W,Y_{W})$ carries the structure of a quasi $V$-module if
and only if $\bar{U}=\{Y_{W}^{0}(u,x)\;|\; u\in U\}$ is a quasi local
subspace of $\E(W)$ and there exists a linear map
$\psi$ from $V$ to $\<\bar{U}\>$ such that
\begin{eqnarray}
& &\psi({\bf 1})=1_{W},\label{e5.27}\\
& &\psi( Y(u,x_{0})v)=\Y_{1}( Y_{W}^{0}(u,x),x_{0})\psi(v)
\label{e5.28}
\end{eqnarray}
for $u\in U, \; v\in V$. Furthermore, such an extension is unique.
\et

\begin{proof} If $Y_{W}$ is an extension of $Y_{W}^{0}$ such that
$(W,Y_{W})$ carries the structure of a quasi $V$-module,
then by Theorem \ref{trepresentation},
$\{ Y_{W}(v,x)\;|\; v\in V\}$ is a $\Y_{1}$-closed quasi local
subspace of $\E(W)$ and $Y_{W}$ is a vertex algebra homomorphism.
Since $V$ is generated by $U$, $Y_{W}(V)$ is generated by
$\bar{U}$, i.e., $Y_{W}(V)=\<\bar{U}\>$. Taking $\psi=Y_{W}$ we have
(\ref{e5.27}) and (\ref{e5.28}).

On the other hand, assume that $\bar{U}$ is a quasi local
subspace of $\E(W)$ and $\psi$ is a linear map from $V$ to
$\<\bar{U}\>$ satisfying (\ref{e5.27}) and (\ref{e5.28}).
Since $V$ and $\<\bar{U}\>$ are vertex algebras and $U$ generates
$V$, by Proposition 5.7.9 of [LL] $\psi$ is a vertex algebra homomorphism.
In view of Theorem \ref{trepresentation}, $(W,Y_{W})$ carries the
structure of a quasi $V$-module with
$Y_{W}(v,x)=\psi_{x}(v)$ for $v\in V$.
For $u\in U$, using (\ref{e5.27}) and (\ref{e5.28}) we have
\begin{eqnarray}
Y_{W}(u,x)=\psi_{x}(u)=\psi_{x}(u_{-1}{\bf 1})
=Y_{W}^{0}(u,x)_{-1}1_{W}=Y_{W}^{0}(u,x).
\end{eqnarray}
Thus $Y_{W}$ extends $Y_{W}^{0}$.

Since $U$ generates $V$, any linear map from $U$ to a
vertex algebra $V'$ extends to at most
one vertex algebra homomorphism from
$V$ to $V'$. Then the uniqueness follows.
\end{proof}

The following is an analogue of
a result in \cite{li-local} (cf. [LL], Theorem 3.6.3):

\bp{p-quasi-module-weak}
Let $V$ be a vertex algebra and let $W$ be a vector space
equipped with a linear map $Y_{W}$ from $V$ to $\Hom (W,W((x)))$
such that $Y_{W}({\bf 1},x)=1_{W}$ and for $u,v\in V$ there exists
a nonzero polynomial $f(x_{1},x_{2})$ such that for $w\in W$
\begin{eqnarray}
& &f(x_{0}+x_{2},x_{2})(x_{0}+x_{2})^{l}Y_{W}(u,x_{0}+x_{2})Y_{W}(v,x_{2})w
\nonumber\\
&=&f(x_{0}+x_{2},x_{2})(x_{0}+x_{2})^{l}Y_{W}(Y(u,x_{0})v,x_{2})w
\end{eqnarray}
for some nonnegative integer $l$.
Then $(W,Y_{W})$ carries the structure of a quasi $V$-module.
\ep

\begin{proof} First we prove that
\begin{eqnarray}\label{edproperty-module}
Y_{W}(\D (u),x)=\frac{d}{dx}Y_{W}(u,x)\;\;\;\mbox{ for }u\in V.
\end{eqnarray}
Let $u\in V,\; w\in W$. Then there exist a nonzero polynomial
$f(x_{1},x_{2})$ and a nonnegative integer $l$ such that
\begin{eqnarray*}
& &f(x_{0}+x_{2},x_{2})(x_{0}+x_{2})^{l}Y_{W}(u,x_{0}+x_{2})Y_{W}({\bf 1},x_{2})w
\nonumber\\
&=&f(x_{0}+x_{2},x_{2})(x_{0}+x_{2})^{l}Y_{W}(Y(u,x_{0}){\bf 1},x_{2})w.
\end{eqnarray*}
With $Y_{W}({\bf 1},x)=1_{W}$ we have
\begin{eqnarray*}
f(x_{0}+x_{2},x_{2})(x_{0}+x_{2})^{l}Y_{W}(u,x_{0}+x_{2})w
=f(x_{0}+x_{2},x_{2})(x_{0}+x_{2})^{l}Y_{W}(Y(u,x_{0}){\bf 1},x_{2})w.
\end{eqnarray*}
We may replace $l$ with a bigger integer if necessary so that
$x^{l}Y_{W}(u,x)w\in W[[x]]$. Then
$$(x_{0}+x_{2})^{l}Y_{W}(u,x_{0}+x_{2})w
=(x_{0}+x_{2})^{l}Y_{W}(u,x_{2}+x_{0})w.$$
We thus have
\begin{eqnarray*}
f(x_{0}+x_{2},x_{2})(x_{0}+x_{2})^{l}Y_{W}(u,x_{2}+x_{0})w
=f(x_{0}+x_{2},x_{2})(x_{0}+x_{2})^{l}Y_{W}(Y(u,x_{0}){\bf 1},x_{2})w.
\end{eqnarray*}
In view of the cancelation law (Lemma \ref{lcancellation})
we have
\begin{eqnarray*}
Y_{W}(u,x_{2}+x_{0})w
=Y_{W}(Y(u,x_{0}){\bf 1},x_{2})w,
\end{eqnarray*}
which implies $\frac{d}{dx}Y_{W}(u,x)w=Y_{W}(\D (u),x)w$.
This proves (\ref{edproperty-module}).

Let $u,v\in V$. Using the skew symmetry of $V$ and 
(\ref{edproperty-module}) we have
\begin{eqnarray}\label{eproof-skew-3}
Y_{W}(Y(u,x_{0})v,x_{2})
=Y_{W}(e^{x_{0}\D}Y(v,-x_{0})u,x_{2})
=Y_{W}(Y(v,-x_{0})u,x_{2}+x_{0}).
\end{eqnarray}
Let $w\in W$. There exist nonzero polynomials $f(x_{1},x_{2})$ and $g(x_{1},x_{2})$
depending only on $u,v$ and nonnegative integers $l,l'$ such that
\begin{eqnarray}
& &(x_{0}+x_{2})^{l}f(x_{0}+x_{2},x_{2})Y_{W}(u,x_{0}+x_{2})Y_{W}(v,x_{2})w\nonumber\\
&=&(x_{0}+x_{2})^{l}f(x_{0}+x_{2},x_{2})Y_{W}(Y(u,x_{0})v,x_{2})w,\label{efguv-1}\\
& &(-x_{0}+x_{1})^{l'}g(x_{1},-x_{0}+x_{1})Y_{W}(v,-x_{0}+x_{1})Y_{W}(u,x_{1})w\nonumber\\
&=&(-x_{0}+x_{1})^{l'}g(x_{1},-x_{0}+x_{1})Y_{W}(Y(v,-x_{0})u,x_{1})w.\label{efgvu-1}
\end{eqnarray}
Replacing $l$ with a larger one if necessary
so that $x^{l}Y_{W}(u,x)w\in W[[x]]$. Using (\ref{efgvu-1}) we see that
the expression 
$$x_{1}^{l}(-x_{0}+x_{1})^{l'}g(x_{1},-x_{0}+x_{1})Y_{W}(Y(v,-x_{0})u,x_{1})w$$
lies in $W[[x_{0},x_{1}]][x_{0}^{-1}]$.
In view of the delta function substitution rule (\ref{e-delta-2}) we have
\begin{eqnarray}\label{esubstitution-2}
& &x_{0}^{-1}\delta\left(\frac{x_{2}-x_{1}}{-x_{0}}\right)
\left(x_{1}^{l}f(x_{1},x_{2})(-x_{0}+x_{1})^{l'}g(x_{1},-x_{0}+x_{1})
Y_{W}(Y(v,-x_{0})u,x_{1})w\right)\nonumber\\
&=&x_{0}^{-1}\delta\left(\frac{x_{2}-x_{1}}{-x_{0}}\right)
\left(x_{1}^{l}f(x_{1},x_{2})(-x_{0}+x_{1})^{l'}g(x_{1},-x_{0}+x_{1})
Y_{W}(Y(v,-x_{0})u,x_{1})w\right)|_{x_{1}=x_{2}+x_{0}}\nonumber\\
&=&x_{0}^{-1}\delta\left(\frac{x_{2}-x_{1}}{-x_{0}}\right)
\left((x_{2}+x_{0})^{l}f(x_{2}+x_{0},x_{2})x_{2}^{l'}g(x_{1},x_{2})
Y_{W}(Y(v,-x_{0})u,x_{2}+x_{0})w\right).\nonumber\\
& &
\end{eqnarray}
Using (\ref{efguv-1}), (\ref{efgvu-1}), (\ref{eproof-skew-3})
and (\ref{esubstitution-2}) in a sequence we get
\begin{eqnarray}
& &x_{0}^{-1}\delta\left(\frac{x_{1}-x_{2}}{x_{0}}\right)
x_{1}^{l}f(x_{1},x_{2})x_{2}^{l'}g(x_{1},x_{2})Y_{W}(u,x_{1})Y_{W}(v,x_{2})w
\nonumber\\
& &-x_{0}^{-1}\delta\left(\frac{x_{2}-x_{1}}{-x_{0}}\right)
x_{1}^{l}f(x_{1},x_{2})x_{2}^{l'}g(x_{1},x_{2})Y_{W}(v,x_{2})Y_{W}(u,x_{1})w
\nonumber\\
&=&x_{0}^{-1}\delta\left(\frac{x_{1}-x_{2}}{x_{0}}\right)
x_{2}^{l'}g(x_{1},x_{2})(x_{0}+x_{2})^{l}f(x_{0}+x_{2},x_{2})Y_{W}(u,x_{0}+x_{2})Y_{W}(v,x_{2})w
\nonumber\\
& &-x_{0}^{-1}\delta\left(\frac{x_{2}-x_{1}}{-x_{0}}\right)
x_{1}^{l}f(x_{1},x_{2})(-x_{0}+x_{1})^{l'}g(x_{1},-x_{0}+x_{1})Y_{W}(v,-x_{0}+x_{1})Y_{W}(u,x_{1})w
\nonumber\\
&=&x_{0}^{-1}\delta\left(\frac{x_{1}-x_{2}}{x_{0}}\right)
\left(x_{2}^{l'}g(x_{1},x_{2})(x_{0}+x_{2})^{l}f(x_{0}+x_{2},x_{2})Y_{W}(Y(u,x_{0})v,x_{2})w\right)
\nonumber\\
& &-x_{0}^{-1}\delta\left(\frac{x_{2}-x_{1}}{-x_{0}}\right)
\left(x_{1}^{l}f(x_{1},x_{2})(-x_{0}+x_{1})^{l'}g(x_{1},-x_{0}+x_{1})
Y_{W}(Y(v,-x_{0})u,x_{1})w\right)
\nonumber\\
&=&x_{0}^{-1}\delta\left(\frac{x_{1}-x_{2}}{x_{0}}\right)
\left(x_{2}^{l'}g(x_{1},x_{2})(x_{0}+x_{2})^{l}f(x_{0}+x_{2},x_{2})
Y_{W}(Y(v,-x_{0})u,x_{2}+x_{0})w\right)
\nonumber\\
& &-x_{0}^{-1}\delta\left(\frac{x_{2}-x_{1}}{-x_{0}}\right)
\left(x_{1}^{l}f(x_{1},x_{2})(x_{0}+x_{2})^{l'}g(x_{1},-x_{0}+x_{1})
Y_{W}(Y(v,-x_{0})u,x_{1})w\right)
\nonumber\\
&=&x_{0}^{-1}\delta\left(\frac{x_{1}-x_{2}}{x_{0}}\right)
\left(x_{2}^{l'}g(x_{1},x_{2})(x_{0}+x_{2})^{l}f(x_{0}+x_{2},x_{2})
Y_{W}(Y(v,-x_{0})u,x_{2}+x_{0})w\right)
\nonumber\\
& &-x_{0}^{-1}\delta\left(\frac{x_{2}-x_{1}}{-x_{0}}\right)
\left(x_{2}^{l'}g(x_{1},x_{2})(x_{0}+x_{2})^{l}f(x_{0}+x_{2},x_{2})
Y_{W}(Y(v,-x_{0})u,x_{2}+x_{0})w\right)
\nonumber\\
&=&x_{1}^{-1}\delta\left(\frac{x_{2}+x_{0}}{x_{1}}\right)
x_{1}^{l}f(x_{1},x_{2})x_{2}^{l'}g(x_{1},x_{2})Y_{W}(Y(u,x_{0})v,x_{2})w.
\end{eqnarray}
Multiplying by $x_{1}^{-l}x_{2}^{-l'}$ we obtain the quasi Jacobi identity
for $(u,v)$. This completes the proof.
\end{proof}

The following result is an analogue of a theorem  of [LL]
(see also [R]):

\bt{tall-vsas}
Let $W$ be a vector space and let $U$ be a subspace of $\E(W)$
in which every pair is compatible. Assume that $U$ is $\Y_{1}$-closed 
and it contains $1_{W}$. 
Then $(U,\Y_{1},1_{W})$ is a vertex algebra if and only if
$U$ is quasi local.
\et

\begin{proof} In view of Theorem \ref{tkey-main}
we only need to prove that
if $(U,\Y_{1},1_{W})$ is a vertex algebra then $U$ is quasi local.
For $a(x),b(x)\in U$,
let $f(x_{1},x_{2})$ be a nonzero polynomial
such that
$$f(x_{1},x_{2})a(x_{1})b(x_{2})\in \Hom (W,W((x_{1},x_{2}))).$$
For $w\in W$, let $l$ be a nonnegative integer such that
$$x_{1}^{l}f(x_{1},x_{2})a(x_{1})b(x_{2})w\in 
W[[x_{1},x_{2}]][x_{2}^{-1}].$$
By Proposition \ref{pexistence} we have
\begin{eqnarray}
& &f(x_{0}+x_{2},x_{2})(x_{0}+x_{2})^{l}
\left(\Y(a(x_{2}),x_{0})b(x_{2})\right)w
\nonumber\\
&=&f(x_{0}+x_{2},x_{2})(x_{0}+x_{2})^{l}a(x_{0}+x_{2})b(x_{2})w.
\end{eqnarray}
In view of Proposition \ref{p-quasi-module-weak},
$W$ is a quasi $U$-module with $Y_{W}(a(x),x_{0})=a(x_{0})$ for
$a(x)\in U$. It follows from the quasi Jacobi identity
that any $a(x), b(x)$ in $U$ are quasi local.
This proves that $U$ is quasi local.
\end{proof}

An important problem could be the classification and construction of
irreducible quasi modules for a vertex algebra. In the following
we present certain results in this direction.

\bl{lpolynomial2}
Let $V$ be a vertex algebra and $(W,Y_{W})$ a quasi $V$-module.
Let $a,b,c\in V$ and let 
$$0\ne f(x_{1},x_{2}), g(x_{1},x_{2}), h(x_{1},x_{2})\in
\C[x_{1},x_{2}]$$
 such that
\begin{eqnarray}
& &f(x_{1},x_{2})[Y_{W}(a,x_{1}),Y_{W}(b,x_{2})]=0,\\
& &g(x_{1},x_{2})[Y_{W}(a,x_{1}),Y_{W}(c,x_{2})]=0,\\
& &h(x_{1},x_{2})[Y_{W}(b,x_{1}),Y_{W}(c,x_{2})]=0.
\end{eqnarray}
Then  for any $n\in \Z$, there exists a nonnegative integer $k$ such
that
\begin{eqnarray}
f(x_{1},x_{2})^{k}g(x_{1},x_{2})[Y_{W}(a,x_{1}),Y_{W}(b_{n}c,x_{2})]=0.
\end{eqnarray}
\el

\begin{proof} In view of Theorem \ref{trepresentation}, 
for any $u,v\in V$, $Y_{W}(u,x)$ and $Y_{W}(v,x)$ are quasi local and
\begin{eqnarray}
Y_{W}(u_{m}v,x)=Y_{W}(u,x)_{m}Y_{W}(v,x)\;\;\;\mbox{ for }m\in \Z.
\end{eqnarray}
{}From Proposition \ref{pgamma-locality-key} 
there exists a nonnegative integer $k$ such that
\begin{eqnarray}
f(x_{1},x_{2})^{k}g(x_{1},x_{2})
\left[Y_{W}(a,x_{1}),Y_{W}(b,x_{2})_{n}Y_{W}(c,x_{2})\right]=0.
\end{eqnarray}
Thus we obtain
\begin{eqnarray}
f(x_{1},x_{2})^{k}g(x_{1},x_{2})
[Y_{W}(a,x_{1}),Y_{W}(b_{n}c,x_{2})]=0,
\end{eqnarray}
completing the proof.
\end{proof}

{}From Proposition \ref{lpolynomial2} and induction we immediately
have:

\bp{ppolynomial}
Let $V$ be a vertex algebra and $(W,Y_{W})$ a quasi $V$-module.
Let $U$ be a generating subspace of $V$ and let 
$0\ne f(x_{1},x_{2})\in \C[x_{1},x_{2}]$
such that for any $u,u'\in U$ there exists a nonnegative integer $k$
such that
\begin{eqnarray}
f(x_{1},x_{2})^{k}[Y_{W}(u,x_{1}),Y_{W}(u',x_{2})]=0.
\end{eqnarray}
Then for any $v,v'\in V$, there exists a nonnegative integer $r$
such that
\begin{eqnarray}
f(x_{1},x_{2})^{r}[Y_{W}(v,x_{1}),Y_{W}(v',x_{2})]=0.
\end{eqnarray}
\ep

As immediate consequences we have:

\bc{cpolynomial}
Let $W$ be any vector space and let $U$ be a quasi local subspace of
$\E(W)$. Assume that $f(x_{1},x_{2})$ is a nonzero polynomial such that
for any $a(x),b(x)\in U$ there exists a nonnegative integer $k$
such that
\begin{eqnarray}
f(x_{1},x_{2})^{k}[a(x_{1}),b(x_{2})]=0.
\end{eqnarray}
Let $V=\<U\>$ be the vertex algebra generated by $U$.
Then for any $\alpha (x),\beta(x)\in V$, 
there exists a nonnegative integer $r$
such that
\begin{eqnarray}
f(x_{1},x_{2})^{r}[\alpha(x_{1}),\beta(x_{2})]=0.
\end{eqnarray}
\ec

\bc{cpolynomial-finite-generating}
Let $V$ be a finitely generated vertex algebra and let $(W,Y_{W})$ be
any quasi $V$-module. Then there exists a nonzero polynomial
$f(x_{1},x_{2})$ such that for any $u,v\in V$,
\begin{eqnarray}
f(x_{1},x_{2})^{r}[Y_{W}(u,x_{1}),Y_{W}(v,x_{2})]=0
\end{eqnarray}
for some nonnegative integer $r$.
\ec

In view of Corollary \ref{cpolynomial-finite-generating}
we would like to associate a ``minimal'' polynomial $f(x_{1},x_{2})$
to a quasi module $W$ for a vertex algebra $V$.
Since $\C[x_{1},x_{2}]$ is not a principal domain,
we have a technical problem uniquely to define such a term.
Now we consider a special case.
Assume that there exists a nonzero homogeneous polynomial $f(x_{1},x_{2})$
such that the assertion of Corollary
\ref{cpolynomial-finite-generating} holds.
Notice that a nonzero polynomial $g(x_{1},x_{2})$ is 
homogeneous if and only if $g(x_{1},x_{2})=x_{2}^{k}p(x_{1}/x_{2})$
for some nonzero polynomial $p(x)$ of degree $k$.

\bd{dminimal}
{\em Let $V$ be a vertex algebra and $W$ a quasi $V$-module
satisfying that there exists a nonzero polynomial $p(x)$ 
such that for any $u,v\in V$,
\begin{eqnarray}
p(x_{1}/x_{2})^{r}[Y_{W}(u,x_{1}),Y_{W}(v,x_{2})]=0
\end{eqnarray}
for some nonnegative integer $r$.
The monic polynomial $p(x)$ of the least degree
is called the {\em minimal polynomial} of $W$.}
\ed

\section{$\Gamma$-vertex algebras and (quasi) modules}

Motivated by Theorem \ref{tkey-main} and by the notion of
$\Gamma$-conformal algebra in [G-K-K], in this section we formulate
and study a notion of $\Gamma$-vertex algebra, where $\Gamma$ is an
arbitrary group. We show that a $\Gamma$-vertex algebra amounts to an
(ordinary) vertex algebra equipped with a compatible $\Gamma$-module
structure in a certain sense.  We also formulate a notion of quasi
module for a $\Gamma$-vertex algebra.  In terms of these notions, we
strengthen Theorem \ref{tkey-main}, showing that for any vector space
$W$, any subgroup $\Gamma$ of $\C^{\times}$ and for any
$\Y_{\gamma}$-closed quasi local subspace $U$ of $\E(W)$ containing
$1_{W}$, $(U,\{\Y_{\alpha}\},1_{W})$ carries the structure of a
$\Gamma$-vertex algebra with $W$ as a quasi module.

Let $\Gamma$ be any (possibly nonabelian and infinite) group. 
We shall use $1$ for the identity element of $\Gamma$.
Let
\begin{eqnarray}
\phi: \Gamma \rightarrow \C^{\times}
\end{eqnarray}
be a group homomorphism.
We fix the pair $(\Gamma,\phi)$ throughout this section.

An important example for such a pair $(\Gamma,\phi)$ is the one with
$\Gamma$ a subgroup of $\C^{\times}$ and with $\phi$ the identity map.
As a convention, whenever $\Gamma$ is given as a subgroup
of $\C^{\times}$ {\em we always assume that $\phi$ is the identity map}.

With the group homomorphism $\phi$, any $\C$-vector space $W$
naturally becomes a $\Gamma$-module by defining
\begin{eqnarray}
\alpha w=\phi(\alpha)w\;\;\;\mbox{ for }\alpha\in \Gamma,\; w\in W.
\end{eqnarray}
In particular, for formal variable $x$, with $\C[[x,x^{-1}]]$ as 
a $\C$-vector space we have
\begin{eqnarray}
\alpha x=\phi(\alpha) x\;\;\;\mbox{ for }\alpha\in \Gamma.
\end{eqnarray}
We have
\begin{eqnarray}
\alpha \beta x=\beta \alpha x\;\;\;\mbox{ for }\alpha,\beta\in \Gamma,
\end{eqnarray}
since $\phi(\alpha\beta)=\phi(\alpha)\phi(\beta)=\phi(\beta)\phi(\alpha)=\phi(\beta\alpha)$.
Due to our convention, whenever $\Gamma\subset \C^{\times}$,
the new $\Gamma$-action through $\phi$ agrees with
the old $\Gamma$-action by scalar multiplication.

The following notion naturally arises from Theorem \ref{tkey-main}:

\bd{dgammava}
{\em A {\em $\Gamma$-vertex algebra} is a vector space $V$ equipped
with linear maps
\begin{eqnarray}
Y_{\alpha}: V&\rightarrow & \Hom (V,V((x)))\subset (\End
V)[[x,x^{-1}]]\nonumber\\
v&\mapsto& Y_{\alpha}(v,x)
\end{eqnarray}
for $\alpha\in \Gamma$ and equipped with a distinguished vector ${\bf
1}$, called the {\em vacuum vector}, such that all the following
axioms hold: For $\alpha\in \Gamma,\;v\in V$,
\begin{eqnarray}
& &Y_{\alpha}({\bf 1},x)=1,\label{egamma-vacuum}\\
& &Y_{\alpha}(v,x){\bf 1}=V[[x]]\;\;\mbox{ and }\;
\lim_{x\rightarrow 0}Y_{1}(v,x){\bf 1}=v\label{egamma-creation}
\end{eqnarray}
and for $u,v,w\in V,\; \alpha,\beta\in \Gamma$, 
\begin{eqnarray}\label{egamma-jacobi}
& &x_{0}^{-1}\delta\left(\frac{x_{1}-\beta^{-1}\alpha x_{2}}{x_{0}}\right)
Y_{\alpha}(u,x_{1})Y_{\beta}(v,x_{2})w\nonumber\\
& &-x_{0}^{-1}\delta\left(\frac{-\beta^{-1}\alpha x_{2}+
x_{1}}{x_{0}}\right)
Y_{\beta}(v,x_{2})Y_{\alpha}(u,x_{1})w\nonumber\\
&=&x_{1}^{-1}\delta\left(\frac{\beta^{-1}\alpha x_{2}+x_{0}}{x_{1}}\right)
Y_{\beta}(Y_{\beta^{-1}\alpha}(u,x_{0})v,x_{2})w.
\end{eqnarray}}
\ed

The $\Gamma$-vertex algebra $V$ is also often denoted by 
$(V,\{Y_{\alpha}\},{\bf 1})$.
For convenience, we refer (\ref{egamma-jacobi}) as
the {\em $\Gamma$-Jacobi identity for the triple $(u,v,w)$.}

\br{rva=gammava}
{\em It is clear that 
the notion of $\Gamma$-vertex algebra with $\Gamma=\{1\}$
reduces to the notion of ordinary vertex algebra.}
\er

Similar to the notion of quasi module for a vertex algebra,
we have the following notion of quasi module for a $\Gamma$-vertex algebra:

\bd{dquasi-modules}
{\em Let $V$ be a $\Gamma$-vertex algebra.
A {\em quasi $V$-module} is a vector space $W$ equipped with linear maps
\begin{eqnarray}
Y^{W}_{\alpha}: V\rightarrow \Hom (W,W((x)))
\end{eqnarray}
for $\alpha\in \Gamma$ such that the following conditions hold:
\begin{eqnarray}
Y^{W}_{\alpha}({\bf 1},x)=1_{W}\;\;\;\mbox{ for }\alpha\in \Gamma
\end{eqnarray}
and for $u,v\in V$, there exists a nonzero polynomial $f(x_{1},x_{2})$
such that
\begin{eqnarray}\label{epjacobi-module}
& &x_{0}^{-1}\delta\left(\frac{x_{1}-\beta^{-1}\alpha x_{2}}{x_{0}}\right)
f(x_{1},x_{2})Y^{W}_{\alpha}(u,x_{1})Y^{W}_{\beta}(v,x_{2})\nonumber\\
& &-x_{0}^{-1}\delta\left(\frac{-\beta^{-1}\alpha x_{2}+
x_{1}}{x_{0}}\right)f(x_{1},x_{2})
Y^{W}_{\beta}(v,x_{2})Y^{W}_{\alpha}(u,x_{1})\nonumber\\
&=&x_{1}^{-1}\delta\left(\frac{\beta^{-1}\alpha x_{2}+x_{0}}{x_{1}}\right)
f(x_{1},x_{2})Y^{W}_{\beta}(Y_{\beta^{-1}\alpha}(u,x_{0})v,x_{2})
\end{eqnarray}
for $\alpha,\beta\in \Gamma$. A quasi $V$-module $W$ is called a {\em module} if
(\ref{epjacobi-module}) holds with $f(x_{1},x_{2})=1$ 
for any $u,v\in V,\; \alpha,\beta\in \Gamma$.}
\ed

Notice that if an (ordinary) vertex algebra $V$ is viewed as
a $\Gamma$-vertex algebra with $\Gamma=\{1\}$, 
the notion of quasi module for $V$ as a $\Gamma$-vertex algebra
exactly gives the notion of quasi module for
$V$ as a vertex algebra.
 
In terms of these notions we have:

\bt{tmain1}
Let $W$ be a vector space, $\Gamma$ a subgroup of $\C^{\times}$ 
and $V$ a $\Y_{\Gamma}$-closed quasi local subspace of $\E(W)$
containing $1_{W}$.
Then $(U,\{\Y_{\alpha}\},1_{W})$ carries 
the structure of a $\Gamma$-vertex algebra 
with $W$ a quasi module with
$Y_{\alpha}^{W}(a(x),x_{0})=a(x_{0})$ for $a(x)\in V,\; \alpha\in \Gamma$. 
Furthermore, for any quasi local subset $S$ of $\E(W)$,
$\<S\>_{\Gamma}$ is
a $\Gamma$-vertex algebra with $W$ as a quasi module.
\et

\begin{proof} Theorem \ref{tkey-main} exactly states that
 $(V,\{\Y_{\alpha}\},1_{W})$ carries 
the structure of a $\Gamma$-vertex algebra.
Let $a(x),b(x)\in V$. There exists a nonzero polynomial $f(x_{1},x_{2})$ such that
\begin{eqnarray}
f(x_{1},x_{2})a(x_{1})b(x_{2})=f(x_{1},x_{2})b(x_{2})a(x_{1}).
\end{eqnarray}
This implies that $f(x_{1},x_{2})a(x_{1})b(x_{2})\in \Hom (W,W((x_{1},x_{2})))$.
Let $w\in W$ and let $l$ be a nonnegative integer such that
$x^{l}a(x)w\in W[[x]]$, so that
$$x_{1}^{l}f(x_{1},x_{2})a(x_{1})b(x_{2})w\in W[[x_{1},x_{2}]][x_{2}^{-1}].$$
For $\alpha,\beta\in \C^{\times}$, we have
\begin{eqnarray}
f(x_{1},x_{2})Y^{W}_{\alpha}(a(x),x_{1})Y^{W}_{\beta}(b(x),x_{2})
&=&f(x_{1},x_{2})a(x_{1})b(x_{2})\nonumber\\
&=&f(x_{1},x_{2})b(x_{2})a(x_{1})\nonumber\\
&=&f(x_{1},x_{2})Y^{W}_{\beta}(b(x),x_{2})Y^{W}_{\alpha}(a(x),x_{1}).
\end{eqnarray}
With
\begin{eqnarray*}
& &Y_{\alpha}^{W}(a(x),x_{0}+\beta^{-1}\alpha x_{2})Y_{\beta}(b(x),x_{2})w
=a(x_{0}+\beta^{-1}\alpha x_{2})b(x_{2})w,\\
& &Y_{\beta}^{W}(\Y_{\beta^{-1}\alpha}(a(x),x_{0})b(x),x_{2})w
=\left(\Y_{\beta^{-1}\alpha}(a(x),x_{0})b(x)\right)w|_{x=x_{2}},
\end{eqnarray*}
{}from Proposition \ref{pexistence}  we also have
\begin{eqnarray}
& &(x_{0}+\beta^{-1}\alpha x_{2})^{l}f(x_{0}+\beta^{-1}\alpha x_{2},x_{2})
Y_{\alpha}^{W}(a(x),x_{0}+\beta^{-1}\alpha x_{2})
Y_{\beta}(b(x),x_{2})w\nonumber\\
&=&(x_{0}+\beta^{-1}\alpha x_{2})^{l}f(x_{0}+\beta^{-1}\alpha x_{2},x_{2})
Y_{\beta}^{W}(\Y_{\beta^{-1}\alpha}(a(x),x_{0})b(x),x_{2})w.
\end{eqnarray}
In view of Lemma \ref{lformaljacobiidentity}, 
the required Jacobi-like identity
holds. By definition $Y_{\alpha}^{W}(1_{W},x_{0})=1_{W}$. 
This proves that $W$ is a quasi module.
\end{proof}

For the rest of this section we shall study the structures of
$\Gamma$-vertex algebras and their (quasi) modules in terms of
ordinary vertex algebras and modules.

\bl{ldproperties}
Let $(V,\{Y_{\alpha}\},{\bf 1})$ be a $\Gamma$-vertex algebra.
For $\alpha\in \Gamma$, define $R_{\alpha}\in \End_{\C}V$ by
\begin{eqnarray}\label{edef-ralpha}
R_{\alpha}(v)=\Res_{x}x^{-1}Y_{\alpha}(v,x){\bf 1}
=\lim_{x\rightarrow 0}Y_{\alpha}(v,x){\bf 1}\;\;\;\mbox{ for }v\in V.
\end{eqnarray}
Define $\D\in \End_{\C}V$ by
\begin{eqnarray}
\D v=\Res_{x}x^{-2}Y_{1}(v,x){\bf 1}\;\;\;\mbox{ for }v\in V.
\end{eqnarray}
Then
\begin{eqnarray}
& & Y_{\alpha}(\D v,x)={d\over dx}Y_{\alpha}(v,x),\\
& &[\D,Y_{\alpha}(v,x)]=\alpha {d\over dx}Y_{\alpha}(v,x),\\
& & R_{\alpha} Y_{\beta}(v,x)=Y_{\alpha\beta}(v,x)R_{\alpha}
\;\;\;\mbox{ for }\alpha,\beta\in \Gamma,\; v\in V.
\label{etranslation}
\end{eqnarray}
\el

\begin{proof}
Taking $\Res_{x_{1}}$ of the $\Gamma$-Jacobi identity for the triple
$(u,v,{\bf 1})$ with $\beta=1$ and using the property that
$Y_{\alpha}(u,x_{1}){\bf 1}\in V[[x_{1}]]$ we get
\begin{eqnarray}
Y_{1}(Y_{\alpha}(u,x_{0})v,x_{2}){\bf 1}
&=&Y_{\alpha}(u,x_{0}+\alpha x_{2})Y_{1}(v,x_{2}){\bf 1}\nonumber\\
&=&e^{\alpha x_{2}\frac{\partial}{\partial x_{0}}}
Y_{\alpha}(u,x_{0})Y_{1}(v,x_{2}){\bf 1}.
\end{eqnarray}
Then
\begin{eqnarray}
\D Y_{\alpha}(u,x_{0})v
=\alpha \frac{\partial}{\partial x_{0}}Y_{\alpha}(u,x_{0})v+
Y_{\alpha}(u,x_{0})\D v,
\end{eqnarray}
which shows that 
$[\D,Y_{\alpha}(v,x)]=\alpha \frac{d}{dx}Y_{\alpha}(v,x)$ for $v\in V$.

Taking $\Res_{x_{1}}$ of the $\Gamma$-Jacobi identity for the triple
($u,{\bf 1},v)$ with $\alpha=\beta$, using $Y_{\alpha}({\bf 1},x)=1$
we have
\begin{eqnarray}\label{eweak-assoc-uv1}
Y_{\alpha}(Y_{1}(u,x_{0}){\bf 1},x_{2})v
=Y_{\alpha}(u,x_{2}+x_{0})v
=e^{x_{0}\frac{\partial}{\partial x_{2}}}
Y_{\alpha}(u,x_{2})v.
\end{eqnarray}
This gives $Y_{\alpha}(\D u,x)={d\over dx}Y_{\alpha}(v,x)$.

Similarly, from the $\Gamma$-Jacobi identity 
for the triple $(u,v,{\bf 1})$ we have
\begin{eqnarray}
Y_{\beta}(Y_{\beta^{-1}\alpha}(u,x_{0})v,x_{2}){\bf 1}
=Y_{\alpha}(u,x_{0}+\beta^{-1}\alpha x_{2})Y_{\beta}(v,x_{2}){\bf 1}.
\end{eqnarray}
Setting $x_{2}=0$ we get
\begin{eqnarray}
R_{\beta}Y_{\beta^{-1}\alpha}(u,x_{0})v=Y_{\alpha}(u,x_{0})R_{\beta}v,
\end{eqnarray}
proving (\ref{etranslation}).
\end{proof}

In [G-K-K], $\Gamma$-conformal algebras are classified as Lie algebras
equipped with an action of $\Gamma$ by automorphisms of the Lie
algebras.  In the same spirit, the following theorem classifies
$\Gamma$-vertex algebras in terms of (ordinary) vertex algebras
equipped with a ``compatible'' $\Gamma$-action:

\bt{tsingle-mult}
Let $(V,\{Y_{\alpha}\},{\bf 1})$ be a $\Gamma$-vertex algebra
and let $R_{\alpha}\in \Aut_{\C}V$ for $\alpha\in \Gamma$
be defined as in Lemma \ref{ldproperties}.
Then $(V,Y_{1},{\bf 1})$ is a vertex algebra and
the map $\alpha\in \Gamma \rightarrow R_{\alpha}\in \End_{\C}V$ 
is a representation of $\Gamma$ on $V$ with $R_{\alpha}({\bf 1})={\bf 1}$ 
for $\alpha\in \Gamma$.
Furthermore, 
\begin{eqnarray}\label{e6.25}
Y_{\alpha}(v,x)=R_{\alpha}Y_{1}(v,x)R_{\alpha^{-1}}=
Y_{1}(R_{\alpha}v,\alpha^{-1}x)
\;\;\;\mbox{ for }\alpha\in \Gamma,\; v\in V.
\end{eqnarray}
On the other hand, let $V$ be an (ordinary)  vertex algebra. 
Assume that $V$ is a
$\Gamma$-module with $\alpha\in \Gamma$ 
acting as $R_{\alpha}\in \Aut_{\C} V$ such that
$R_{\alpha}({\bf 1})={\bf 1}$  and 
\begin{eqnarray}\label{e-abstract-conjugation}
R_{\alpha}Y(v,x)R_{\alpha}^{-1}=Y(R_{\alpha}(v),\alpha^{-1}x)
\;\;\;\mbox{ for }\alpha\in \Gamma,\; v\in V.
\end{eqnarray}
For $\alpha\in \Gamma$, define $Y_{\alpha}\in (\End V)[[x,x^{-1}]]$ by
\begin{eqnarray}
Y_{\alpha}(v,x)=R_{\alpha}Y(v,x)R_{\alpha}^{-1}=Y(R_{\alpha}(v),\alpha^{-1}x)
\;\;\;\mbox{ for }v\in V.
\end{eqnarray}
Then $(V,\{Y_{\alpha}\}, {\bf 1})$ is a $\Gamma$-vertex algebra.
\et

\begin{proof} 
Let $(V,\{Y_{\alpha}\},{\bf 1})$ be a $\Gamma$-vertex algebra. 
Clearly, $(V,Y_{1},{\bf 1})$ is an (ordinary)  vertex algebra.
Using the vacuum property (\ref{egamma-vacuum}), the creation property
(\ref{egamma-creation}) and the definition of $R_{\alpha}$ we immediately have
that $R_{\alpha}({\bf 1})={\bf 1}$ for $\alpha\in \Gamma$ and $R_{1}=1_{V}$.
Taking $\Res_{x_{1}}$ of the $\Gamma$-Jacobi identity for the triple $(u,{\bf 1},{\bf 1})$, 
we get
$$Y_{\alpha}(u,x_{0}+\beta^{-1}\alpha x_{2})Y_{\beta}({\bf 1},x_{2}){\bf 1}
=Y_{\beta}(Y_{\beta^{-1}\alpha}(u,x_{0}){\bf 1},x_{2}){\bf 1}.$$
With $Y_{\beta}({\bf 1},x)=1$, setting $x_{0}=x_{2}=0$  we obtain
\begin{eqnarray}
R_{\alpha}(u)=R_{\beta}R_{\beta^{-1}\alpha}(u).
\end{eqnarray}
Thus, the map $\alpha\mapsto R_{\alpha}$ is a representation of 
$\Gamma$ on $V$.

Taking $\Res_{x_{1}}$ of the $\Gamma$-Jacobi identity for the triple ($u,{\bf 1},v)$
with $\beta=1$ we have
$$Y_{1}(Y_{\alpha}(u,x_{0}){\bf 1},x_{2})v
=Y_{\alpha}(u,\alpha x_{2}+x_{0})v.$$
Setting $x_{0}=0$, we get
$$Y_{1}(R_{\alpha}(u),x_{2})v=Y_{\alpha}(u,\alpha x_{2})v,$$
which gives
\begin{eqnarray}\label{e5.29}
Y_{\alpha}(u,x)=Y_{1}(R_{\alpha}(u),\alpha^{-1}x)\;\;\;\mbox{ for }u\in V.
\end{eqnarray}
{}From Lemma \ref{ldproperties} we also have
$R_{\alpha}Y_{1}(v,x)=Y_{\alpha}(v,x)R_{\alpha}$ for $\alpha\in
\Gamma,\; v\in V$.

On the other hand, let $V$ be a vertex algebra equipped with a
$\Gamma$-module structure with $\alpha$ acting as $R_{\alpha}\in
\Aut_{\C}V$ $(\alpha\in \Gamma)$ with the given properties.  For
$\alpha\in \Gamma,\; v\in V$, we have
\begin{eqnarray}
& &Y_{\alpha}({\bf 1},x)=Y(R_{\alpha}{\bf 1},\alpha^{-1}x)
=Y({\bf 1},\alpha^{-1}x)=1,\\
& &Y_{\alpha}(v,x){\bf 1}=Y(R_{\alpha}v,\alpha^{-1}x){\bf 1}
=e^{\alpha^{-1}x\D}R_{\alpha}v\in V[[x]],\\
& &Y_{1}(v,x){\bf 1}=e^{x\D}R_{1}v=e^{x\D}v.
\end{eqnarray}
For $u,v\in V,\; \alpha,\beta\in \Gamma$, we have
\begin{eqnarray}
& &x_{0}^{-1}\delta\left(\frac{\alpha^{-1}x_{1}-\beta^{-1}x_{2}}
{x_{0}}\right)
Y(R_{\alpha}u,\alpha^{-1}x_{1})Y(R_{\beta}v,\beta^{-1}x_{2})\nonumber\\
& &-x_{0}^{-1}\delta\left(
\frac{-\beta^{-1}x_{2}+\alpha^{-1}x_{1}}{x_{0}}\right)
Y(R_{\beta}v,\beta^{-1}x_{2})Y(R_{\alpha}u,\alpha^{-1}x_{1})\nonumber\\
&=&\alpha x_{1}^{-1}\delta\left(
\frac{\beta^{-1}x_{2}+x_{0}}{\alpha^{-1}x_{1}}\right)
Y(Y(R_{\alpha}u,x_{0})R_{\beta}v,\beta^{-1}x_{2}).
\end{eqnarray}
By rewriting the three delta functions we get
\begin{eqnarray}
& &(\alpha x_{0})^{-1}\delta\left(\frac{x_{1}-\alpha\beta^{-1}x_{2}}
{\alpha x_{0}}\right)
Y(R_{\alpha}u,\alpha^{-1}x_{1})Y(R_{\beta}v,\beta^{-1}x_{2})\nonumber\\
& &-(\alpha x_{0})^{-1}\delta\left(
\frac{-\alpha\beta^{-1}x_{2}+x_{1}}{\alpha x_{0}}\right)
Y(R_{\beta}v,\beta^{-1}x_{2})Y(R_{\alpha}u,\alpha^{-1}x_{1})\nonumber\\
&=&x_{1}^{-1}\delta\left(
\frac{\alpha\beta^{-1}x_{2}+\alpha x_{0}}{x_{1}}\right)
Y(Y(R_{\alpha}u,x_{0})R_{\beta}v,\beta^{-1}x_{2}).
\end{eqnarray}
In view of (\ref{e-abstract-conjugation}), we have
$$Y(Y(R_{\alpha}u,x_{0})R_{\beta}v,\beta^{-1}x_{2})=
Y_{\beta}(Y_{\beta^{-1}\alpha}(u,\alpha x_{0})v,x_{2}).$$
Replacing $x_{0}\rightarrow \alpha^{-1}x_{0}$ we get
the desired $\Gamma$-Jacobi identity, proving that 
$(V,\{Y_{\alpha}\}, {\bf 1})$ is a $\Gamma$-vertex algebra.
\end{proof}

\br{rgammavertexalgebra} 
{\em In view of Theorem \ref{tsingle-mult},
one can define $\Gamma$-vertex algebras simply as (ordinary) vertex
algebras equipped with a compatible $\Gamma$-module structure
(cf. (\ref{e-abstract-conjugation})). }  
\er

In Theorem \ref{tsingle-mult}, if the group homomorphism 
$\phi$ from $\Gamma$ to $\C^{\times}$ is
trivial, then the second part of
(\ref{e6.25}) and (\ref{e-abstract-conjugation})
exactly amounts to that $R_{\alpha}$ is an automorphism of vertex
algebra $V$. In view of this, we immediately have:

\bc{c-trivial-case}
Let $\Gamma$ be any group equipped with
the trivial group homomorphism from $\Gamma$ to $\C^{\times}$.
Then a $\Gamma$-vertex algebra
exactly amounts to a vertex algebra equipped with
an action of $\Gamma$ by automorphisms of vertex algebra $V$.
\ec

\bd{dzgradedgammava}
{\em A $\Gamma$-vertex algebra $V$ equipped with a $\Z$-grading
$V=\coprod_{n\in \Z}V_{(n)}$ is called a {\em $\Z$-graded
$\Gamma$-vertex algebra} if ${\bf 1}\in V_{(0)}$ and 
if for any $u\in V_{(m)},\; \alpha\in
\Gamma,\; m,n,r\in \Z$,
\begin{eqnarray}
u_{(\alpha,n)}V_{(r)}\subset V_{(r+m-n-1)},
\end{eqnarray}
where $Y_{\alpha}(u,x)=\sum_{n\in \Z}u_{(\alpha,n)}x^{-n-1}$.}
\ed

The following result states that
$\Z$-graded $\Gamma$-vertex algebras can be characterized 
in terms of $\Z$-graded vertex algebras equipped with
an action of $\Gamma$ by grading-preserving automorphisms:

\bp{pgraded-case}
Let $V=\coprod_{n\in \Z}V_{(n)}$ be a $\Z$-graded vertex algebra 
in the sense that $V$ is a vertex algebra equipped with
a $\Z$-grading $V=\coprod_{n\in \Z}V_{(n)}$ such that
\begin{eqnarray}
& &{\bf 1}\in V_{(0)},\\
& &u_{m}V_{(n)}\subset V_{(n+k-m-1)}
\;\;\;\mbox{ for }u\in V_{(k)},\; k,m,n\in \Z.\label{ezgradedva}
\end{eqnarray}
Denote by $L(0)$ the degree operator, 
i.e., $L(0)v=nv$ for $v\in V_{(n)}$ with $n\in \Z$.
Let $\pi$ be a representation of $\Gamma$ on $V$ by grading-preserving
automorphisms of $V$ as a vertex algebra.
For $\alpha\in \Gamma,\; v\in V$, set
\begin{eqnarray}
Y_{\alpha}(v,x)=Y(\alpha^{-L(0)}\pi_{\alpha}v,\alpha^{-1}x) 
\in \Hom (V,V((x))).
\end{eqnarray}
Then $(V,\{Y_{\alpha}\},{\bf 1})$ is a $\Z$-graded $\Gamma$-vertex algebra.
On the other hand, let $(V,\{Y_{\alpha}\},{\bf 1})$ be 
a $\Z$-graded $\Gamma$-vertex algebra with
the degree operator denoted by $L(0)$. 
For $\alpha\in \Gamma$, set
\begin{eqnarray}
R^{o}_{\alpha}=\alpha^{L(0)}R_{\alpha}\in \Aut_{\C}V.
\end{eqnarray}
Then $R^{o}_{\alpha}$ is a grading-preserving automorphism 
of $V$ as a vertex algebra
and $R^{o}$ is a representation on $V$.
\ep

\begin{proof} With (\ref{ezgradedva}), from [FHL], 
for $z\in \C^{\times},\; v\in V$, we have
\begin{eqnarray}
z^{L(0)}Y(v,x)z^{-L(0)}=Y(z^{L(0)},z x),
\end{eqnarray}
so that
\begin{eqnarray}
\alpha^{-L(0)}\pi_{\alpha}Y(v,x)\pi_{\alpha}^{-1}\alpha^{L(0)}
=\alpha^{-L(0)}Y(\pi_{\alpha}v,x)\alpha^{L(0)}
=Y(\alpha^{-L(0)}\pi_{\alpha}v,\alpha^{-1}x).
\end{eqnarray}
With ${\bf 1}\in V_{(0)}$, we have
$\alpha^{-L(0)}\pi_{\alpha}{\bf 1}={\bf 1}$ for $\alpha\in \Gamma$.
Notice that $\pi_{\alpha}$ for $\alpha\in \Gamma$
being grading-preserving amounts to
$\alpha^{L(0)}\pi_{\beta}=\pi_{\beta}\alpha^{L(0)}$ for
$\alpha,\beta\in \Gamma$, so that
$\alpha\mapsto \alpha^{-L(0)}\pi_{\alpha}$ is a representation of
$\Gamma$ on $V$.
It follows immediately from the second part of
Theorem \ref{tsingle-mult} that $(V,\{Y_{\alpha}\},{\bf 1})$ 
is a $\Gamma$-vertex algebra. Clearly, it is a  $\Z$-graded
$\Gamma$-vertex algebra.

On the other hand, let $(V,\{Y_{\alpha}\},{\bf 1})$ 
be a $\Z$-graded $\Gamma$-vertex algebra
with $V=\coprod_{n\in \Z}V_{(n)}$.
For $\alpha\in \Gamma,\; v\in V$, we have
\begin{eqnarray}
& &R^{o}_{\alpha}Y_{1}(v,x)(R^{o}_{\alpha})^{-1}
=\alpha^{L(0)}R_{\alpha}Y_{1}(v,x)R_{\alpha}^{-1}\alpha^{-L(0)}
=Y_{1}(\alpha^{L(0)}R_{\alpha},x)=Y_{1}(R^{o}_{\alpha}v,x),\ \ \ \ \ \\
& &R^{o}_{\alpha}({\bf 1})=\alpha^{L(0)}R_{\alpha}({\bf 1})={\bf 1},
\end{eqnarray}
using the assumption ${\bf 1}\in V_{(0)}$.
Thus, for each $\alpha\in \Gamma$,
$R_{\alpha}^{o}$ is an automorphism of vertex algebra $V$.
Furthermore, since $V$ is a $\Z$-graded $\Gamma$-vertex algebra,
for $v\in V_{(m)},\; \alpha\in \Gamma$, we have
\begin{eqnarray*}
R_{\alpha}(v)=v_{(\alpha,-1)}{\bf 1}\in V_{(m)},
\end{eqnarray*}
recalling (\ref{edef-ralpha}) for the definition of $R_{\alpha}$.
Thus $R_{\alpha}$ preserves the $\Z$-grading of $V$. Consequently,
$R_{\alpha}$ commutes with $z^{L(0)}$ for any $z\in \C^{\times}$. 
Now, it follows that
$R^{o}$ is a representation of $\Gamma$ on $V$.
\end{proof}

As an immediate consequence we have:

\bc{cvoa-gammava}
Let $(V,Y, {\bf 1},\omega)$ be a vertex operator algebra and let $\pi$ be 
a group homomorphism {}from $\Gamma$ to $\Aut V$ (the automorphism group of $V$ 
as a vertex operator algebra). Then
$(V,\{Y_{\alpha}\},{\bf 1})$ is a $\Gamma$-vertex algebra, where
\begin{eqnarray}
Y_{\alpha}(v,x)=Y(\alpha^{-L(0)}\pi_{\alpha}(v),\alpha^{-1}x)
\;\;\;\mbox{ for }\alpha\in \Gamma,\; v\in V.
\end{eqnarray}
\ec

The following is an analogue of a well known result in vertex algebra theory
(cf. [LL]):

\bp{pdiff-def}
Let $V$ be a vector space equipped with linear maps $Y_{\alpha}$ from $V$
to $\Hom (V,V((x)))$ for $\alpha\in \Gamma$, a distinguished vector ${\bf 1}\in V$, 
a linear operator $\D$ on $V$ such that $\D {\bf 1}=0$
and a representation $R$ of $\Gamma$ on $V$ such that $R_{\alpha}{\bf 1}={\bf 1}$ 
for $\alpha\in \Gamma$ and such that
\begin{eqnarray}
& &Y_{\alpha}({\bf 1},x)=1,\\
& &Y_{\alpha}(v,x){\bf 1}\in V[[x]]\;\;\mbox{ and }\;
\lim_{x\rightarrow 0}Y_{\alpha}(v,x){\bf 1}=R_{\alpha}v,\label{ethm-creation}\\
& &[\D, Y_{\alpha}(v,x)]=\alpha {d\over dx}Y_{\alpha}(v,x),\label{ethm-d-bracket}\\
& &R_{\alpha}Y_{\beta}(v,x)=Y_{\alpha \beta}(v,x)R_{\alpha}\label{ethm-conjugation}
\end{eqnarray}
and for $u,v\in V,\; \alpha,\beta\in\Gamma$ there is a nonnegative integer $k$ such that
\begin{eqnarray}\label{ealpha-beta-locality}
(x_{1}-\alpha\beta^{-1}x_{2})^{k}[Y_{\alpha}(u,x_{1}),Y_{\beta}(v,x_{2})]=0.
\end{eqnarray} 
Then $V$ is a $\Gamma$-vertex algebra.
\ep

\begin{proof} First, from [Li2], [LL] (cf. [DL], [FHL]) 
$(V,Y_{1},{\bf 1})$ is a vertex algebra.
With the assumption (\ref{ethm-conjugation}), we have 
$Y_{\alpha}(v,x)=R_{\alpha}Y_{1}(v,x)R_{\alpha^{-1}}$. 
{}From (\ref{ethm-creation}) and (\ref{ethm-d-bracket}) we have
\begin{eqnarray}
& &Y_{\alpha}(u,x){\bf 1}=e^{\alpha^{-1}x\D}R_{\alpha}u,\\
& &e^{x\D}Y_{\alpha}(u,x_{0})e^{-x\D}=Y_{\alpha}(u,x_{0}+\alpha x).
\end{eqnarray}
Let $u,v\in V$ and let $k$ be a nonnegative integer such that
(\ref{ealpha-beta-locality}) holds and such that
$x^{k}Y_{\beta}(v,x)R_{\alpha}u\in V[[x]]$, so that
$$(x_{1}-\alpha\beta^{-1})^{k}Y_{\beta}(v,x_{2}-\alpha^{-1}\beta x_{1})u
=(x_{1}-\alpha\beta^{-1})^{k}Y_{\beta}(v,-\alpha^{-1}\beta x_{1}+x_{2})u.$$
Then we have
\begin{eqnarray}
& &(x_{1}-\alpha\beta^{-1}x_{2})^{k}Y_{\alpha}(u,x_{1})Y_{\beta}(v,x_{2}){\bf 1}\nonumber\\
&=&(x_{1}-\alpha\beta^{-1}x_{2})^{k}Y_{\beta}(v,x_{2})Y_{\alpha}(u,x_{1}){\bf 1}\nonumber\\
&=&(x_{1}-\alpha\beta^{-1}x_{2})^{k}Y_{\beta}(v,x_{2})e^{\alpha^{-1}x_{1}\D}R_{\alpha}u\nonumber\\
&=&e^{\alpha^{-1}x_{1}\D}(x_{1}-\alpha\beta^{-1}x_{2})^{k}Y_{\beta}(v,x_{2}-\alpha^{-1}\beta x_{1})
R_{\alpha}u\nonumber\\
&=&e^{\alpha^{-1}x_{1}\D}(x_{1}-\alpha\beta^{-1}x_{2})^{k}Y_{\beta}(v,-\alpha^{-1}\beta x_{1}+x_{2})
R_{\alpha}u.
\end{eqnarray}
Setting $x_{2}=0$ we get
$$x_{1}^{k}Y_{\alpha}(u,x_{1})R_{\beta}v
=x_{1}^{k}e^{\alpha^{-1}x_{1}\D}Y_{\beta}(v,-\alpha^{-1}\beta x_{1})R_{\alpha}u,$$
which immediately gives the following skew symmetry
\begin{eqnarray}
Y_{\alpha}(u,x_{1})R_{\beta}v=e^{\alpha^{-1}x_{1}\D}Y_{\beta}(v,-\alpha^{-1}\beta x_{1})R_{\alpha}u.
\end{eqnarray}
Taking $\beta=1$ and then using skew symmetry for $(V,Y_{1},{\bf 1})$ we get
\begin{eqnarray}
Y_{\alpha}(u,x_{1})v=e^{\alpha^{-1}x_{1}\D}Y_{1}(v,-\alpha^{-1} x_{1})R_{\alpha}u
=Y_{1}(R_{\alpha}u,\alpha^{-1}x_{1})v.
\end{eqnarray}
Now it follows from the second assertion of Theorem \ref{tsingle-mult} that
$(V,\{Y_{\alpha}\},{\bf 1})$ is a $\Gamma$-vertex algebra.
\end{proof}

{}From the proofs of Theorem \ref{tsingle-mult} and 
Proposition \ref{pgraded-case} 
we immediately have:

\bp{pvoa-gammava-module}
Let $(V,\{Y_{\alpha}\}, {\bf 1})$ be a $\Gamma$-vertex algebra 
and let $(W,\{Y^{W}_{\alpha}\})$ be a quasi module 
for $(V,\{Y_{\alpha}\}, {\bf 1})$. Then $(W,Y^{W}_{1})$ is 
a quasi module for the vertex algebra $(V,Y_{1},{\bf 1})$ 
and the following relation
holds for $\alpha\in \Gamma,\;v\in V,\; w\in W$:
\begin{eqnarray}\label{e6.55module}
Y^{W}_{\alpha}(w,x)=Y^{W}_{1}(R_{\alpha}v,\alpha^{-1}x).
\end{eqnarray}
On the other hand, let $(W,Y_{W})$ be a quasi module 
for the vertex algebra $(V,Y_{1},{\bf 1})$.
For $\alpha\in \Gamma,\; v\in V$, set
\begin{eqnarray}
Y^{W}_{\alpha}(v,x)=Y_{W}(R_{\alpha}v,\alpha^{-1}x).
\end{eqnarray}
Then $(W,\{Y^{W}_{\alpha}\})$ is a quasi module for the $\Gamma$-vertex
algebra $(V,\{Y_{\alpha}\}, {\bf 1})$.
\ep

\begin{proof} For the first part, the first assertion is clear
and the second assertion (\ref{e6.55module}) follows from the same proof of
(\ref{e5.29}). The assertion  of the second part
follows from the second part of the proof of Theorem \ref{tsingle-mult}.
\end{proof}

\section{Examples of quasi modules for vertex algebras}
In this section we give three families of examples of quasi modules for 
vertex algebras; the first involves
twisted affine Lie algebras $\hat{\g}[\sigma]$ with respect to 
an automorphism $\sigma$ of $\g$ of finite order,
the second involves certain quantum Heisenberg Lie algebras
and the third involves certain Lie algebras on quantum torus as in [BGT].

First, following [LL] we recall the $\N$-graded vertex algebras
associated with affine Lie algebras. In fact, all the three families of examples
are related to such vertex algebras.
Let $\g$ be a (possibly infinite-dimensional) Lie algebra equipped with 
a nondegenerate symmetric invariant bilinear form $\<\cdot,\cdot\>$. 
To the pair $(\g, \<\cdot,\cdot\>)$ 
we associate the {\em untwisted affine Lie algebra},
with the underlying vector space
\begin{eqnarray}
\hat{\g}=\g \otimes \C[t,t^{-1}]\oplus \C {\bf k},
\end{eqnarray}
equipped with the bracket relations
\begin{eqnarray}\label{eunteisted-affine-comm}
[a\otimes t^{m},b\otimes t^{n}]
=[a,b]\otimes t^{m+n}+m\<a,b\>\delta_{m+n,0}{\bf k}
\end{eqnarray}
for $a,b\in \g,\; m,n\in \Z$, together with the condition that
${\bf k}$ is a nonzero central element of $\hat{\g}$.
For $n\in \Z$, set
\begin{eqnarray}
\hat{\g}_{(0)}=\g\oplus \C{\bf k},\;\; \hat{\g}_{(n)}=\g\otimes t^{-n}\;\;\;\mbox{ for }n\ne 0.
\end{eqnarray}
Then $\hat{\g}=\coprod_{n\in \Z}\hat{\g}_{(n)}$ becomes a $\Z$-graded Lie algebra.
Set
\begin{eqnarray}
\hat{\g}_{(\pm) }=\coprod_{n\ge 1}\hat{\g}_{(\pm n)}=\g \otimes t^{\mp 1}\C[t^{\mp 1}],
\end{eqnarray}
which are graded Lie subalgebras of $\hat{\g}$.

For $a\in \g$, form the generating function
\begin{eqnarray}
a(x)=\sum_{n\in \Z}(a\otimes t^{n})x^{-n-1}\in \hat{\g}[[x,x^{-1}]].
\end{eqnarray}
For $a\in \g,\; n\in \Z$, we use $a(n)$ for the operator
on a $\hat{\g}$-module, corresponding to $a\otimes t^{n}$.
For a $\hat{\g}$-module $W$ and for $a\in \g$, we set
\begin{eqnarray}
a_{W}(x)=\sum_{n\in \Z}a(n)x^{-n-1}\in (\End W)[[x,x^{-1}]].
\end{eqnarray}

A $\hat{\g}$-module $W$ on which ${\bf k}$ acts as a scalar $\ell\in
\C$ is said to be of {\em level} $\ell$. 
A $\hat{\g}$-module $W$ is said to be {\em restricted} if for any
$a\in \g$ and $w\in W$, $a(n)w=0$ for $n$ sufficiently large.
Then a $\hat{\g}$-module $W$ is a restricted module if and only if
$a_{W}(x)\in \E(W)$ for $a\in \g$.

For any complex number $\ell$, let $\C_{\ell}$ be the one-dimensional
module for the Lie subalgebra $\g\otimes \C[t] \oplus\C {\bf k}$ with
$\g\otimes \C[t]$ acting trivially and with ${\bf k}$ acting as the scalar $\ell$.
Then form the induced $\hat{\g}$-module
\begin{eqnarray}
V_{\hat{\g}}(\ell,0)=U(\hat{\g})\otimes_{U(\g\otimes \C[t] \oplus\C {\bf k})} \C_{\ell},
\end{eqnarray}
which is a restricted $\hat{\g}$-module of level $\ell$ with $1\otimes 1$
as a generator. Set
\begin{eqnarray}
{\bf 1}=1\otimes 1\in V_{\hat{\g}}(\ell,0).
\end{eqnarray}
The induced $\hat{\g}$-module $V_{\hat{\g}}(\ell,0)$ is naturally an $\N$-graded module with
$\deg \C_{\ell}=0$.
In view of the Poincare-Birkhoff-Witt theorem we have
\begin{eqnarray}
V_{\hat{\g}}(\ell,0)=U(\hat{\g}_{(+)})
\end{eqnarray}
as an $\N$-graded  vector space. Identify $\g$ as a subspace of $V_{\hat{\g}}(\ell,0)$ through the linear map
$a\mapsto a(-1){\bf 1}$. 

The following results are well known (see [FZ], [Lia], cf. \cite{li-local}, [LL]
for the first assertion and see
\cite{li-local}, [LL], cf. [FZ] for the second assertion):

\bt{taffine-old}
There exists a unique vertex algebra structure on
$V_{\hat{\g}}(\ell,0)$ with ${\bf 1}$ as vacuum vector and
with $Y(a,x)=a(x)$ for $a\in \g$.
The vertex algebra $V_{\hat{\g}}(\ell,0)$ equipped with the grading
is an $\N$-graded vertex algebra 
with $V_{\hat{\g}}(\ell,0)_{(0)}=\C {\bf 1}$ and with $V_{\hat{\g}}(\ell,0)_{(1)}=\g$
generates  $V_{\hat{\g}}(\ell,0)$ as a vertex algebra. Furthermore,
on any restricted $\hat{\g}$-module $W$ of level $\ell$ there exists a 
unique $V_{\hat{\g}}(\ell,0)$-module structure 
with $Y_{W}(a,x)=a_{W}(x)$ for $a\in \g=V_{\hat{\g}}(\ell,0)_{(1)}$.
\et

Our first family of examples is about constructing quasi modules for 
$V_{\hat{\g}}(\ell,0)$ using restricted modules of level $\ell$ 
for a twisted affine Lie algebra.

Let $\sigma$ be an automorphism of $\g$ of (finite) order $T$, 
which preserves the bilinear form $\<\cdot,\cdot\>$.
Set
\begin{eqnarray}
\omega_{T}=\exp \left(\frac{2\pi \sqrt{-1}}{T}\right),
\end{eqnarray}
the principal $T$th root of unity, and set
\begin{eqnarray}
\Gamma_{T}=\{ \omega_{T}^{r}\;|\; r=0,\dots,T-1\},
\end{eqnarray}
the group of $T$th roots of unity. 
For $r\in \Z$, set
\begin{eqnarray}
\g_{r}=\{ a\in \g\;|\; \sigma (a)=\omega_{T}^{r}a\}.
\end{eqnarray}
Then 
\begin{eqnarray}
& &\g_{r}=\g_{s}\;\;\;\mbox{ if }r\equiv s \;({\rm mod}\; T),\\
& &\g=\g_{0}\oplus \g_{1}\oplus \cdots \oplus \g_{T-1}.
\end{eqnarray}

\br{rautomorphism}
{\em It is easy to see that $\sigma$ gives rise to an automorphism, also denoted by $\sigma$,
of the $\N$-graded vertex algebra $V_{\hat{\g}}(\ell,0)$.}
\er

Associated to the triple $(\g,\<\cdot,\cdot\>,\sigma)$ we have a twisted affine Lie algebra (see [K1])
\begin{eqnarray}
\hat{\g}^{\sigma}=\coprod_{i=0}^{T-1}\g_{i}\otimes
t^{i}\C[t^{T},t^{-T}]\oplus \C {\bf k},
\end{eqnarray}
which is a Lie subalgebra of the untwisted affine algebra $\hat{\g}$.
We also have the following variant of the twisted affine Lie algebra (cf. [FLM])
\begin{eqnarray}
\hat{\g}[\sigma]
=\coprod_{i=0}^{T-1}\g_{i}\otimes t^{\frac{i}{T}}\C[t,t^{-1}]\oplus \C {\bf k},
\end{eqnarray}
where the bracket is given by (\ref{eunteisted-affine-comm})
with $m,n\in \frac{1}{T}\Z$. The notions of module of level $\ell$ and restricted module
are defined in the obvious way.

\br{rtwisted-affine}
{\em It is well known (and easy to see) that the linear map defined by
\begin{eqnarray}
a\otimes t^{n}\mapsto a\otimes t^{\frac{n}{T}},\;\;\;\; T{\bf k}\mapsto {\bf k}
\end{eqnarray}
is a Lie algebra isomorphism from $\hat{\g}^{\sigma}$ to $\hat{\g}[\sigma]$. 
Consequently, a (restricted) $\hat{\g}^{\sigma}$-module of level $\ell/T$
exactly amounts to a (restricted) $\hat{\g}[\sigma]$-module of level $\ell$. }
\er

For $a\in \g$, set
\begin{eqnarray}
a(x)^{\sigma}=\sum_{i=0}^{T-1}\sigma^{i}(a)\otimes 
x^{-1}\delta\left(\frac{\omega^{-i} t}{x}\right)
=\sum_{n\in \Z}\sum_{i=0}^{T-1}\omega_{T}^{-in}\left(\sigma^{i}(a)\otimes t^{n}\right)x^{-n-1}
\in \hat{\g}[[x,x^{-1}]].
\end{eqnarray}
If $a\in \g_{s}$ with $s\in \Z$, i.e., $\sigma (a)=\omega_{T}^{s}a$, we have
\begin{eqnarray}
a(x)^{\sigma}=T\sum_{n\in \Z}(a\otimes t^{s+nT})x^{-s-nT-1}.
\end{eqnarray}
Let $a\in \g_{r},\; b\in \g_{s}$ with $r,s\in \Z$. 
Since $\<a,b\>=0$ if $r+s\notin T\Z$, we have
$$\<a,b\>\delta_{\bar{r}+\bar{s},\bar{0}}=\<a,b\>,$$
where $\bar{r}=r+T\Z\in \Z/T\Z$.
Then
\begin{eqnarray}\label{easigmabsigma-comm}
& &[a(x_{1})^{\sigma},b(x_{2})^{\sigma}]\nonumber\\
&=&T^{2}\sum_{m,n\in \Z}[a\otimes t^{r+mT},b\otimes t^{s+nt}] 
x_{1}^{-(r+mT)-1}x_{2}^{-(s+nT)-1}\nonumber\\
&=&T^{2}\sum_{m,n\in \Z}\left([a,b]\otimes t^{r+s+(m+n)T}
+(r+mT)\<a,b\>\delta_{r+s+mT+nT,0}{\bf k}\right)
x_{1}^{-(r+mT)-1}x_{2}^{-(s+nT)-1}\nonumber\\
&=&T^{2}\sum_{m,n\in \Z}[a,b]\otimes t^{r+s+(m+n)T}
x_{2}^{-(r+s+mT+nT)-1}x_{1}^{-(r+mT)-1}x_{2}^{r+mT}\nonumber\\
& &+T^{2}\sum_{m\in \Z}(r+mT)\<a,b\>\delta_{\bar{r}+\bar{s},\bar{0}}{\bf k}
x_{1}^{-(r+mT)-1}x_{2}^{r+mT-1}\nonumber\\
&=&T[a,b](x_{2})^{\sigma}x_{1}^{-1}\delta\left(\left(\frac{x_{2}}{x_{1}}\right)^{T}\right)
\left(\frac{x_{2}}{x_{1}}\right)^{r}
+T^{2}\<a,b\> {\bf k}\frac{\partial}{\partial x_{2}}
x_{1}^{-1}\delta\left(\left(\frac{x_{2}}{x_{1}}\right)^{T}\right)
\left(\frac{x_{2}}{x_{1}}\right)^{r}\nonumber\\
&=&\sum_{i=0}^{T-1}\left([a,b](x_{2})^{\sigma}\omega_{T}^{-ri}
x_{1}^{-1}\delta\left(\frac{\omega_{T}^{i}x_{2}}{x_{1}}\right)
+\<a,b\> T{\bf k}\omega_{T}^{-ri}\frac{\partial}{\partial x_{2}}
x_{1}^{-1}\delta\left(\frac{\omega_{T}^{i}x_{2}}{x_{1}}\right)\right),
\end{eqnarray}
where we are using the fact
\begin{eqnarray}
Tx_{1}^{-1}\delta\left(\left(\frac{x_{2}}{x_{1}}\right)^{T}\right)
\left(\frac{x_{2}}{x_{1}}\right)^{r}
=\sum_{i=0}^{T-1}x_{1}^{-1}\delta\left(\frac{\omega_{T}^{i}x_{2}}{x_{1}}\right)
\left(\frac{x_{2}}{x_{1}}\right)^{r}
=\sum_{i=0}^{T-1}\omega_{T}^{-ri}x_{1}^{-1}\delta\left(\frac{\omega_{T}^{i}x_{2}}{x_{1}}\right).
\end{eqnarray}
{}From this we have
\begin{eqnarray*}
\left(\prod_{i=0}^{T-1}(x_{1}-\omega_{T}^{i}x_{2})^{2}\right)
[a(x_{1})^{\sigma},b(x_{2})^{\sigma}]=0.
\end{eqnarray*}
That is,
\begin{eqnarray}\label{eab-locality-relation}
(x_{1}^{T}-x_{2}^{T})^{2}[a(x_{1})^{\sigma},b(x_{2})^{\sigma}]=0.
\end{eqnarray}

Let $W$ be a restricted $\hat{\g}^{\sigma}$-module of level $\ell/T$, or equivalently,
a restricted $\hat{\g}[\sigma]$-module of level $\ell$. For $a\in \g$,
$a(x)^{\sigma}$ acts on $W$ giving rise to an element of $\E(W)$,
which we denote by $a_{W}(x)^{\sigma}$, Set
\begin{eqnarray}
S_{W}=\{ a_{W}(x)^{\sigma}\;|\; a\in \g\}\subset \E(W).
\end{eqnarray}
{}From (\ref{eab-locality-relation}), $S_{W}$ is a $\Gamma_{T}$-local subspace.
For $\alpha \in \Gamma_{T}$, we have
\begin{eqnarray}
R_{\alpha}a_{W}(x)^{\sigma}=T\sum_{n\in \Z}(a\otimes t^{s+nT})(\alpha x)^{-s-nT-1}
=\alpha^{-s-1} a_{W}(x)^{\sigma}.
\end{eqnarray}
Then $S_{W}$ is also stale under the action of $\Gamma_{T}$.
By Theorems \ref{tkey-main} and  \ref{tvertex-algebra-case}
$S_{W}$ generates a vertex algebra 
$V_{W}=\<S_{W}\>_{\Gamma_{T}}$ with $W$ 
as a quasi module. 
By Lemma \ref{lfor-application}, $V_{W}$ as a vertex algebra is generated by $S_{W}$.
In view of (\ref{easigmabsigma-comm}) and 
Proposition \ref{pdecomposition} we have
\begin{eqnarray}
& &a_{W}(x)^{\sigma}_{0}b_{W}(x)^{\sigma}=[a,b]_{W}(x)^{\sigma},\\
& &a_{W}(x)^{\sigma}_{1}b_{W}(x)^{\sigma}=\ell \<a,b\> 1_{W},\\
& &a_{W}(x)^{\sigma}_{n}b_{W}(x)^{\sigma}=0\;\;\;\mbox{ for }n\ge 2.
\end{eqnarray}
It follows (cf. [LL], Remark 6.2.17) that $V_{W}$ is a
$\hat{\g}$-module of level $\ell$ with
$a(x_{1})$ acting as $\Y_{1}(a_{W}(x)^{\sigma},x_{1})$ for $a\in \g$.
Since $S_{W}$ generates $V_{W}$ as a vertex algebra,
$V_{W}$ is a $\hat{\g}$-module generated by $1_{W}$.
It follows that there exists a $\hat{\g}$-module homomorphism
$\psi$ from $V_{\hat{\g}}(\ell,0)$ onto $V_{W}$, extending
the canonical map from $\g$ to $S_{W}$. With $\g$ as a generating subspace
of $V_{\hat{\g}}(\ell,0)$, by Theorem \ref{tmodule-extension}
$\psi$ is a vertex algebra homomorphism and
$W$ is naturally a quasi module for $V_{\hat{\g}}(\ell,0)$.
Summarizing these we have: 

\bp{ptwisted-module} 
Let $\ell$ be any complex number and let $W$ be
any restricted $\hat{\g}[\sigma]$-module of level $\ell$, or
equivalently, a restricted $\hat{\g}^{\sigma}$-module of level
$\ell/T$. Then there exists a unique quasi module structure $Y_{W}$ on
$W$ for $V_{\hat{\g}}(\ell,0)$ such that
$Y_{W}(a,x)=a_{W}(x)^{\sigma}$ for $a\in \g$, with $p(x)=x^{T}-1$ as the
minimal polynomial.  
\ep

\br{rtwisted-quasi-modules} 
{\em It was proved in \cite{li-twisted}
that giving a restricted $\hat{\g}[\sigma]$-module structure of level
$\ell$ on a vector space $W$ is canonically equivalent to
giving $\sigma$-twisted module structure on $W$ for the
vertex algebra $V_{\hat{\g}}(\ell,0)$. In view of this and Proposition
\ref{ptwisted-module}, one naturally expects that for a general vertex
operator algebra $V$, there should exist a canonical
connection between twisted $V$-modules and quasi
$V$-modules.  Indeed, such a connection exists and it will be
given in a sequel.}
\er

\br{rduality}
{\em Let $W$ be a restricted $\hat{\g}$-module of level $\ell$.
With $\hat{\g}^{\sigma}$ as a subalgebra by definition,
$W$ is naturally a restricted $\hat{\g}^{\sigma}$-module of level $\ell$.
Then $W$ is naturally a quasi module for $V_{\hat{\g}}(\ell T,0)$.}
\er

Now let us play a slightly different game.
Fix a positive integer $k$.
For $a\in \g,\; r\in \Z$, set
\begin{eqnarray}
E(a,r,x)=k\sum_{n\in \Z}(a\otimes t^{r+nk})x^{-r-nk-1}\in \hat{\g}[[x,x^{-1}]].
\end{eqnarray}
Clearly, if $r\equiv s$ (mod $k$) we have
\begin{eqnarray}
E(a,r,x)=E(a,s,x).
\end{eqnarray}

For $a,b\in \g,\; r,s\in \Z$, using (\ref{eunteisted-affine-comm}) we have
\begin{eqnarray}\label{efull-subalgebra-commutator}
& &[E(a,r,x_{1}),E(b,s,x_{2})]\nonumber\\
&=&k^{2}\sum_{m,n\in \Z}[a\otimes t^{r+mk},b\otimes t^{s+nk}]
x_{1}^{-r-mk-1}x_{2}^{-s-nk-1}\nonumber\\
&=&k^{2}\sum_{m,n\in \Z}\left( [a,b]\otimes t^{r+s+(m+n)k}
+(r+mk)\<a,b\>{\bf k}\delta_{r+mk+s+nk,0}\right)
x_{1}^{-r-mk-1}x_{2}^{-s-nk-1}\nonumber\\
&=&\sum_{i=0}^{k-1}\left(E([a,b],r+s,x_{2})\omega_{k}^{-ri}
x_{1}^{-1}\delta\left(\frac{\omega_{k}^{i}x_{2}}{x_{1}}\right)
+\<a,b\> \delta_{\bar{r}+\bar{s},\bar{0}}k{\bf k}\omega_{k}^{-ri}
\frac{\partial}{\partial x_{2}}
x_{1}^{-1}\delta\left(\frac{\omega_{k}^{i}x_{2}}{x_{1}}\right)\right).\nonumber\\
\end{eqnarray}
Let $W$ be any restricted $\hat{\g}$-module of level $\ell$.
Then $E(a,r,x)\in \E(W)$ for $a\in \g,\; r\in \Z$.
Set
\begin{eqnarray}
E_{W}={\rm span}\{ E(a,r,x)\;|\;a\in \g,\; r\in \Z\}\subset \E(W).
\end{eqnarray}
With (\ref{efull-subalgebra-commutator}), we see that
$E_{W}$ is a $\Gamma_{k}$-local subspace of $\E(W)$.
For $\alpha\in \Gamma_{k}$, we have
\begin{eqnarray}
R_{\alpha}E(a,r,x)=E(a,r,\alpha x)=\alpha^{-r-1}E(a,r,x)
\;\;\;\mbox{ for }r\in \Z,\; a\in \g.
\end{eqnarray}
Then $E_{W}$ is closed under the action of $\Gamma_{k}$.
By Theorems \ref{tkey-main} and  \ref{tvertex-algebra-case}
$E_{W}$ generates a vertex algebra 
$V_{W}=\<E_{W}\>_{\Gamma_{k}}$ with $W$ 
as a quasi module. 
By Lemma \ref{lfor-application}, $V_{W}$ as a vertex algebra 
is generated by $E_{W}$.
Combining (\ref{efull-subalgebra-commutator}) 
with Proposition \ref{pdecomposition} we have
\begin{eqnarray}
& &E(a,r,x)_{0}E(b,s,x)=E([a,b],r+s,x),\\
& &E(a,r,x)_{1}E(b,s,x)=\delta_{\bar{r}+\bar{s},\bar{0}}k\<a,b\>\ell,\\
& &E(a,r,x)_{n}E(b,s,x)=0\;\;\;\mbox{ for }n\ge 2.
\end{eqnarray}

We next relate $V_{W}$ to a vertex algebra associated with an affine
Lie algebra.
Notice that $\C[\Z/k\Z]$ is a commutative associative algebra
with a symmetric invariant bilinear form defined by
$\<e^{\bar{r}},e^{\bar{s}} \>=\delta_{\bar{r}+\bar{s},\bar{0}}$ 
for $r,s\in \Z$.
Then $\g\otimes \C[\Z/k\Z]$ is naturally a Lie algebra
which is equipped with the symmetric invariant bilinear form defined by
\begin{eqnarray}
\< a\otimes e^{\bar{r}},b\otimes e^{\bar{s}}\>
=\delta_{\bar{r}+\bar{s},\bar{0}}\<a,b\>
\;\;\;\mbox{ for }r,s\in \Z,\; a,b\in \g.
\end{eqnarray}
Associated with the pair $(\g[\Z/k\Z],\<\cdot,\cdot\>)$ 
we have the following untwisted affine Lie algebra
\begin{eqnarray}
\widehat{\g[\Z/k\Z]}
=\left(\g\otimes \C[\Z/k\Z]\right) \otimes \C[t,t^{-1}]\oplus 
\C {\bf k}.
\end{eqnarray}
Using the same arguments for Proposition \ref{ptwisted-module}
we immediately have:

\bp{pexample2}
Let $W$ be an (irreducible) restricted $\hat{\g}$-module of level $\ell$ and 
$k$ a positive integer.
For $r\in \Z,\; a\in \g$, set
\begin{eqnarray}
E(a,r,x)=\sum_{n\in \Z}a(r+nk)x^{-r-nk-1}\in \E(W).
\end{eqnarray}
There exists a unique (irreducible) quasi module structure on $W$
for the vertex algebra
$V_{\widehat{\g[\Z/k\Z]}}(k\ell,0)$ with 
$Y_{W}(a\otimes e^{\bar{r}},x)=E(a,r,x)$ for $a\in \g,\; r\in \Z$,
with $x^{k}-1$ as the minimal polynomial.
\ep

\br{rsymmetric-orbifold}
{\em Since $\C[\Z/k\Z]=\C^{\otimes k}$ as a commutative associative algebra,
we have that $\g[\Z/k\Z]=\g^{\otimes k}$.
(One can find an explicit isomorphism.)
Using the bilinear form one can see that
$V_{\widehat{\g[\Z/k\Z]}}(k\ell,0)=V_{\hat{\g}}(\ell,0)^{\otimes k}$.
In view of the first example, quasi modules 
for $V_{\hat{\g}}(\ell,0)^{\otimes k}$
with minimal polynomial $x^{k}-1$ are closely related twisted modules.
Then $V_{\hat{\g}}(\ell,0)$-modules are closely related
to twisted modules for $V_{\hat{\g}}(\ell,0)^{\otimes k}$.
This agrees with a general conceptual result obtained in [BDM]. }
\er

Next, we present our second family of examples, which
involves certain ``quantum'' Heisenberg Lie algebras.
Let $\h$ be a vector space equipped with a nondegenerate symmetric bilinear 
form $\<\cdot,\cdot\>$ and
let $q$ be any nonzero complex number. Consider the following  Lie algebra 
\begin{eqnarray}\label{equantum-heisenberg-space}
\hat{\h}_{q}=\h\otimes \C[t,t^{-1}]\oplus \C c
\end{eqnarray}
with the following commutation relation:
\begin{eqnarray}
& &[c, \hat{\h}_{q}]=0,\\
& &[a\otimes t^{m},b\otimes t^{n}]=m\<a,b\>\delta_{m+n,0}(q^{m}+q^{-m})c
\label{equantum-heisenberg-relation}
\end{eqnarray}
for $a,b\in \h,\;m,n\in \Z$.
As with affine Lie algebras, for $a\in \h,\; m\in \Z$,
we use $a(m)$ to represent the corresponding operator of $a\otimes t^{m}$ 
on any $\hat{\h}_{q}$-module. The notions of level and restricted module
are defined in the obvious ways.

For $a\in \h$, form the generating function
\begin{eqnarray}
a(x)=\sum_{n\in \Z}a(m)x^{-m-1}.
\end{eqnarray}
In terms of the generating functions we have
\begin{eqnarray}
[a(x_{1}),b(x_{2})]=\<a,b\>\frac{\partial}{\partial x_{2}}
\left( x_{1}^{-1}\delta\left(\frac{qx_{2}}{x_{1}}\right)
+x_{1}^{-1}\delta\left(\frac{q^{-1}x_{2}}{x_{1}}\right)\right)c.
\end{eqnarray}
Consequently,
\begin{eqnarray}
(x_{1}-qx_{2})^{2}(x_{1}-q^{-1}x_{2})^{2}[a(x_{1}),b(x_{2})]=0.
\end{eqnarray}
For $a\in \h,\; m\in \Z$, set
\begin{eqnarray}
a(m,x)=q^{m}a(q^{m}x)\in \hat{\h}_{q}[[x,x^{-1}]].
\end{eqnarray}
For $a,b\in \h,\; m,n\in \Z$, we have
\begin{eqnarray}\label{eambn-comm}
[a(m,x_{1}),b(n,x_{2})]&=&q^{m+n}[a(q^{m}x_{1}),b(q^{n}x_{2})]\nonumber\\
&=&\<a,b\>  \frac{\partial}{\partial x_{2}}
\left(x_{1}^{-1}\delta\left(\frac{q^{n+1}x_{2}}{q^{m}x_{1}}\right)
+x_{1}^{-1}\delta\left(\frac{q^{n-1}x_{2}}{q^{m}x_{1}}\right)\right)c.
\end{eqnarray}
Set
\begin{eqnarray}
\Gamma_{(q)}=\{ q^{m}\;|\; m\in \Z\}\subset \C^{\times}.
\end{eqnarray}
Let $W$ be any restricted $\hat{\h}_{q}$-module of level $\ell\in \C$.
For $a\in \h,\; m\in \Z$, $a(m,x)$ acting on $W$ gives rise to
an element of $\E(W)$, which we denote by $a_{W}(m,x)$.
Set
\begin{eqnarray}
U_{W}={\rm span}\{ a_{W}(m,x)\;|\; a\in \h,\; m\in \Z\}\subset \E(W).
\end{eqnarray}
{}From (\ref{eambn-comm}) we see that $U_{W}$ is a $\Gamma_{(q)}$-local 
subspace of $\E(W)$.
Furthermore, for $a\in \g,\; m,n\in \Z,$ we have
\begin{eqnarray}
R_{q^{n}}a_{W}(m,x)=a_{W}(m,q^{n}x)=q^{m}a_{W}(q^{m+n}x)=q^{-n}a_{W}(m+n,x).
\end{eqnarray}
Thus $U_{W}$ is closed under the action of $\Gamma_{(q)}$.
In view of Theorems \ref{tkey-main} and  \ref{tvertex-algebra-case}
and Lemma \ref{lfor-application}, 
$U_{W}$ generates a vertex algebra $V_{W}$ with $W$ as a
quasi module. 

Assume that $q$ is not a root of unity, so that
\begin{eqnarray}
q^{m}\ne q^{n}\;\;\;\mbox{ whenever }m\ne n.
\end{eqnarray}
In view of (\ref{eambn-comm}) and Proposition \ref{pdecomposition},
the following relations hold for $a,b\in \h,\; m,n\in \Z$:
\begin{eqnarray}
& &a_{W}(m,x)_{1}b_{W}(n,x)=\<a,b\>\ell (\delta_{m,n+1}+\delta_{m,n-1}){\bf 1}\\
& &a_{W}(m,x)_{r}b_{W}(n,x)=0\;\;\;\mbox{ for }r\ge 0,\; r\ne 1.
\end{eqnarray}

On the other hand, consider the space $\h[\Z]=\h\otimes \C[\Z]$ 
equipped with the symmetric bilinear form $\<\cdot,\cdot\>$ defined by
\begin{eqnarray}\label{equantum-form}
\<a\otimes e^{m},b\otimes e^{n}\>=\<a,b\> (\delta_{m,n+1}+\delta_{n,m+1})
\;\;\;\mbox{ for }a,b\in \h,\; m,n\in \Z.
\end{eqnarray}
Then consider the affine Lie algebra $\widehat{\h[\Z]}$
associated with the pair $(\h[\Z],\<\cdot,\cdot\>)$, where 
$\h[\Z]$ is viewed as an abelian Lie algebra. Let $V_{\widehat{\h[\Z]}}(\ell,0)$ 
be the vertex algebra
associated with the affine Lie algebra $\widehat{\h[\Z]}$ of level $\ell$. 
This is an $\N$-graded 
vertex algebra with $V_{\widehat{\h[\Z]}}(\ell,0)_{(0)}=\C {\bf 1}$ and
$V_{\widehat{\h[\Z]}}(\ell,0)_{(1)}=\h[\Z]$ as its generating subspace.
Using the same arguments as we used for Proposition \ref{ptwisted-module} 
we immediately have:

\bp{pquantum-Heisenberg}
Let $\h$ be a vector space equipped with a symmetric bilinear form $\<\cdot,\cdot\>$,
let $q$ be a nonzero complex number which is not a root of unity and
let $\hat{\h}_{q}$ be the quantum Lie algebra defined in
(\ref{equantum-heisenberg-space})--(\ref{equantum-heisenberg-relation}).
Let $\h[\Z]$ be the abelian Lie algebra equipped with the symmetric bilinear form
defined in (\ref{equantum-form})
and $\widehat{\h[\Z]}$ be the corresponding affine Lie algebra.
Then for any restricted $\hat{\h}_{q}$-module $W$ of level $\ell$, there exists a unique
quasi $V_{\widehat{\h[\Z]}}(\ell,0)$-module structure on $W$ 
with $Y_{W}(a\otimes e^{m},x)=a_{W}(m,x)$ for $a\in \g,\; m\in \Z$.
\ep

\br{rroot-unity}
{\em If $q$ is a root of unity, one can appropriately modify
the current analysis to work out the corresponding results.}
\er

Finally, we present our third family of examples. This family involves
certain Lie algebras of quantum torus which have been studied in [BGT] and [G-K-L].
First, following [BGT] we present the Lie algebras.

Let $q$ be a nonzero complex number. Consider the following twisted group
algebra of $\Z^{2}$
\begin{eqnarray}
\C_{q}[t_{0}^{\pm 1},t_{1}^{\pm 1}]=\C[t_{0}^{\pm 1},t_{1}^{\pm 1}]
\end{eqnarray}
as a vector space, where
\begin{eqnarray}
t_{1}t_{0}=qt_{0}t_{1}.
\end{eqnarray}
For $m,n,r,s\in \Z$, we have
\begin{eqnarray}
(t_{0}^{m}t_{1}^{n})(t_{0}^{r}t_{1}^{s})=q^{nr}t_{0}^{m+r}t_{1}^{n+s}.
\end{eqnarray}
Let $A$ be an associative algebra equipped with a nondegenerate 
symmetric invariant bilinear form $\<\cdot,\cdot\>$ in the sense that
\begin{eqnarray}
\<ab,c\>=\<a,bc\>\;\;\;\mbox{ for }a,b,c\in A. 
\end{eqnarray}
Consider the tensor product associative algebra
\begin{eqnarray}
A_{q}[t_{0}^{\pm 1},t_{1}^{\pm 1}]
=A\otimes \C_{q}[t_{0}^{\pm 1},t_{1}^{\pm 1}].
\end{eqnarray}
Naturally, $A_{q}[t_{0}^{\pm 1},t_{1}^{\pm 1}]$ 
is a Lie algebra with the commutator map as the Lie
bracket, which is denoted by $[\cdot,\cdot]_{loop}$.
For $a,b\in A,\; m,n,r,s\in \Z$, we have
\begin{eqnarray}
[a\otimes t_{0}^{m}t_{1}^{n}, b\otimes t_{0}^{r}t_{1}^{s}]_{loop}
=q^{nr}(ab\otimes t_{0}^{m+r}t_{1}^{n+s})
-q^{ms}(ba\otimes t_{0}^{m+r}t_{1}^{n+s}).
\end{eqnarray}
We have the following two-dimensional central extension of the Lie
algebra $A_{q}[t_{0}^{\pm 1},t_{1}^{\pm 1}]$:
\begin{eqnarray}\label{eEALA}
\hat{A_{q}}[t_{0}^{\pm 1},t_{1}^{\pm 1}]
=A_{q}[t_{0}^{\pm 1},t_{1}^{\pm 1}]\oplus \C c_{0}\oplus \C c_{1}=A\otimes
\C_{q}[t_{0}^{\pm 1},t_{1}^{\pm 1}]\oplus \C c_{0}\oplus \C c_{1},
\end{eqnarray}
where $c_{0},c_{1}$ are central and
\begin{eqnarray}
[a\otimes t_{0}^{m}t_{1}^{n}, b\otimes t_{0}^{r}t_{1}^{s}]
=[a\otimes t_{0}^{m}t_{1}^{n}, b\otimes t_{0}^{r}t_{1}^{s}]_{loop}
+\<a,b\>q^{nr}\delta_{m+r,0}\delta_{n+s,0}(mc_{0}+nc_{1}).
\end{eqnarray}
Thus
\begin{eqnarray}
& &[a\otimes t_{0}^{m}t_{1}^{n}, b\otimes t_{0}^{r}t_{1}^{s}]\nonumber\\
&=&q^{nr}(ab\otimes t_{0}^{m+r}t_{1}^{n+s})
-q^{ms}(ba\otimes t_{0}^{m+r}t_{1}^{n+s})
+\<a,b\>q^{nr}\delta_{m+r,0}\delta_{n+s,0}(mc_{0}+nc_{1}).
\end{eqnarray}

For $a\in A,\; n\in \Z$, set
\begin{eqnarray}
X(a,n,x)=\sum_{m\in \Z}(a\otimes t_{0}^{m}t_{1}^{n})x^{-m-1}
\in \hat{A_{q}}[t_{0}^{\pm 1},t_{1}^{\pm 1}][[x,x^{-1}]].
\end{eqnarray}
Then
\begin{eqnarray}
& &[X(a,n,x_{1}), X(b,s,x_{2})]\nonumber\\
&=&\sum_{m,r\in \Z}[a\otimes t_{0}^{m}t_{1}^{n},b\otimes
t_{0}^{r}t_{1}^{s}] x_{1}^{-m-1}x_{2}^{-r-1}\nonumber\\
&=&\sum_{m,r\in \Z}\left(q^{nr}(ab\otimes t_{0}^{m+r}t_{1}^{n+s})
-q^{ms}(ba\otimes t_{0}^{m+r}t_{1}^{n+s})\right)x_{1}^{-m-1}x_{2}^{-r-1}
\nonumber\\
& &+\sum_{m,r\in \Z} \<a,b\>q^{nr}\delta_{m+r,0}
\delta_{n+s,0}(mc_{0}+nc_{1})x_{1}^{-m-1}x_{2}^{-r-1}\nonumber\\
&=&\sum_{m,r\in \Z}(ab\otimes t_{0}^{m+r}t_{1}^{n+s})
(q^{-n}x_{2})^{-m-r-1}x_{1}^{-m-1}(q^{-n}x_{2})^{m}q^{-n}\nonumber\\
& &-\sum_{m,r\in \Z}(ba\otimes t_{0}^{m+r}t_{1}^{n+s})x_{2}^{-m-r-1}
x_{1}^{-m-1}(q^{s}x_{2})^{m}\nonumber\\
& &+\sum_{m\in \Z}\<a,b\> \delta_{n+s,0}
q^{-nm}(mc_{0}+nc_{1})x_{1}^{-m-1}x_{2}^{m-1}
\nonumber\\
&=&q^{-n}X(ab,n+s,q^{-n}x_{2})
x_{1}^{-1}\delta\left(\frac{q^{-n}x_{2}}{x_{1}}\right)
-X(ba,n+s,x_{2})x_{1}^{-1}\delta\left(\frac{q^{s}x_{2}}{x_{1}}\right)
\nonumber\\
& &+\<a,b\>\delta_{n+s,0}\frac{\partial}{\partial x_{2}}
x_{1}^{-1}\delta\left(\frac{q^{-n}x_{2}}{x_{1}}\right)c_{0}
+\<a,b\>\delta_{n+s,0}n \left(c_{1}x_{2}^{-1}\right)
x_{1}^{-1}\delta\left(\frac{q^{-n}x_{2}}{x_{1}}\right).
\end{eqnarray}

Furthermore, for $a\in A,\; m,n\in \Z$, we set
\begin{eqnarray}
\bar{X}(a,m,n,x)=q^{n}X(a,m,q^{n}x).
\end{eqnarray}
Then
\begin{eqnarray}\label{e-quantum-torus-comm}
& &[\bar{X}(a,n,m,x_{1}),\bar{X}(b,s,r,x_{2})]\nonumber\\
&=&q^{m+r}[X(a,n, q^{m}x_{1}),X(b,s,q^{r}x_{2})]\nonumber\\
&=&q^{m+r}q^{-n}X(ab,n+s,q^{-n+r}x_{2})q^{-m}
x_{1}^{-1}\delta\left(\frac{q^{-n+r}x_{2}}{q^{m}x_{1}}\right)
\nonumber\\
& &-q^{m+r}X(ba,n+s,q^{r}x_{2})q^{-m}
x_{1}^{-1}\delta\left(\frac{q^{s+r}x_{2}}{q^{m}x_{1}}\right)\nonumber\\
& &+q^{m+r}\<a,b\>\delta_{n+s,0}c_{0}q^{-r}\frac{\partial}{\partial x_{2}}
q^{-m}x_{1}^{-1}\delta\left(\frac{q^{-n+r}x_{2}}{q^{m}x_{1}}\right)
\nonumber\\
& &+q^{m+r}\<a,b\>\delta_{n+s,0}c_{1}q^{-r}x_{2}^{-1}\cdot 
q^{-m}x_{1}^{-1}\delta\left(\frac{q^{-n+r}x_{2}}{q^{m}x_{1}}\right)
\nonumber\\
&=&q^{r-n}X(ab,n+s,q^{-n+r}x_{2})
x_{1}^{-1}\delta\left(\frac{q^{r-n}x_{2}}{q^{m}x_{1}}\right)
-q^{r}X(ba,n+s,q^{r}x_{2})
x_{1}^{-1}\delta\left(\frac{q^{s+r}x_{2}}{q^{m}x_{1}}\right)\nonumber\\
& &+\<a,b\>\delta_{n+s,0}c_{0}\frac{\partial}{\partial x_{2}}
x_{1}^{-1}\delta\left(\frac{q^{r-n}x_{2}}{q^{m}x_{1}}\right)
+\<a,b\>\delta_{n+s,0}\left(c_{1}x_{2}^{-1}\right)
x_{1}^{-1}\delta\left(\frac{q^{r-n}x_{2}}{q^{m}x_{1}}\right)
\nonumber\\
&=&\bar{X}(ab,n+s,r-n,x_{2})
x_{1}^{-1}\delta\left(\frac{q^{r-n}x_{2}}{q^{m}x_{1}}\right)
-\bar{X}(ba,n+s,r,x_{2})
x_{1}^{-1}\delta\left(\frac{q^{s+r}x_{2}}{q^{m}x_{1}}\right)\nonumber\\
& &+\<a,b\>\delta_{n+s,0}\left(c_{1}x_{2}^{-1}\right)
x_{1}^{-1}\delta\left(\frac{q^{r-n}x_{2}}{q^{m}x_{1}}\right)
+\<a,b\>\delta_{n+s,0}c_{0}\frac{\partial}{\partial x_{2}}
x_{1}^{-1}\delta\left(\frac{q^{r-n}x_{2}}{q^{m}x_{1}}\right).
\end{eqnarray}
Set
\begin{eqnarray}
c_{1}(x)=c_{1}x^{-1}.
\end{eqnarray}

For convenience let us just use $\hat{A_{q}}$ for the Lie algebra 
$\hat{A_{q}}[t_{0}^{\pm 1},t_{1}^{\pm 1}]$.
An $\hat{A_{q}}$-module 
on which $c_{0}$ acts as a scalar $\ell\in \C$ 
is said to be {\em of level} $\ell$.
An $\hat{A_{q}}$-module $W$ is said to be {\em restricted} if
for any $a\in A,\; n\in \Z$ and for any $w\in W$,
$X(a,n,x)w\in W((x))$. That is, $X(a,n,x)$ acting on $W$ is an element
of $\E(W)$. Denote by $X_{W}(a,n,x)$ the corresponding element of
$\E(W)$ and similarly for $\bar{X}_{W}(a,n,m,x)$.

Let $W$ be a restricted $\hat{A_{q}}$-module of level $\ell$.
Set
\begin{eqnarray}
U_{W}={\rm span }\{ c_{1}(x),\;\bar{X}_{W}(a,n,m,x)\;|\; a\in A,\; m,n\in \Z\}
\subset \E(W).
\end{eqnarray}
In view of (\ref{e-quantum-torus-comm}), $U_{W}$ is a $\Gamma_{(q)}$-local
subspace of $\E(W)$, where 
\begin{eqnarray*}
\Gamma_{(q)}=\{ q^{m}\;|\; m\in \Z\}\subset \C^{\times}.
\end{eqnarray*}
We also have
\begin{eqnarray}
& &R_{q^{n}}c_{1}(x)=c_{1}(x),\\
& &R_{q^{r}}\bar{X}_{W}(a,n,m,x)=\bar{X}_{W}(a,n,m,q^{r}x)
=q^{m}X_{W}(a,n,q^{m+r}x)\nonumber\\
& &=q^{-r}\bar{X}_{W}(a,n,m+r,x)
\end{eqnarray}
for $a\in A,\; m,n,r\in \Z$. Then $U_{W}$ is stable under the action of
of $\Gamma_{(q)}$.
In view of Theorem \ref{tmain1} and Lemma \ref{lfor-application}, 
$U_{W}$ generates a vertex algebra $V_{W}$ with $W$ as a
quasi module. 

Next we determine the vertex algebra structure of $V_{W}$.
The fact is that the structure of $V_{W}$depends on
if $q$ is a root of unity.
Here we shall just consider the case that $q$ is not a root of unity. That is,
for any $m,n\in \Z$, $q^{m}=q^{n}$ if and only if $m=n$.
In view of (\ref{e-quantum-torus-comm}) and Proposition \ref{pdecomposition}
we have the following relations in $V_{W}$ 
for $a,b\in A,\; m,n,r,s\in \Z$:
\begin{eqnarray}
& &\bar{X}_{W}(a,n,m,x)_{0}\bar{X}(b,s,r,x)
=\delta_{m,r-n}\bar{X}_{W}(ab,n+s,r-n,x)\nonumber\\
& &\hspace{2cm}-\delta_{m,s+r}\bar{X}_{W}(ba,n+s,r,x) 
+\<a,b\>\delta_{n+s,0}\delta_{m,r-n}c_{1}(x),\label{e7.68} \\
& &\bar{X}_{W}(a,n,m,x)_{1}\bar{X}(b,s,r,x)
=\<a,b\>\delta_{n+s,0}\delta_{m,r-n}\ell,\\
& &\bar{X}_{W}(a,n,m,x)_{k}\bar{X}(b,s,r,x)
=0\hspace{2cm}\mbox{for }k\ge 2.\label{e7.70}
\end{eqnarray}

To determine the vertex algebra $V_{W}$,
first we construct an affine Lie algebra.

\bl{lLie-algebra}
We define a non-associative algebra with the underlying vector space
\begin{eqnarray}
\C_{*}[t_{0}^{\pm 1},t_{1}^{\pm 1}]
=\C[t_{0}^{\pm 1},t_{1}^{\pm 1}],
\end{eqnarray}
equipped with the multiplication
\begin{eqnarray}\label{enonassoc-product}
(t_{0}^{n}t_{1}^{m})(t_{0}^{s}t_{1}^{r})
=\delta_{n+m,r}t_{0}^{n+s}t_{1}^{m}
\end{eqnarray}
for $m,n,r,s\in \Z$. Endow this non-associative algebra with a bilinear form defined by
\begin{eqnarray}\label{e-nonassoc-form}
\< t_{0}^{n}t_{1}^{m},t_{0}^{s}t_{1}^{r}\>=\delta_{n+s,0}\delta_{m+n,r}
\end{eqnarray}
for $m,n,r,s\in \Z$. Then
this algebra is associative and the bilinear form is 
symmetric and associative.
\el

\begin{proof} This algebra is associative as
\begin{eqnarray*}
& &(t_{0}^{k}t_{1}^{l})\left((t_{0}^{n}t_{1}^{m})(t_{0}^{s}t_{1}^{r})\right)
=\delta_{n+m,r}(t_{0}^{k}t_{1}^{l})(t_{0}^{n+s}t_{1}^{m})
=\delta_{n+m,r}\delta_{k+l,m}t_{0}^{k+n+s}t_{1}^{l},\\
& &\left((t_{0}^{k}t_{1}^{l})(t_{0}^{n}t_{1}^{m})\right)
(t_{0}^{s}t_{1}^{r})
=\delta_{k+l,m}(t_{0}^{k+n}t_{1}^{l})(t_{0}^{s}t_{1}^{r})
=\delta_{k+l,m}\delta_{k+n+l,r}t_{0}^{k+n+s}t_{1}^{l}\\
& &\hspace{2cm}=\delta_{k+l,m}\delta_{m+n,r}t_{0}^{k+n+s}t_{1}^{l}.
\end{eqnarray*}
The bilinear form is symmetric since
$$\delta_{n+s,0}\delta_{m+n,r}=\delta_{n+s,0}\delta_{r-n,m}
=\delta_{n+s,0}\delta_{r+s,m}$$
and it is invariant because
\begin{eqnarray*}
\< t_{0}^{n}t_{1}^{m},t_{0}^{s}t_{1}^{r}\>
=\phi \left( (t_{0}^{n}t_{1}^{m})(t_{0}^{s}t_{1}^{r})\right),
\end{eqnarray*}
where $\phi$ is the linear functional on 
$\C_{*}[t_{0}^{\pm 1},t_{1}^{\pm 1}]$ defined by
$$\phi(t_{0}^{n}t_{1}^{m})=\delta_{n,0}\;\;\;\mbox{ for }
m,n\in \Z.$$
This completes the proof.
\end{proof}

Now, we have an associative algebra 
$A\otimes \C_{*}[t_{0}^{\pm 1},t_{1}^{\pm 1}]$
equipped with a symmetric invariant bilinear form $\<\cdot,\cdot\>$.
Clearly, the bilinear form is also invariant
when $A\otimes \C_{*}[t_{0}^{\pm 1},t_{1}^{\pm 1}]$ is viewed as a Lie
algebra. For convenience, we set
\begin{eqnarray}\label{eA*}
A_{*}=A\otimes \C_{*}[t_{0}^{\pm 1},t_{1}^{\pm 1}].
\end{eqnarray}

Form the affine Lie algebra
\begin{eqnarray}
\widehat{A_{*}}
=(A\otimes \C_{*}[t_{0}^{\pm 1},t_{1}^{\pm 1}])\otimes
\C[t,t^{-1}]\oplus \C {\bf k}.
\end{eqnarray}
Denote by $V_{\widehat{A_{*}}}(\ell,0)$ the $\N$-graded vertex algebra
associated with the affine Lie algebra $\widehat{A_{*}}$ of level $\ell$.
Comparing (\ref{e7.68})--(\ref{e7.70}) with
(\ref{enonassoc-product}) and (\ref{e-nonassoc-form})
we see that $V_{W}$ is an $\hat{A}_{*}$-module of level $\ell$.
Using the same arguments as we used for Proposition \ref{ptwisted-module} 
(also for Proposition \ref{pquantum-Heisenberg})
we immediately have:

\bp{pextendedaffine}
Let $A$ be an associative algebra equipped with a symmetric invariant
bilinear form $\<\cdot,\cdot\>$ and let
$q$ be a nonzero complex number which is not a root of unity.
Let $\hat{A}_{q}$ be the quantum torus Lie algebra defined in (\ref{eEALA}).
Let $A_{*}$ be the Lie algebra defined in (\ref{eA*}) equipped with 
the symmetric invariant bilinear form.
Let $W$ be any restricted $\hat{A}_{q}$-module of level $\ell\in \C$,
on which the central element $c_{1}$ acts as zero. 
Then there exists a unique quasi
$V_{\widehat{A_{*}}}(\ell,0)$-module structure on $W$
with $Y_{W}(a\otimes t_{0}^{m}t_{1}^{n},x)=\bar{X}(a,m,n,x)$ 
for $a\in A,\; m,n\in \Z$.
\ep

\br{rproblem2}
{\em One can also study the case with $q$ a root of unity
and we leave this exercise to an interested reader.
If $q$ is a primitive $k$th root of unity, then
\begin{eqnarray}
\bar{X}(a,m,rk+n,x)=\bar{X}(a,m,n,x)
\end{eqnarray}
for $a\in A,\; m,n,r\in \Z$. 
In this case, the corresponding vertex algebra is smaller.}
\er

\end{document}